# A quantitative pairwise comparison-based constraint handling technique for constrained optimization

Ting Huang[1,2,3], Qiang Zhang[1,2,3*], Witold Pedrycz[4], Shanlin Yang[1,2,3]

[1]*School of Management, Hefei University of Technology, Hefei, Box 270, Hefei 230009, Anhui, P.R. China;*

[2]*Key Laboratory of Process Optimization and Intelligent Decision-making, Ministry of Education, Hefei, Box 270, Hefei 230009, Anhui, P.R. China;*

[3]*Ministry of Education Engineering Research Center for Intelligent Decision-Making & Information System Technologies, Hefei 230009, Anhui, P.R. China;*

[4]*Department of Electrical & Computer Engineering, University of Alberta, Edmonton, AB T6G 2R3, Canada*

**Abstract**

This study proposes a new constraint handling technique for assisting metaheuristic optimization algorithms to solve constrained optimization problems more effectively and efficiently. Given any two solutions of any constrained optimization problems, they are first mapped into a two-dimensional Cartesian coordinate system with their objective function value differences and constraint violation differences as the two axes. To the best of our knowledge, we are the first to deal with constraints by building such a Cartesian coordinate system. Then, the Cartesian coordinate system is divided into a series of grids by assigning ranks to different intervals of differences. In this way, a pairwise comparison criterion is derived with the use of the fused ranks, which achieves non-hierarchical comparison neither preferring objective function values nor constraint violations, resulting in more accurate evaluation compared with existing techniques. Moreover, an evaluation function that is equivalent to the pairwise comparison criterion is proposed, which further improves computational efficiency. The effectiveness and efficiency of the proposed constraint handling technique are verified on two well-known public datasets, that is, CEC 2006 and CEC 2017. The results demonstrate that metaheuristic optimization algorithms with using the proposed constraint handling technique can converge to a feasible optimal solution faster and more reliably. Experimental analysis on the parameters involved reveal guidance for their optimal settings.

**Summary of Contribution**

A quantitative pairwise comparison-based constraint handling technique is proposed based on a Cartesian mapping with objective function value difference and constraint violation difference as axes, which assists metaheuristic optimization algorithms to solve constrained optimization problems more effectively and efficiently. The proposed constraint handling technique is a general one that can be easily used by existing metaheuristics optimization algorithms owing to that its parameter setting has a solid theoretical basis. Besides, the Cartesian mapping provides a theoretical framework to balance objective function values and constraint violations in solving constrained optimization problems. Under the theoretical framework, we are expected to deduce numerous quantitative pairwise comparison-based constraint handling technique.

---

[*] Corresponding author. Tel: +86 0551 62901501; fax: +86 0551 62905263.

*E-mail address:* qiang_zhang@hfut.edu.cn. (Q. Zhang).

# 1. Introduction

Real-world optimization problems usually contain constraints that should be satisfied when an optimization algorithm is employed to obtain an optimal solution. Various metaheuristic optimization algorithms such as tabu search (Glover, 1990, 1989; Kulturel-Konak et al., 2004; Toth and Vigo, 2003), variable neighborhood search (Mladenović and Hansen, 1997; Pei et al., 2019), evolutionary algorithms (Holland, 1992; Nagata and Kobayashi, 2013; Storn and Price, 1997), swarm intelligence algorithms (Ang et al., 2020; Kennedy and Eberhart, 1995), simulated annealing algorithms (Burke et al., 2009; Steinbrunn et al., 1997), and others (Yang, 2010) have been proposed to solve these constrained optimization problems (COPs). However, the presence of constraints brings challenges to metaheuristic optimization algorithms since it may cause the feasible region to become very narrow, or separate. To this end, certain constraint handling techniques (CHTs) used to guide metaheuristic optimization algorithms close to the feasible region containing the optimal solution are indispensable. In this study, we focus on developing a CHT for metaheuristic optimization algorithms to solve COPs effectively and efficiently.

Overall, four types of CHTs have been developed, that is, penalty function-based methods (Coit et al., 1996; Di Pillo and Grippo, 1989), pairwise comparison-based methods (Deb, 2000; Runarsson and Yao, 2000; Takahama and Sakai, 2006), multi-objective optimization-based methods (Coello, 2000), and hybrid methods (Wu et al., 2022). A recent review on CHTs can be referred to (Rahimi et al., 2023). Penalty function-based methods convert constraints into a penalty function and integrate it into the objective function such that the considered COP is changed into an unconstrained optimization problem. Given multiple alternative solutions of a COP, pairwise comparison-based methods adopt certain comparison criteria to realize the superiority evaluation of these solutions. Multi-objective optimization-based methods convert a COP into a multi-objective optimization problem by considering constraints as one or multiple objective functions. Hybrid methods combine several different CHTs aiming to integrate their respective advantages. Among the four types of CHTs, the first two are actually evaluation-based methods that try to balance the contradiction between objective function values and constraint violations through certain evaluation function/criteria such that the alternatives approaching the feasible optimal solution are selected. The evaluation function is a quantitative form, while the evaluation criteria is a qualitative one.

Two important aspects, namely accuracy and computational cost of evaluation, are important for evaluation-based CHTs, where the former ensures that better alternatives can always be selected, and the latter ensures the computational efficiency of evaluation. Unfortunately, existing methods are difficult to achieve both aspects. Due to the use of the evaluation function, penalty function-based methods achieve efficient computation, but it is difficult to guarantee the evaluation accuracy because it is difficult to determine the introduced problem-dependent parameters, namely penalty factors. As for pairwise comparison-based methods, although they are able to accurately select alternatives, they are computationally inefficient because pairwise comparisons are required.

In existing pairwise comparison-based methods, various comparison criteria, including stochastic ranking (Runarsson and Yao, 2000), feasibility rules (Deb, 2000), and $\varepsilon$ level comparison (Takahama and Sakai, 2006), have been proposed. In general, two comparison indices, that is, objective function values and constraint violations, are involved in these comparison criteria. The random ranking uses either objective function values or constraint violations as the comparison index in each comparison, while the other two use both of them and thus are able to obtain better evaluation accuracy because of the complete exploitation of information. More specifically, the feasibility rules and $\varepsilon$ level comparison use a two-level hierarchical comparison to evaluate the alternatives, where constraint violations are always preferred. Although such a preference can make a metaheuristic optimization algorithm search for a feasible solution as much as possible, it may also make the obtained feasible solution not the optimal one.

According to the analysis completed above, if we can design a non-hierarchical comparison criterion neither preferring objective function values nor constraint violations, and use a quantitative evaluation function to implement such a qualitative evaluation criterion, we may obtain a CHT enjoying both high evaluation accuracy and high computational efficiency. To this end, in this study we propose a difference rank-based comparison criterion realizing a non-hierarchical comparison. Furthermore, a quantitative evaluation function is developed and proven to be equivalent to the proposed comparison criterion.

In the difference rank-based comparison criterion, given any two alternative solutions ($\mathbf{x}_1$, $\mathbf{x}_2$) are mapped into a two-dimensional Cartesian coordinate system with objective function value difference and constraint violation difference as its two axes. The two types of difference are further divided into a series of intervals to which assigned are integer difference ranks based on their corresponding difference values. In detail, the larger the difference values will be assigned a greater difference rank. In particular, zero will be assigned in the case when the difference value is equal to zero. When the assigned difference rank is greater than zero, it indicates that the objective function value or the constraint violation of $\mathbf{x}_1$ is greater than that of $\mathbf{x}_2$. It can be seen that the difference rank actually represents the relative performance of $\mathbf{x}_1$ and $\mathbf{x}_2$ in terms of objective functions or constraint violations. In this way, the Cartesian coordinate system is actually divided into small grids that are collections of ($\mathbf{x}_1$, $\mathbf{x}_2$) with the similar relative performance in terms of objective functions and constraint violations. By additively fusing the two difference ranks, we accordingly obtain the relative performance of $\mathbf{x}_1$ and $\mathbf{x}_2$ as a whole. The fusing result is considered as an important asset used for designing the following comparison criteria.

(1) When the fusing result is greater than zero, we definitely consider $\mathbf{x}_1$ to be worse than $\mathbf{x}_2$.
(2) When the fusing result is equal to zero, we consider the pros and cons of $\mathbf{x}_1$ and $\mathbf{x}_2$ cannot be judged. In this case, additional comparison criteria are devised with the help of a line defined in the Cartesian coordinate system. More specifically, (a) when ($\mathbf{x}_1$, $\mathbf{x}_2$) is mapped into the region below the line, we consider $\mathbf{x}_1$ to be better than $\mathbf{x}_2$; (b) when ($\mathbf{x}_1$, $\mathbf{x}_2$) is mapped into the region on the line, we consider $\mathbf{x}_1$ is equivalent to $\mathbf{x}_2$. (c) when ($\mathbf{x}_1$, $\mathbf{x}_2$) is mapped into the region above the line, we consider $\mathbf{x}_1$ is worse than $\mathbf{x}_2$.
(3) When the fusing result is less than zero, we definitely consider $\mathbf{x}_1$ to be better than $\mathbf{x}_2$.

The comparison criteria above are indeed non-hierarchical without preference for neither objective function values nor constrain violations, and thus it is possible to obtain a more accurate evaluation compared to the existing ones. Besides, we found an evaluation function $\pi(\mathbf{x})$ satisfying that $\pi(\mathbf{x}_1) - \pi(\mathbf{x}_2)$ meets the following conditions when given any two solutions $\mathbf{x}_1$ and $\mathbf{x}_2$. Moreover, we theoretically analyzed the parameters involved in $\pi(\mathbf{x})$ and obtain the mechanism for how the parameters affect the evaluation of $\pi(\mathbf{x})$, which has important theoretical implications for parameter settings.

(1) In the case (1) and the sub-case (c) of the case (2), $\pi(\mathbf{x}_1) - \pi(\mathbf{x}_2)$ is always less than zero.
(2) In the case (3) and the sub-case (a) of the case (2), $\pi(\mathbf{x}_1) - \pi(\mathbf{x}_2)$ is always greater than zero.

(3) In the sub-case (b) of the case (2), $\pi(\mathbf{x}_1) - \pi(\mathbf{x}_2)$ is always equal to zero.

Overall, the originality of this study lies in proposing a quantitative pairwise comparison-based CHT for assisting metaheuristic optimization algorithms in solving COP more effectively and efficiently in comparison to existing commonly used approaches. Our main contributions are summarized below.

(1) The proposed quantitative pairwise comparison-based CHT is a general one such that it can be applied to the solve process of COPs using almost any metaheuristic optimization algorithms.

(2) The proposed quantitative pairwise comparison-based CHT not only achieves accurate evaluation and high computational efficiency, but also exhibits a solid theoretical basis for the setting of parameters.

(3) The mapping based on objective function value difference and constraint violation difference provides a theoretical framework for us to balance objective function values and constraint violations in COPs. Under the theoretical framework, we are expected to deduce numerous quantitative pairwise comparison-based CHTs.

The structure of this paper is arranged as follows. In Section 2, we describe in detail the problem this study focuses on. In Section 3 we introduce four classical CHTs related to this study and summarize their advantages and disadvantages. Section 4 presents the proposed quantitative pairwise comparison-based CHT. Computational experiments are included in Section 5. We conclude this study in Section 6.

Throughout the paper we use lowercase italic letters for variables, lowercase bold letters for vectors, and uppercase italic letters for sets.

## 2. Problem description

Without loss of generality, a constrained optimization problem (COP) is defined as

$$\min f(\mathbf{x}) \\ s.t. \begin{cases} g_i(\mathbf{x}) \leq 0, i = 1,...,m \\ h_i(\mathbf{x}) = 0, i = m+1,...,n \end{cases} \quad , \tag{1}$$

where $\mathbf{x} = (x_1, x_2, ..., x_k)$ is the decision vector with $k$ decision variables, $g_i(\mathbf{x}) \leq 0$ denotes the inequality constraint, $h_i(\mathbf{x}) = 0$ denotes the equality constraint, $m$ is the number of inequality constraints, and $(n - m)$ is the number of equality constraints. Each decision variable $x_j$ ($j = 1, 2, ..., k$) is bounded by a lower limit $l_j$ and an upper limit $u_j$, namely $l_j \leq x_j \leq u_j$. Let $S = \{\mathbf{x} \mid l_j \leq x_j \leq u_j, j = 1, 2, ..., k\}$ denote the search region and $F = \{\mathbf{x} \mid \mathbf{x} \in S \wedge g_i(\mathbf{x}) \leq 0 \wedge h_i(\mathbf{x}) = 0, i = 1, 2, ..., n\}$ the feasible region.

To handle equality constraints, they are generally transformed into inequality ones by introducing a very small tolerance value $\delta > 0$. Accordingly, Equation (1) is converted to Equation (2),

$$\min f(\mathbf{x}) \\ s.t. \begin{cases} g_i(\mathbf{x}) \leq 0, i = 1,...,m \\ |h_i(\mathbf{x})| - \delta \leq 0, i = m+1,...,n \end{cases} \quad . \tag{2}$$

In Equation (2), the degree of constraint violation of $\mathbf{x}$ on the $i$th constraint is formulated as

$$v_i(\mathbf{x}) = \begin{cases} \max\{g_i(\mathbf{x}), 0\}, i = 1, ..., m \\ \max\{|h_i(\mathbf{x})| - \delta, 0\}, i = m+1, ..., n \end{cases}. \tag{3}$$

Subsequently, the degree of constraint violation of **x** on all constraints is calculated as the sum of $v_i(\mathbf{x})$ ($i$ = 1, 2, …, $n$), that is,

$$\vartheta(\mathbf{x}) = \sum_{i=1}^{n} v_i(\mathbf{x}). \tag{4}$$

If $\vartheta(\mathbf{x}) = 0$, **x** is called a feasible solution, otherwise we regard it as infeasible one.

All metaheuristic optimization algorithms for solving aforesaid COPs try to search for a feasible optimal solution $\mathbf{x} \in S$ whose $f(\mathbf{x})$ is as small as possible while keeping $\vartheta(\mathbf{x}) = 0$. The search process is generally carried out iteratively. After each iteration, new solutions are generated based on the information provided by old ones. Subsequently, CHTs provide a guide for metaheuristic optimization algorithms to select better solutions between the newly-generated ones and the old ones so that they eventually converge to the feasible optimal solution. Therefore, CHTs are important and indispensable in solving COPs.

It is obvious that evaluation criteria defined by CHTs for selecting better solutions are very important for optimization algorithms. If only $f(\mathbf{x})$ is considered, metaheuristic optimization algorithms may not be able to search for a feasible solution. If only $\vartheta(\mathbf{x})$ is considered, metaheuristic optimization algorithms may not be able to search for an optimal solution. Consequently, a reasonably efficient balance between $f(\mathbf{x})$ and $\vartheta(\mathbf{x})$ is the key to a CHT selecting a better solution, and it is also the basic problem this paper focus on.

**3. Classical constraint handling techniques**

Virous CHTs have been proposed in previous studies, such as penalty functions, feasibility rules, stochastic ranking, and the $\varepsilon$ level comparison. This section introduces these classical CHTs and summarize their advantages and disadvantages.

3.1. Penalty functions

Given a solution **x** of a COP, a quantitative function denoted by $\tilde{f}(\mathbf{x})$ for the evaluation of the pros and cons of **x** is defined as an additive combination of $f(\mathbf{x})$ and $p(\mathbf{x})$, that is,

$$\tilde{f}(\mathbf{x}) = f(\mathbf{x}) + p(\mathbf{x}), \tag{5}$$

where

$$p(\mathbf{x}) = \sum_{i=1}^{n} r_i \cdot v_i(\mathbf{x}) \tag{6}$$

refers to a penalty function, and $r_i \geq 0$ is user-defined parameters called penalty coefficients.

It can be seen that $p(\mathbf{x})$ is a parameterized function of $\vartheta(\mathbf{x})$, and parameters $r_i$ achieve a balance between $f(\mathbf{x})$ and $\vartheta(\mathbf{x})$. The larger parameters $r_i$ are set, the greater the effect of $\vartheta(\mathbf{x})$ relative to $f(\mathbf{x})$, and thus appropriate value of $r_i$ (not too big nor too small) is the key to using penalty functions. Unfortunately, it is difficult to determine $r_i$ in practical applications because it is usually problem-dependent.

3.2. Feasibility rules

Given two arbitrary solutions $\mathbf{x}_1$ and $\mathbf{x}_2$ of a COP, the pros and cons between $\mathbf{x}_1$ and $\mathbf{x}_2$ is evaluated by the following three pairwise comparison-based rules.

(1) If $\vartheta(\mathbf{x}_1) \leq 0$, $\vartheta(\mathbf{x}_2) \leq 0$, and $f(\mathbf{x}_1) < f(\mathbf{x}_2)$, $\mathbf{x}_1$ is preferred to $\mathbf{x}_2$;

(2) If $\vartheta(\mathbf{x}_1) \leq 0$ and $\vartheta(\mathbf{x}_2) > 0$, $\mathbf{x}_1$ is preferred to $\mathbf{x}_2$;

(3) If $\vartheta(\mathbf{x}_1) > 0$, $\vartheta(\mathbf{x}_2) > 0$, and $\vartheta(\mathbf{x}_1) < \vartheta(\mathbf{x}_2)$, $\mathbf{x}_1$ is preferred to $\mathbf{x}_2$.

It can be seen that the three rules are defined on the three cases (a)-(c) as shown in Fig. 1, and each rule achieves a balance between $f(\mathbf{x})$ and $\vartheta(\mathbf{x})$ by a $\vartheta(\mathbf{x})$-prioritized two-level pairwise comparison. However, when solutions are scattered near the boundary of the feasible and infeasible regions as shown in the case (d), the rule (2) may produce incorrect evaluation since those slightly infeasible solutions that are close to the feasible region and have better $f(\mathbf{x})$ (e.g., $\mathbf{x}_2$) are arbitrarily judged to be worse than any feasible solution (e.g., $\mathbf{x}_1$).

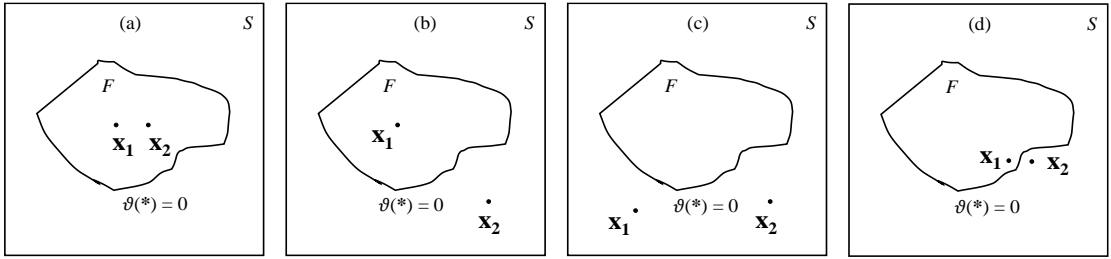

Fig. 1. Four distributions of the two compared solutions by feasibility rules: (a) both solutions are in the feasible region; (b) one is inside the feasible region and the other is outside the feasible region; (c) both solutions are outside the feasible region; (d) the two solutions are located on either side of the feasible region boundary.

3.3. The $\varepsilon$ level comparison

Given two arbitrary solutions $\mathbf{x}_1$ and $\mathbf{x}_2$ of a COP, the pros and cons between $\mathbf{x}_1$ and $\mathbf{x}_2$ is evaluated by the following four pairwise comparison-based rules.

(1) If $\vartheta(\mathbf{x}_1) \leq \varepsilon$, $\vartheta(\mathbf{x}_2) \leq \varepsilon$, and $f(\mathbf{x}_1) < f(\mathbf{x}_2)$, $\mathbf{x}_1$ is preferred to $\mathbf{x}_2$;

(2) If $\vartheta(\mathbf{x}_1) \leq \varepsilon$, $\vartheta(\mathbf{x}_2) > \varepsilon$, $\mathbf{x}_1$ is preferred to $\mathbf{x}_2$;

(3) If $\vartheta(\mathbf{x}_1) > \varepsilon$, $\vartheta(\mathbf{x}_2) > \varepsilon$, and $\vartheta(\mathbf{x}_1) < \vartheta(\mathbf{x}_2)$, $\mathbf{x}_1$ is preferred to $\mathbf{x}_2$.

(4) If $\vartheta(\mathbf{x}_1) = \vartheta(\mathbf{x}_2)$, and $f(\mathbf{x}_1) < f(\mathbf{x}_2)$, $\mathbf{x}_1$ is preferred to $\mathbf{x}_2$.

where $\varepsilon \in [0, \infty]$ denotes the threshold regarding the satisfaction of $\vartheta(\mathbf{x})$. A smaller $\varepsilon$ indicates that a CHT prioritizes $\vartheta(\mathbf{x})$ when performing pairwise comparisons of solutions. In the case of $\varepsilon = 0$, the $\varepsilon$ level comparison is equivalent to the classical feasibility rules in which $\vartheta(\mathbf{x})$ precedes $f(\mathbf{x})$. When $\varepsilon = \infty$, the $\varepsilon$ level comparison depends entirely on the comparison of $f(\mathbf{x})$.

The four rules are defined on the four cases (a)-(d) as shown in Fig. 2, and each rule achieves a balance between $f(\mathbf{x})$ and $\vartheta(\mathbf{x})$ by a $\vartheta(\mathbf{x})$-prioritized two-level pairwise comparison. We can see that those slightly infeasible solutions are considered as feasible ones by the $\varepsilon$ level comparison, which avoids incorrect evaluations to some extent.

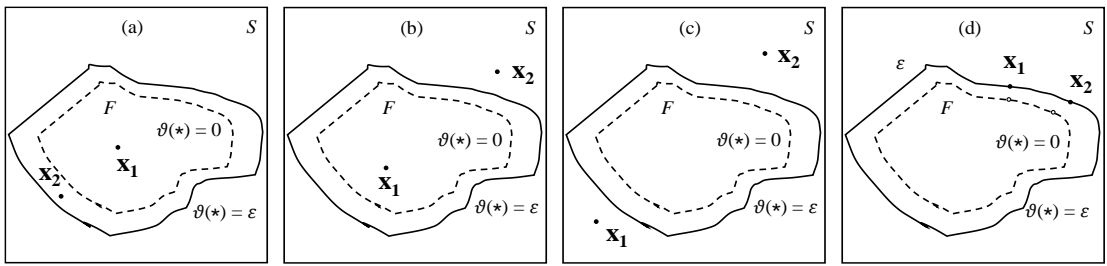

Fig. 2. Four distributions of the two compared solutions by the $\varepsilon$ level comparison: (a) one is inside the constrain violation line $\vartheta(*) = 0$ and the other is between the constrain violation lines

$\vartheta(*) = 0$ and $\vartheta(*) = \varepsilon$; (b) one is inside the constrain violation line $\vartheta(*) = 0$ and the other is outside the constrain violation line $\vartheta(*) = \varepsilon$; (c) both solutions are outside the constrain violation line $\vartheta(*) = \varepsilon$; (d) both solutions are on a certain constrain violation line.

3.4. Stochastic ranking

Given a solution sequence $X = <\mathbf{x}_1, \mathbf{x}_2, …, \mathbf{x}_l>$ of a COP, $X$ is traversed $\bar{t}$ times by a bubble-sort-like procedure to produces a new solution sequence $Y = <\mathbf{y}_1, \mathbf{y}_2, …, \mathbf{y}_l>$. In each traverse, two adjacent solutions $\mathbf{x}_s$ and $\mathbf{x}_{s+1}$ of $X$ are compared either by $f(\mathbf{x})$ or by $\vartheta(\mathbf{x})$, and the probability of the former is defined as $0 \leq p \leq 1$. The algorithm of the stochastic ranking is specified as Algorithm 1.

---

**Algorithm 1: Stochastic ranking**

---

**Input**: $X = (\mathbf{x}_1, \mathbf{x}_2, …, \mathbf{x}_l)$, $\bar{t}$, $p$
**Output**: $Y = (\mathbf{y}_1, \mathbf{y}_2, …, \mathbf{y}_l)$
1: $Y \leftarrow X$
2: **for** $t = 1$ **to** $\bar{t}$ **do**
3:    **for** $s = 1$ **to** $l - 1$ **do**
4:       $\bar{p} \leftarrow$ random()
5:       **if** $\bar{p} < p$ or $\vartheta(\mathbf{x}_s) = \vartheta(\mathbf{x}_{s+1}) = 0$ **then**
6:          **if** $f(\mathbf{x}_s) > f(\mathbf{x}_{s+1})$ **then**
7:             $\mathbf{y}_s \leftarrow \mathbf{x}_{s+1}$
8:             $\mathbf{y}_{s+1} \leftarrow \mathbf{x}_s$
9:          **end if**
10:      **else**
11:         **if** $\vartheta(\mathbf{x}_s) > \vartheta(\mathbf{x}_{s+1})$ **then**
12:            $\mathbf{y}_s \leftarrow \mathbf{x}_{s+1}$
13:            $\mathbf{y}_{s+1} \leftarrow \mathbf{x}_s$
14:         **end if**
15:      **end if**
16:    **end for**
17:    $X \leftarrow Y$
18: **end for**

---

As shown in Algorithm 1, a balance between $f(\mathbf{x})$ and $\vartheta(\mathbf{x})$ is achieved by the definition of $p$, and the priority of $f(\mathbf{x})$ and $\vartheta(\mathbf{x})$ is determined by the value of $p$. When $0 \leq p < 0.5$, $\vartheta(\mathbf{x})$ takes precedence over $f(\mathbf{x})$. When $p = 0.5$, $f(\mathbf{x})$ and $\vartheta(\mathbf{x})$ have equal priority. When $0.5 < p \leq 1$, $f(\mathbf{x})$ takes precedence over $\vartheta(\mathbf{x})$.

3.5. Performance analysis

The four CHTs above can actually be further divided into two types according to their evaluation formalisms, that is, quantitative CHTs (penalty functions) and qualitative pairwise comparison-based CHTs (feasibility rules, the $\varepsilon$ level comparison, and stochastic ranking). In what followings, we analyze the performance of the four CHTs above in terms of two aspects, that is, evaluation accuracy and computational efficiency.

Regarding evaluation accuracy, the following results can be obtained.

(1) The evaluation accuracy of penalty functions completely depends on the setting of parameters $r_i$. However, it is difficult to determine parameters $r_i$ such that the evaluation accuracy of penalty functions can't be guaranteed.

(2) The evaluation of feasibility rules is basically accurate but produces inaccurate results when

those slightly infeasible solutions are involved.

(3) The evaluation of the $\varepsilon$ level comparison is basically accurate and depends on the setting of the parameter $\varepsilon$.

(4) The evaluation accuracy of stochastic ranking is random and depends on the setting of the parameter $p$.

On the one hand, due to the use of pairwise comparison strategy, where $f(\mathbf{x})$ and $\vartheta(\mathbf{x})$ are used as two separate comparison indices, evaluation accuracy of qualitative pairwise comparison-based CHTs is generally higher than that of quantitative ones. On the other hand, since the comparison criterion of qualitative pairwise comparison-based CHTs is hierarchical, that is, $\vartheta(\mathbf{x})$ is always preferred, inaccuracies may arise in some special cases. Besides, the evaluation accuracy of existing CHTs depends on the parameter setting, except feasibility rules because it has no parameters. As such, certain dynamic or adaptive parameter setting strategies are often used in the iterative process of metaheuristic optimization algorithms to improve the evaluation accuracy. Nevertheless, these strategies include some pieces of experiences gained from experiments, however they lack rigorous theoretical analysis.

As for computational efficiency, the results shown in Table 1 can be obtained. Penalty functions realize a linear time complexity for evaluating $l$ solutions of a COP with the help of a quantitative function, while those of the other three CHTs are quadratic time complexity due to the adoption of qualitative pairwise comparison.

Table 1. Time complexity of different CHTs

| Penalty functions | Feasibility rules | The $\varepsilon$ level comparison | Stochastic ranking |
| --- | --- | --- | --- |
| $O(l)$ | $O((l \times (l-1)/2))$ | $O((l \times (l-1)/2))$ | $O(\bar{t} \times l)$ |

*Note*: The computational complexity of $\vartheta(\mathbf{x})$ is ignored because all contrasting CHTs involve the computation of $\vartheta(\mathbf{x})$.

Overall, quantitative CHTs contribute to higher computational efficiency compared to qualitative ones, and qualitative CHTs contribute to higher evaluation accuracy compared to quantitative ones. Therefore, this study tries to improve existing qualitative CHTs and quantify the improved one to obtain both high computational efficiency and evaluation accuracy.

**4. Proposed quantitative pairwise comparison-based constraint handling technique**

This section introduces the proposed CHT. The two core components of the proposed CHT, that is, the difference rank-based comparison criterion and the corresponding quantization function, are described in detail.

4.1. Difference rank-based comparison criterion

In this subsection, we first define the basic concepts of difference division and rank. After that, the determination of a function involved in the definition of constrain violation difference is discussed. Finally, the comparison criterion based on difference rank is presented.

4.1.1. Difference division and rank

Given a variable $\varDelta \in [-b, b]$ denotes the degree of difference in a certain aspect between two compared subjects, we define an ordered set $Q_\alpha = <q_1, q_2, \ldots, q_\alpha>$ for a division of $\varDelta$, where

$$0 < q_1 < q_2 < \cdots < q_\alpha < b , \tag{7}$$

and $\alpha \in N$. Specifically, $\alpha = 0$ indicates $Q_\alpha = \emptyset$. Accordingly, a division of $\varDelta$ with the use of $Q_\alpha$ is

defined as an order set $D_\alpha$ below

$$D_\alpha = \begin{cases} \langle [-b,-q_\alpha),[-q_\alpha,-q_{\alpha-1}),\cdots,[-q_1,0),0,(0,q_1],\cdots,(q_{\alpha-1},q_\alpha],(q_\alpha,b] \rangle, & \text{if } \alpha \in Z^+ \\ \langle [-b,0),0,(0,b] \rangle, & \text{if } \alpha = 0 \end{cases}. \quad (8)$$

Each element in $D_\alpha$ is assigned a difference rank such that an ordered set $R_\alpha$ corresponding to $D_\alpha$ is obtained as follows,

$$R_\alpha = \begin{cases} \langle -(\alpha+1),-\alpha,\cdots,-1,0,1,\cdots,\alpha,\alpha+1 \rangle, & \text{if } \alpha \in Z^+ \\ \langle -1,0,1 \rangle, & \text{if } \alpha = 0 \end{cases}. \quad (9)$$

For notation convenience, an element in $D_\alpha$ with a difference rank $\chi \in R_\alpha$ is denoted as $d_\chi$. For example, $d_0 = 0$, $d_1 = (0, p_1]$, and so on.

Given any two solutions $\mathbf{x}_1$ and $\mathbf{x}_2$ of a COP, let

$$y = f(\mathbf{x}_1) - f(\mathbf{x}_2) \quad (10)$$

denote the degree of difference between $\mathbf{x}_1$ and $\mathbf{x}_2$ in terms of objective function values. If $y < 0$, $\mathbf{x}_1$ is better than $\mathbf{x}_2$ in terms of objective function values. Suppose that

$$\overline{f} = \max_{\mathbf{x} \in S}\{f(\mathbf{x})\}, \quad \underline{f} = \min_{\mathbf{x} \in S}\{f(\mathbf{x})\}, \quad (11)$$

then we can have that

$$y \in \left[ -(\overline{f} - \underline{f}), \overline{f} - \underline{f} \right]. \quad (12)$$

Besides, let

$$z = \theta(\mathbf{x}_1) - \theta(\mathbf{x}_2) \quad (13)$$

denote degree of difference between $\mathbf{x}_1$ and $\mathbf{x}_2$ in terms of constraint violations, where $\theta$ denote a function that characterizes constraint violations. The determination of $\theta$ will be discussed in the following sub-section in detail. If $z < 0$, $\mathbf{x}_1$ is better than $\mathbf{x}_2$ in terms of constraint violations. Suppose that

$$\overline{\theta} = \max_{\mathbf{x} \in S}\{\theta(\mathbf{x})\}, \quad \underline{\theta} = \min_{\mathbf{x} \in S}\{\theta(\mathbf{x})\}, \quad (14)$$

then we can have that

$$z \in \left[ -(\overline{\theta} - \underline{\theta}), \overline{\theta} - \underline{\theta} \right]. \quad (15)$$

With the definitions of $y$ and $z$, let $D_\alpha^y$ and $D_\beta^z$ denote two divisions of $y$ and $z$, respectively, where $\alpha$ and $\beta$ are two user-defined parameters. Furthermore, let $R_\alpha^y$ and $R_\beta^z$ denote two sets of difference ranks corresponding to $D_\alpha^y$ and $D_\beta^z$, where their elements are denoted as $\chi$ and $\lambda$, that is, $\chi \in R_\alpha^y$ and $\lambda \in R_\beta^z$. In this way, given any two solutions $\mathbf{x}_1$ and $\mathbf{x}_2$ of a COP, we can calculate two difference ranks, namely $\chi$ and $\lambda$, for the comparison in terms of objective function values and constraint violations, respectively. To provide a comprehensive comparison between $\mathbf{x}_1$ and $\mathbf{x}_2$, the Cartesian product ($\times$) of $y$ and $z$ are further defined as

$$s = y \times z, \tag{16}$$

and the range of $s$ can be easily derived as

$$s \in G = \left[-(\overline{f}-\underline{f}), \overline{f}-\underline{f}\right] \times \left[-(\overline{\theta}-\underline{\theta}), \overline{\theta}-\underline{\theta}\right]. \tag{17}$$

Based on the operation of Equation (16), we actually map the pair ($\mathbf{x}_1$, $\mathbf{x}_2$) into a two-dimensional Cartesian coordinate system with $y$ and $z$ as its two axes, where $G$ is the set of all mapped points $s$.

We can regard $s$ as a variable fusing two types of difference, namely the ones regarding objective function values and constraint violations, and the division of $s$ is defined as follows,

$$D_\gamma^s = \left\{ d_\varphi^s = \bigcup_{\chi+\lambda=\varphi} d_\chi^y \times d_\lambda^z \middle| d_\chi^y \in D_\alpha^y, d_\lambda^z \in D_\beta^z \right\}, \tag{18}$$

where

$$\varphi \in R_\gamma^s = \left\{\chi + \lambda \middle| \chi \in R_\alpha^y, \lambda \in R_\beta^z\right\} = \left\{-(\alpha+\beta+2), \cdots, (\alpha+\beta+2)\right\} \tag{19}$$

denotes the composite difference rank that fuses $\chi$ and $\lambda$.

An illustrative figure (Fig. 3) is presented for explaining the aforesaid mapping and division in detail. Different settings of the two parameters, $\alpha$ and $\beta$, lead to different mapping results, which is further illustrated in Fig. 4.

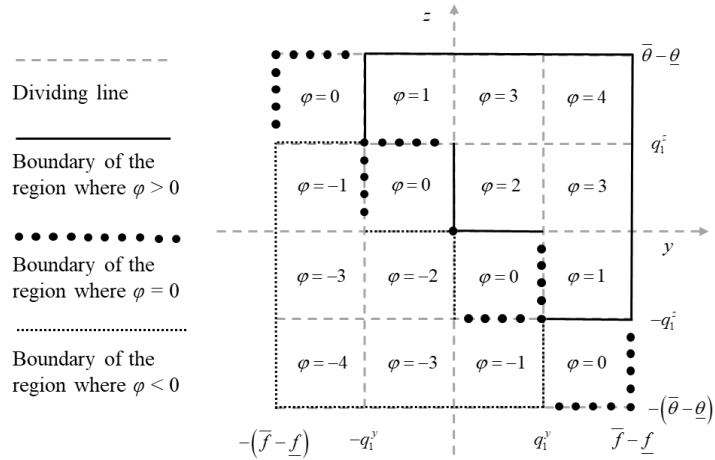

Fig. 3. Mapping and division results of pair ($\mathbf{x}_1$, $\mathbf{x}_2$)

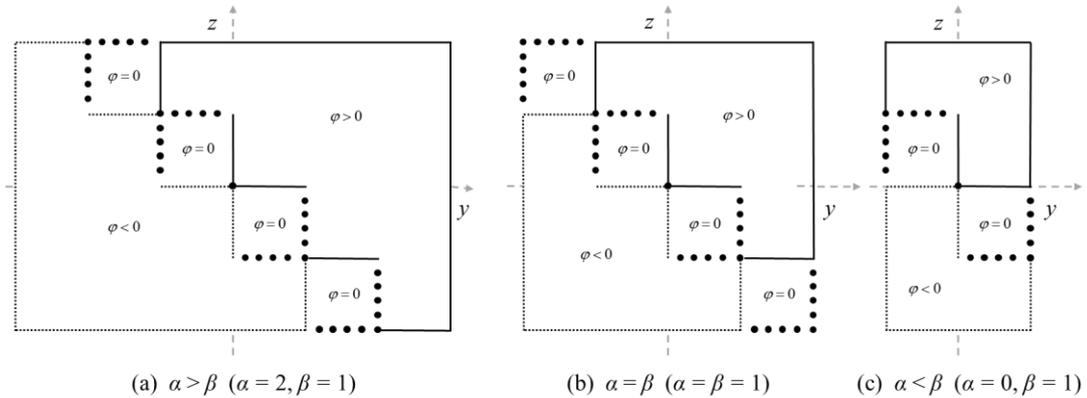

(a) $\alpha > \beta$ ($\alpha = 2, \beta = 1$)  (b) $\alpha = \beta$ ($\alpha = \beta = 1$)  (c) $\alpha < \beta$ ($\alpha = 0, \beta = 1$)

Fig. 4. Mapping and division results under different settings of $\alpha$ and $\beta$

4.1.2. Determination of $\theta$

Formally, the analysis above can be described as a mapping $s: S \times S \mapsto G$. Given $F \subset S$ is the feasible region, we can divide $S \times S$ into four complementarily subsets as follows.

$$T_1 = F \times F, \quad T_2 = F \times (S - F), \quad T_3 = (S - F) \times F, \quad T_4 = (S - F) \times (S - F), \tag{20}$$

where $T_1$ denotes the set of the pair $(\mathbf{x}_1, \mathbf{x}_2)$ where $\mathbf{x}_1$ and $\mathbf{x}_2$ are both feasible solutions, $T_2$ denotes the set of the pair $(\mathbf{x}_1, \mathbf{x}_2)$ where $\mathbf{x}_1$ is a feasible solution while $\mathbf{x}_2$ is a feasible solution, and so on. Let

$$G_1 = s(T_1), \quad G_2 = s(T_2), \quad G_3 = s(T_3), \quad G_4 = s(T_4), \tag{21}$$

it is necessary for us to find a specific mapping $s$ so that the following two properties hold,

**Property 1:** $G_1 \cup G_2 \cup G_3 \cup G_4 = G$;

**Property 2:** $G_1 \cap G_2 \cap G_3 \cap G_4 = \varnothing$.

The satisfaction of Property 1 and Property 2 ensures that the mapping $s$ is able to distinguish feasible solutions from infeasible ones, which is very important for tche design of subsequent comparison criteria.

Let

$$\theta(\mathbf{x}) = \begin{cases} -\overline{\vartheta}, & \text{if } \vartheta(\mathbf{x}) = 0 \\ \vartheta(\mathbf{x}), & \text{if } \vartheta(\mathbf{x}) \neq 0 \end{cases} \tag{22}$$

with

$$\overline{\vartheta} = \max_{\mathbf{x} \in S} \{\vartheta(\mathbf{x})\}, \tag{23}$$

then we can have that

$$z \in \left[-2\overline{\vartheta}, -\overline{\vartheta}\right) \cup \left(-\overline{\vartheta}, \overline{\vartheta}\right) \cup \left(\overline{\vartheta}, 2\overline{\vartheta}\right]. \tag{24}$$

**Lemma 1**: Given $\mathbf{x}_1 \in S$ and $\mathbf{x}_2 \in S$, $y = f(\mathbf{x}_1) - f(\mathbf{x}_2) \in \left[-(\overline{f} - \underline{f}), \overline{f} - \underline{f}\right]$, $z = \theta(\mathbf{x}_1) - \theta(\mathbf{x}_2) \in \left[-2\overline{\vartheta}, -\overline{\vartheta}\right) \cup \left(-\overline{\vartheta}, \overline{\vartheta}\right) \cup \left(\overline{\vartheta}, 2\overline{\vartheta}\right]$, $G = \left[-(\overline{f} - \underline{f}), \overline{f} - \underline{f}\right] \times \left[-2\overline{\vartheta}, -\overline{\vartheta}\right) \cup \left(-\overline{\vartheta}, \overline{\vartheta}\right) \cup \left(\overline{\vartheta}, 2\overline{\vartheta}\right]$, where $f$ denotes the objective function and $\theta$ is specified by Equation (22), the mapping $s: S \times S \mapsto G$ with $s = y \times z$ satisfies **Property 1** and **Property 2**.

In order to save space, the proofs of Lemma 1 and all subsequent lemmas are presented in Section II of Supplementary Material. The above analysis is illustrated as Fig. 5.

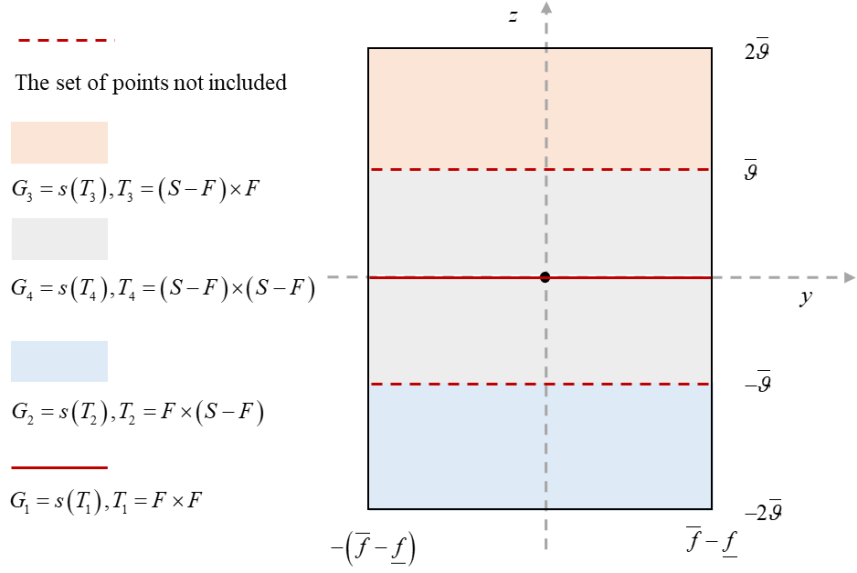

Fig. 5. Mapping results of pair ($\mathbf{x}_1$, $\mathbf{x}_2$) with the use of $\theta(\mathbf{x})$ specified by Equation (22).

4.1.3. Comparison criterion

As mentioned before, we first map any two solutions $\mathbf{x}_1$ and $\mathbf{x}_2$ of a COP into a two-dimensional Cartesian coordinate system with $y = f(\mathbf{x}_1) - f(\mathbf{x}_2)$ (difference of objective function values) and $z = \theta(\mathbf{x}_1) - \theta(\mathbf{x}_2)$ (difference of constraint violations) as its two axes, which is denoted as a mapping $s$: $S \times S \mapsto G$ with $s = y \times z$. Then, $G$ is divided into a series of grids by two divisions $D_\alpha^y$ and $D_\beta^z$ along $y$ and $z$, and each element in $D_\alpha^y$ and $D_\beta^z$ are assigned a difference rank $\chi$ and $\lambda$, respectively. Finally, the composite difference rank $\varphi$ corresponding to a grid is obtained by additively fusing $\chi$ and $\lambda$, that is, $\varphi = \chi + \lambda$. The obtained $\varphi$ actually characterizes the relative performance of the mapped $\mathbf{x}_1$ and $\mathbf{x}_2$ in the corresponding grid, and thus it is considered as an important asset used to design the following comparison criterion.

(1) When $\varphi > 0$: The following two sub-cases are involved: a) Both of $\chi$ and $\lambda$ are greater than zero, which means that $\mathbf{x}_1$ is worse than $\mathbf{x}_2$ in both objective function values and constraint violations; b) Only one of $\chi$ and $\lambda$ is greater than zero while the other is less than or equal to zero, and the absolute value of the former is greater than the latter, which means that $\mathbf{x}_1$ is worse than $\mathbf{x}_2$ in one aspect (either objective function values or constraint violations) and better than or equivalent to $\mathbf{x}_2$ in another aspect, moreover the degree of the former is greater than that of the latter. In both of the two sub-cases we definitely consider $\mathbf{x}_1$ to be worse than $\mathbf{x}_2$.

(2) When $\varphi = 0$: The following two sub-cases are involved: a) Both of $\chi$ and $\lambda$ are equal to zero, which means that $\mathbf{x}_1$ is actually the same as $\mathbf{x}_2$; b) One of $\chi$ and $\lambda$ is greater than zero while the other is less than zero, and their absolute values are the same, which means that $\mathbf{x}_1$ is worse than $\mathbf{x}_2$ in one aspect and better than $\mathbf{x}_2$ in another aspect, moreover the degree of the former is equal to that of the latter. In both of the two sub-cases, we consider the pros and cons of $\mathbf{x}_1$ and $\mathbf{x}_2$ cannot be judged by $\varphi$, and the further comparison criterion is discussed in detail below.

(3) When $\varphi < 0$: Similar to the analysis of case (1), we definitely consider $\mathbf{x}_1$ to be better than $\mathbf{x}_2$.

In the case of $\varphi = 0$, a line $L$ is first defined as follows,

$$L: z = -e \cdot y, y \in \left[-\left(\overline{f} - \underline{f}\right), \overline{f} - \underline{f}\right] \text{ and } z \in \left[-2\overline{\vartheta}, -\overline{\vartheta}\right) \cup \left(-\overline{\vartheta}, \overline{\vartheta}\right) \cup \left(\overline{\vartheta}, 2\overline{\vartheta}\right] \quad (25)$$

where

$$e = \begin{cases} \dfrac{2\overline{\vartheta}}{q_\alpha^y}, & \alpha > \beta \\ \dfrac{2\overline{\vartheta}}{\overline{f} - \underline{f}}, & \alpha = \beta \\ \dfrac{q_\beta^z}{\overline{f} - \underline{f}}, & \alpha < \beta \end{cases} \quad (26)$$

Fig. 6 depicts the line $L$ under the three cases $\alpha = \beta$, $\alpha > \beta$, and $\alpha < \beta$. In what follows, let's analyze the geometric meaning of the line $L$.

(1) In the case of $\alpha = \beta$, the line $L$ passes through three special points, namely $\left(-\left(\overline{f} - \underline{f}\right), 2\overline{\vartheta}\right)$, $(0, 0)$, and $\left(\left(\overline{f} - \underline{f}\right), -2\overline{\vartheta}\right)$. The point $\left(-\left(\overline{f} - \underline{f}\right), 2\overline{\vartheta}\right)$ indicates that $f(\mathbf{x}_1) = \underline{f}$, $f(\mathbf{x}_2) = \overline{f}$, $\theta(\mathbf{x}_1) = \overline{\vartheta}$, and $\theta(\mathbf{x}_2) = -\overline{\vartheta}$, which means that $\mathbf{x}_1$ is maximally better than $\mathbf{x}_2$ in terms of objective function values, and maximally worse than $\mathbf{x}_2$ in terms of constraint violations. The opposite result of the point $\left(-\left(\overline{f} - \underline{f}\right), 2\overline{\vartheta}\right)$ can be deduced at the point $\left(\left(\overline{f} - \underline{f}\right), -2\overline{\vartheta}\right)$. The point $(0, 0)$ indicates that $f(\mathbf{x}_1) = f(\mathbf{x}_2)$ and $\theta(\mathbf{x}_1) = \theta(\mathbf{x}_2)$, which means that $\mathbf{x}_1$ and $\mathbf{x}_2$ are the same solution. Therefore, at all the three points we should actually consider $\mathbf{x}_1$ and $\mathbf{x}_2$ to be equivalent. Similar conclusions can be drawn in other points scattered in the line $L$, and thus the line $L$ actually represents the set of all mapping points where $\mathbf{x}_1$ and $\mathbf{x}_2$ should be considered equivalent.

(2) In the case of $\alpha > \beta$, it can be seen that the rest of regions is the same as that of the case of $\alpha = \beta$ if we exclude the regions where $y < -q_\alpha^y$ and $y > q_\alpha^y$. In fact, the parameter $q_\alpha^y$ can be seen as the maximum tolerance of difference in terms of objective function values ($y$), which means that if the absolute value of $y$ between $\mathbf{x}_1$ and $\mathbf{x}_2$ exceeds $q_\alpha^y$, that is, $y < -q_\alpha^y$ or $y > q_\alpha^y$, then the performance of $\mathbf{x}_1$ over $\mathbf{x}_2$ is determined entirely by their objective function values and has nothing to do with their constraint violations. When $-q_\alpha^y < y < q_\alpha^y$, the determination of the performance of $\mathbf{x}_1$ over $\mathbf{x}_2$ depends on both objective function values and constraint violations, and the geometric meaning of the line $L$ is the same as that in the case of $\alpha = \beta$. From the above analysis it can be seen that the case of $\alpha = \beta$ is a special case of $\alpha > \beta$ where the maximum tolerance of difference in terms of objective function values is set to $\overline{f} - \underline{f}$.

(3) In the case of $\alpha < \beta$, a similar analysis result to that of the case of $\alpha > \beta$ can be obtained if we view the parameter $q_\beta^z$ as the maximum tolerance of difference in terms of constraint violations.

In summary, the line $L$ represents the set of all mapping points where $\mathbf{x}_1$ and $\mathbf{x}_2$ should be considered equivalent, which is an important basis for developing the comparison criteria in the region where $\varphi = 0$. As shown in Fig. 6, the region where $\varphi = 0$ is divided into three parts, namely the blue one, the red one, and the yellow one, by the line $L$. The following comparison criteria are further developed to judge the pros and cons of $\mathbf{x}_1$ and $\mathbf{x}_2$.

(1) If the pair of ($\mathbf{x}_1$, $\mathbf{x}_2$) is mapped into the line $L$, namely the red region, then it should consider that $\mathbf{x}_1$ is equivalent to $\mathbf{x}_2$ based on the aforesaid analysis.

(2) If the pair of ($\mathbf{x}_1$, $\mathbf{x}_2$) is mapped into the region below the line $L$, namely the blue region, then we have

$$z < -e \cdot y, z \cdot y < 0 \Leftrightarrow \begin{cases} \dfrac{|z|}{|y|} > e, y > 0, z < 0 \\ \dfrac{|y|}{|z|} > e, y < 0, z > 0 \end{cases}. \tag{27}$$

Equation (27) indicates that when $\mathbf{x}_1$ is better than $\mathbf{x}_2$ in one aspect but worse in another aspect, the ratio of superiority to inferiority is greater than the slope of the line $L$, namely $e$, which indicates that we should regard $\mathbf{x}_1$ to be better than $\mathbf{x}_2$ when both the two aspects (objective function values and constraint violations) are considered.

(3) If the pair of ($\mathbf{x}_1$, $\mathbf{x}_2$) is mapped into the region above the line $L$, namely the yellow region, we should consider $\mathbf{x}_1$ to be worse than $\mathbf{x}_2$ based on an analysis similar to that of Case (2).

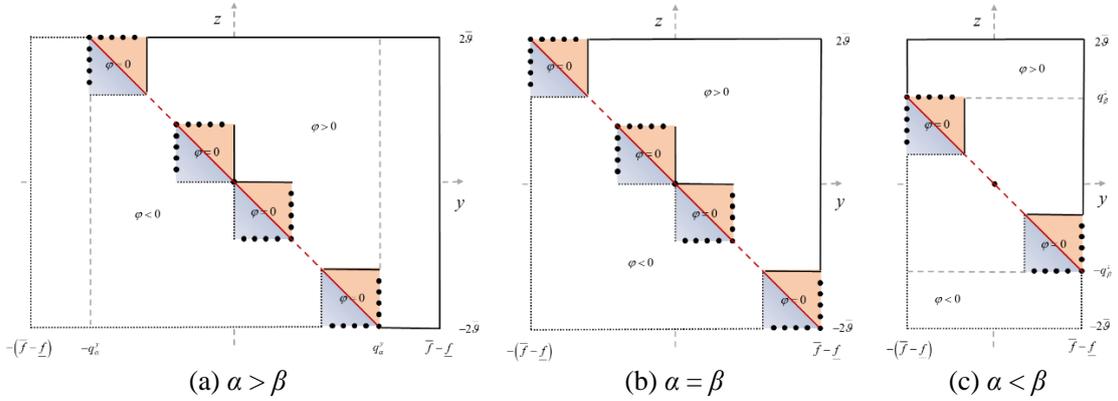

(a) $\alpha > \beta$      (b) $\alpha = \beta$      (c) $\alpha < \beta$

Fig. 6. Images of the line $L$ and corresponding divisions of the region where $\varphi = 0$.

4.2. Quantization function

In this subsection we first deduce the general conditions for the establishment of the quantization function. Then, a feasible quantization function is presented, and its feasibility is rigorously proven. Finally, all the involved parameters in the feasible quantization function are analyzed in detail.

4.2.1. Conditions

Suppose there exists a quantization function $\pi(\mathbf{x})$ that satisfies the comparison criteria proposed in Section 4.2.2, then given two arbitrary solutions $\mathbf{x}_1$ and $\mathbf{x}_2$ of a COP, it should meet the

following conditions

$$\begin{cases} \pi(\mathbf{x}_1) > \pi(\mathbf{x}_2), & \forall \mathbf{x}_1 \in S, \mathbf{x}_2 \in S, \text{ if } \varphi < 0 \\ \pi(\mathbf{x}_1) < \pi(\mathbf{x}_2), & \forall \mathbf{x}_1 \in S, \mathbf{x}_2 \in S, \text{ if } \varphi > 0 \\ \pi(\mathbf{x}_1) = \pi(\mathbf{x}_2), & \forall \mathbf{x}_1 \in S, \mathbf{x}_2 \in S, \text{ if } \varphi = 0 \text{ and } z = -e \cdot y \\ \pi(\mathbf{x}_1) > \pi(\mathbf{x}_2), & \forall \mathbf{x}_1 \in S, \mathbf{x}_2 \in S, \text{ if } \varphi = 0 \text{ and } z < -e \cdot y \\ \pi(\mathbf{x}_1) < \pi(\mathbf{x}_2), & \forall \mathbf{x}_1 \in S, \mathbf{x}_2 \in S, \text{ if } \varphi = 0 \text{ and } z > -e \cdot y \end{cases} \quad (28)$$

To specify the regions where $\varphi < 0$, $\varphi > 0$, and $\varphi = 0$, we need to determine the ordered sets $Q_\alpha^y = \langle q_1^y, q_2^y, \cdots q_\alpha^y \rangle$ and $Q_\beta^z = \langle q_1^z, q_2^z, \cdots q_\beta^z \rangle$ for the divisions of $y$ and $z$, respectively. To start with a simple and convincing case, let consider $\alpha = 1$ and $\beta = 2$. Moreover,

$$Q_1^y = \langle q_1^y = \eta \cdot (\overline{f} - \underline{f}) \rangle, \quad (29)$$

$$Q_2^z = \langle q_1^z = \eta \cdot \xi \cdot \overline{\vartheta}, q_2^z = \xi \cdot \overline{\vartheta} \rangle, \quad (30)$$

where $0 < \eta < 1$, and $0 < \xi \leq 1$ are two user-defined parameters. In this way, $\varphi < 0$, $\varphi > 0$, and $\varphi = 0$ can be specified as follows,

$$\varphi < 0 \Leftrightarrow \begin{cases} (1) \; -2\overline{\vartheta} \leq z < -\overline{\vartheta} \\ (2) \; -\overline{\vartheta} \leq z < -\xi \cdot \overline{\vartheta} \\ (3) \; -\xi \cdot \overline{\vartheta} \leq z < -\xi \cdot \eta \cdot \overline{\vartheta} \text{ and } -(\overline{f} - \underline{f}) \leq y \leq \eta \cdot (\overline{f} - \underline{f}) \\ (4) \; -\xi \cdot \eta \cdot \overline{\vartheta} \leq z \leq 0 \text{ and } -(\overline{f} - \underline{f}) \leq y \leq 0 \text{ and } z + y \neq 0 \\ (5) \; 0 < z \leq \xi \cdot \eta \cdot \overline{\vartheta} \text{ and } -(\overline{f} - \underline{f}) \leq y < -\eta \cdot (\overline{f} - \underline{f}) \end{cases} \quad (31)$$

$$\varphi > 0 \Leftrightarrow \begin{cases} (1) \; \overline{\vartheta} \leq z < 2\overline{\vartheta} \\ (2) \; \xi \cdot \overline{\vartheta} < z \leq \overline{\vartheta} \\ (3) \; \xi \cdot \eta \cdot \overline{\vartheta} < z \leq \xi \cdot \overline{\vartheta} \text{ and } -\eta \cdot (\overline{f} - \underline{f}) \leq y \leq (\overline{f} - \underline{f}) \\ (4) \; \xi \cdot \eta \cdot \overline{\vartheta} \leq z \leq 0 \text{ and } 0 \leq y \leq (\overline{f} - \underline{f}) \text{ and } z + y \neq 0 \\ (5) \; 0 < z \leq -\xi \cdot \eta \cdot \overline{\vartheta} \text{ and } \eta \cdot (\overline{f} - \underline{f}) < y \leq (\overline{f} - \underline{f}) \end{cases} \quad (32)$$

$$\varphi = 0 \Leftrightarrow \begin{cases} (1) \; \xi \cdot \eta \cdot \overline{\vartheta} < z \leq \xi \cdot \overline{\vartheta} \text{ and } -(\overline{f} - \underline{f}) \leq y < -\eta \cdot (\overline{f} - \underline{f}) \\ (2) \; 0 < z \leq \xi \cdot \eta \cdot \overline{\vartheta} \text{ and } -\eta \cdot (\overline{f} - \underline{f}) \leq y < 0 \\ (3) \; -\xi \cdot \eta \cdot \overline{\vartheta} \leq z < 0 \text{ and } 0 < y \leq \eta \cdot (\overline{f} - \underline{f}) \\ (4) \; -\xi \cdot \overline{\vartheta} \leq z < -\xi \cdot \eta \cdot \overline{\vartheta} \text{ and } \eta \cdot (\overline{f} - \underline{f}) < y \leq (\overline{f} - \underline{f}) \end{cases} \quad (33)$$

The corresponding regions on the two-dimensional Cartesian system for each sub-condition in Equations (31)-(33) are depicted in Fig. 7.

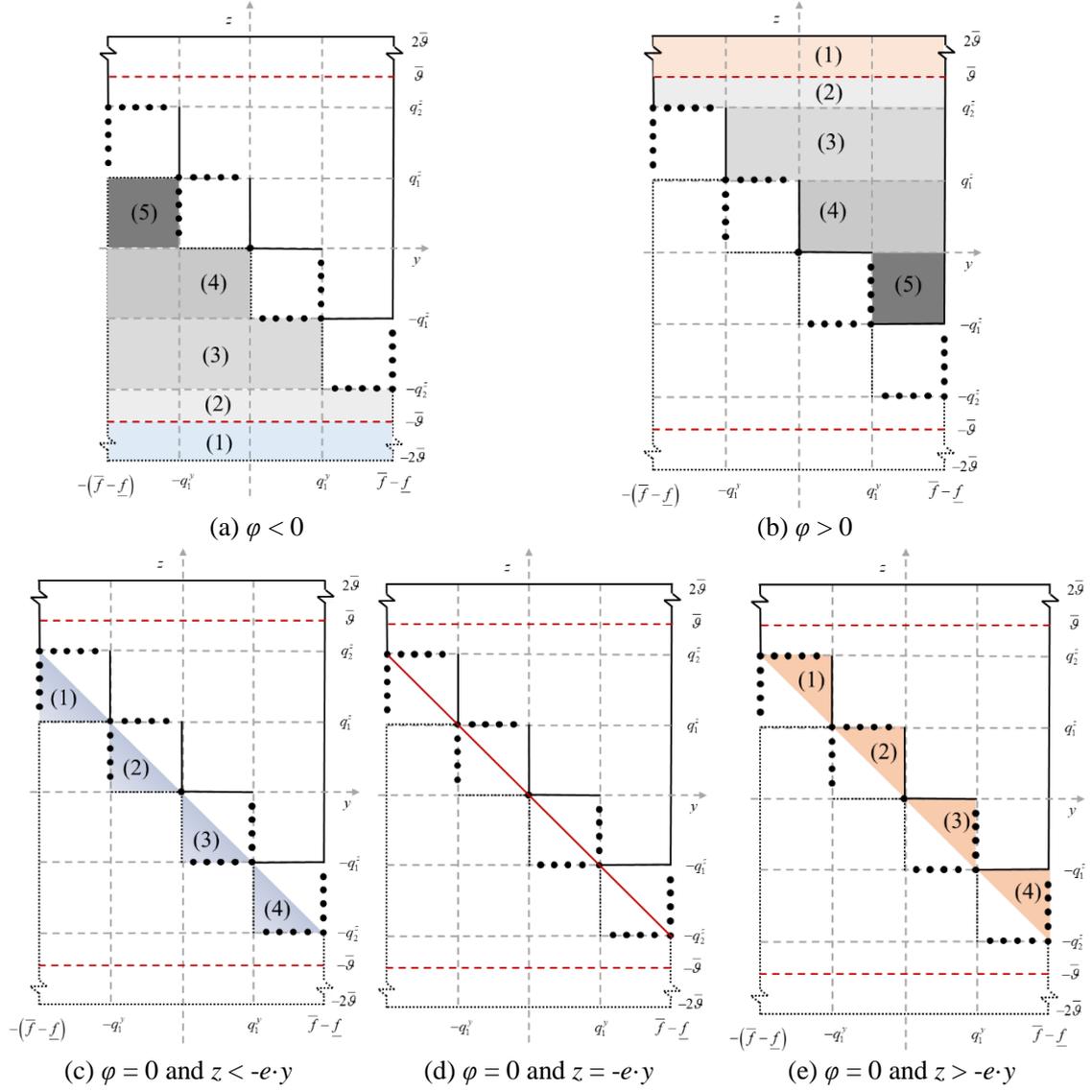

Fig. 7. Corresponding regions for each sub-condition in Equations (31)-(33)

4.2.2. A feasible solution

Here we present a feasible quantization function satisfying Equation (28). Such feasible quantization function has a basic form the same as penalty functions presented in Equation (5),

$$\pi(\mathbf{x}) = -(f(\mathbf{x}) + \sigma(\mathbf{x})), \tag{34}$$

where $\sigma(\mathbf{x})$ has the same role as $p(\mathbf{x})$ defined by Equation (6), but with a different expression.

**Lemma 2.** The following expression of $\sigma(\mathbf{x})$ makes $\pi(\mathbf{x})$ satisfy Equation (31),

$$\sigma(\mathbf{x}) = \mathbb{I}(\rho(\mathbf{x}) > 0) \cdot (\overline{f} - \underline{f}) + \frac{(\overline{f} - \underline{f})}{\xi \cdot \overline{\vartheta}} \theta(\mathbf{x}), \tag{35}$$

where

$$\rho(\mathbf{x}) = \sum_{i=1}^{n} \mathbb{I}(v_i(\mathbf{x}) > 0), \tag{36}$$

and

$$\mathbb{I}(\rho(\mathbf{x}) > 0) = \begin{cases} 1 & \text{if } \rho(\mathbf{x}) > 0 \\ 0 & \text{otherwises} \end{cases}. \tag{37}$$

**Lemma 3.** The expression of $\sigma(\mathbf{x})$ specificed by Equation (35) makes $\pi(\mathbf{x})$ satisfy Equation (32).

**Lemma 4.** The expression of $\sigma(\mathbf{x})$ specificed by Equation (35) makes $\pi(\mathbf{x})$ satisfy Equation (33) with $z < -e \cdot y$.

**Lemma 5.** The expression of $\sigma(\mathbf{x})$ specificed by Equation (35) makes $\pi(\mathbf{x})$ satisfy Equation (33) with $z = -e \cdot y$.

**Lemma 6.** The expression of $\sigma(\mathbf{x})$ specificed by Equation (35) makes $\pi(\mathbf{x})$ satisfy Equation (33) with $z > -e \cdot y$.

Lemmas 2-6 demonstrates that $\pi(\mathbf{x})$ specified by Equations (34)-(37) is a feasible quantization function satisfying Equation (28).

4.2.3. Control of $\xi$

As analyzed in Section 4.2.2, $\xi$ is a user-defined parameter used to specify the ordered set $Q_\beta^z = \langle q_1^z, q_2^z, \cdots q_\beta^z \rangle$ for the division of $z$. Besides, it is also a parameter needed to be set in advance for the use of Equation (35), namely the proposed quantization function. In what follows, we first discuss the effect of changes of $\xi$ on the division, and then propose an adaptive dynamic control method for $\xi$.

The line $L$ passes through the points $\left(-\left(\overline{f}-\underline{f}\right), q_2^z\right)$ and $(0, 0)$. Besides, $q_2^z = \xi \cdot \overline{\vartheta}$, and thus the change of $\xi$ changes the slope of $L$. As shown in Fig. 8, as $\xi$ gradually decreases from 1 to 0 (cannot reach 0), the slope of $L$ gradually tends to $-0$ from $\overline{\vartheta}/-\left(\overline{f}-\underline{f}\right)$, and the area of the region where $\varphi = 0$ also gradually tends to 0. When the area of the region where $\varphi = 0$ is relatively bigger, that is, $\xi$ is set larger, those infeasible solutions whose objective functions values are favored will have a greater chance of being selected as better alternatives.

As we all know, metaheuristic optimization algorithms generally search for optimal solutions in an iterative manner. In the early stage of the iteration, algorithms should relax the violation of the constraint conditions, and at the same time make full use of the information provided by objective function values, which can prevent algorithms from falling into a local optimum prematurely. In the later stage of the iteration, in order to ensure that algorithms can finally find a feasible solution, the violation of the constraint conditions should be treated strictly. To this end, to meet these preferences of algorithms at different stages of the iteration, it is necessary to reduce by $\xi$ during the iterative process.

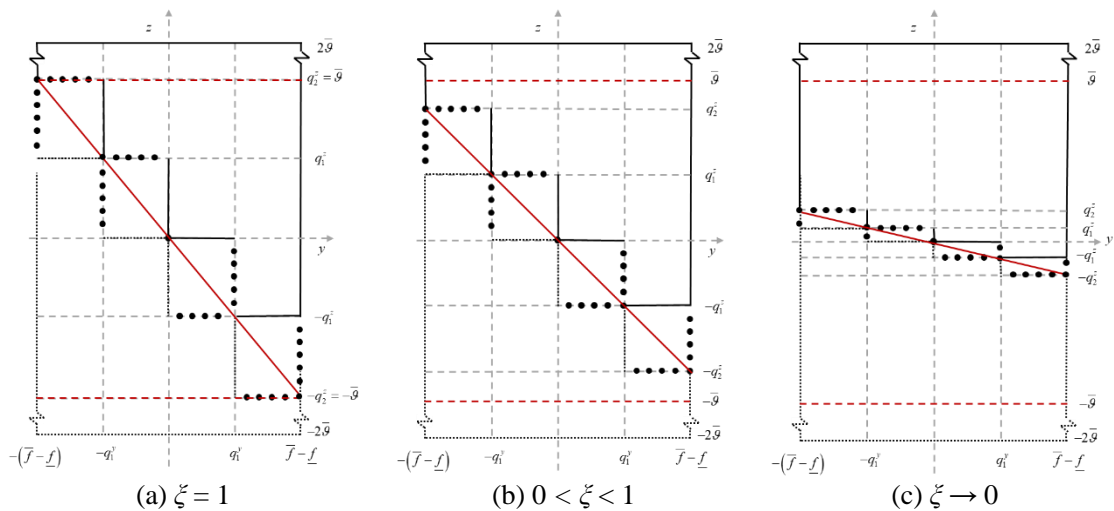

(a) $\xi = 1$  (b) $0 < \xi < 1$  (c) $\xi \to 0$

Fig. 8. Slope of *L* with different setting values of ξ

Let *t* denote the number of the iteration, $\bar{t}$ the maximum number of the iteration, ξ(*t*) the value of ξ at *t*-th iteration, $\bar{\xi}$ and $\underline{\xi}$ the maximum and minimum values of ξ(*t*), the following expression can be regarded as an alternative ones of ξ(*t*).

$$\xi(t) = \bar{\xi} - (\bar{\xi} - \underline{\xi})\left(\frac{t}{\bar{t}}\right)^p, \tag{38}$$

where *p* > 0 is an user-defined parameter used to control the deceleration rate of ξ. The larger *p* is, the faster ξ decreases.

4.2.4. Relaxation of $(\bar{f} - \underline{f})$ and $\bar{\vartheta}$

In Equation (35) we assume that difference between the maximum and minimum values of the objective function, namely $(\bar{f} - \underline{f})$, and the maximum value of constraint violation, namely $\bar{\vartheta}$, are known. Since the determination of $(\bar{f} - \underline{f})$ and $\bar{\vartheta}$ does not involve any constraints, in most cases we can use some existing optimization methods to determine $(\bar{f} - \underline{f})$ and $\bar{\vartheta}$ accurately. Nevertheless, it is a challenge to determine $(\bar{f} - \underline{f})$ and $\bar{\vartheta}$ in practical applications, especially when the objective function and constraint function are very complex. To this end, it is necessary to relax $(\bar{f} - \underline{f})$ and $\bar{\vartheta}$ such that we do not determine them exactly while Equation (35) still meets Equation (28).

Let the setting value of $(\bar{f} - \underline{f})$ is $f$, and the setting value of $\bar{\vartheta}$ is $\vartheta$. By re-examining the proofs of Lemmas 2-6 presented in Section II of Supplementary Material, the following corollaries can be obtained.

**Corollary 1.** If Equation (39) holds, Equation (35) with the setting values $f$ and $\vartheta$ meets Equation (31).

$$\frac{f}{\vartheta} = \frac{\bar{f} - \underline{f}}{\bar{\vartheta}} \tag{39}$$

**Corollary 2.** If Equation (39) holds, Equation (35) with the setting values $f$ and $\vartheta$ meets Equation (32).

**Corollary 3.** If Equation (39) holds, Equation (35) with the setting values $f$ and $\vartheta$ meets Equation (33) with $z < -e \cdot y$.

**Corollary 4.** If Equation (39) holds, Equation (35) with the setting values $f$ and $\vartheta$ satisfies Equation (33) with $z = -e \cdot y$.

**Corollary 5.** If Equation (39) holds, Equation (35) with the setting values $f$ and $\vartheta$ meets Equation (33) with $z > -e \cdot y$.

Corollaries 1-5 demonstrate that Equation (35) with the setting values $f$ and $\vartheta$ meet Equation (28) as well if Equation (39) holds. In this way, we don't need to know exactly $(\overline{f} - \underline{f})$ and $\overline{\vartheta}$, but just the ratio between them to still be able to use Equation (35). However, knowing the ratio is also a challenge in some cases, so we continue to relax $(\overline{f} - \underline{f})$ and $\overline{\vartheta}$.

Let's re-examine the proofs of Lemmas 2-6 again in the case where Equation (40) is established, and the following corollaries can be obtained.

$$\frac{\underline{f}}{\overline{\vartheta}} > \frac{\overline{f} - \underline{f}}{\overline{\vartheta}} \tag{40}$$

**Corollary 6.** If Equation (40) holds, Equation (35) with the setting values $f$ and $\vartheta$ meets the sub-conditions (1)-(4) of Equation (31).

**Corollary 7.** If Equation (40) holds, Equation (35) with the setting values $f$ and $\vartheta$ also meets the sub-conditions (1)-(4) of Equation (32).

Corollaries 6-7 indicate that the evaluation of $\mathbf{x}_1$ and $\mathbf{x}_2$ satisfies the proposed criterion if the pair ($\mathbf{x}_1$, $\mathbf{x}_2$) is mapped into the gray region in Fig. 9. Regarding the discussion on the remaining region, let's define a new line first.

$$\hat{L}: z = -\hat{e} \cdot y, \hat{e} = -\frac{\xi \cdot \vartheta}{f} \tag{41}$$

**Lemma 7.** If $z < -\hat{e} \cdot y$, Equation (35) with the setting values $f$ and $\vartheta$ satisfying Equation (40) meets the sub-condition (5) of Equations (31).

**Lemma 8.** If $z > -\hat{e} \cdot y$, Equation (35) with the setting values $f$ and $\vartheta$ satisfying Equation (40) meets the sub-condition (5) of Equations (32).

**Lemma 9.** If the following two inequalities hold, Equation (35) with the setting values $f$ and $\vartheta$ satisfying Equation (40) meets the conditions of $\varphi = 0$.

$$\begin{cases} z < -\hat{e} \cdot y \text{ or } z > -e \cdot y \text{ and } z > 0 \quad (1) \\ z > -\hat{e} \cdot y \text{ or } z < -e \cdot y \text{ and } z < 0 \quad (2) \end{cases} \tag{42}$$

Lemmas 7-9 demonstrate that the evaluation of $\mathbf{x}_1$ and $\mathbf{x}_2$ satisfies the proposed criterion if the pair ($\mathbf{x}_1$, $\mathbf{x}_2$) is mapped into the white region in Fig. 9. In summary, only when the pair ($\mathbf{x}_1$, $\mathbf{x}_2$) is mapped into the region between the line $L$ and $\hat{L}$, that is, the green region in Fig. 9, Equation (35) with the setting values $f$ and $\vartheta$ satisfying Equation (40) may doesn't work. Assume that a random pair ($\mathbf{x}_1$, $\mathbf{x}_2$) has the same probability of mapping to any point in the two-dimensional Cartesian coordinate system. In this way, we can use the area to approximate the error rate of evaluation. It can be easily obtained the area of the green region as,

$$o = \left(\overline{f} - \underline{f}\right)\left[\xi \cdot \overline{\vartheta} - \frac{\xi \cdot \vartheta}{f}\left(\overline{f} - \underline{f}\right)\right], \tag{43}$$

and the whole area of the mapped region comes as

$$o = 8\overline{\vartheta} \cdot \left(\overline{f} - \underline{f}\right). \tag{44}$$

The error rate of evaluation is estimated as

$$\mu = \frac{o}{o} = \frac{\xi}{8\overline{\vartheta}\left(\overline{f} - \underline{f}\right)}\left(\frac{\overline{\vartheta}}{\overline{f} - \underline{f}} - \frac{\vartheta}{f}\right). \tag{45}$$

From Equation (45) the two conclusions can be obtained:

$$\lim_{\vartheta/f \to \overline{\vartheta}/(\overline{f}-\underline{f})} \mu = 0 \tag{46}$$

$$\lim_{\xi \to 0} \mu = 0 \tag{47}$$

It can be seen from Equation (46) that as long as we set $\vartheta/f$ close enough to $\overline{\vartheta}/(\overline{f} - \underline{f})$, the evaluation error rate $\mu$ will be small enough. Equation (47) indicates that the evaluation error rate $\mu$ gradually decreases and tends to 0 with the use of Equation (38).

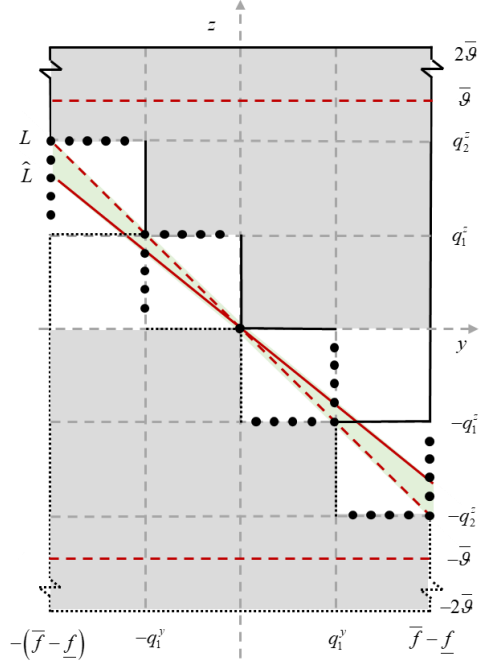

Fig. 9. Division with the two straight lines $L$ and $\hat{L}$

4.2.5. $\alpha$ and $\beta$

In the above analysis, we assume that $\alpha = 1$ and $\beta = 2$. In fact, we can find more cases about $\alpha$ and $\beta$, and the above analysis still holds true.

**Lemma 10.** When $\beta - \alpha = 1$ ($\alpha = 1, 2, \ldots$), and

$$Q_1^y = \left\langle q_1^y = \eta_1 \cdot \left(\overline{f} - \underline{f}\right), \cdots, q_\alpha^y = \eta_\alpha \cdot \left(\overline{f} - \underline{f}\right)\right\rangle, \tag{48}$$

where $0 < \eta_1 < \cdots < \eta_\alpha < 1$, if

$$Q_2^z = \left\langle q_1^z = \eta_1 \cdot \xi \cdot \bar{\vartheta}, \cdots q_\alpha^z = \eta_\alpha \cdot \xi \cdot \bar{\vartheta}, q_\beta^z = \xi \cdot \bar{\vartheta} \right\rangle, \qquad (49)$$

then the expression of $\sigma(\mathbf{x})$ specified by Equation (35) makes $\pi(\mathbf{x})$ satisfy Equation (28) as well.

The division of the mapping region in the case of $\beta - \alpha = 1$ ($\alpha = 1$) is shown in Fig. 10.(a). The conditions specfied by Equation (28) for the establishment of $\pi(\mathbf{x})$ are derived from the comparison criterion for the three regions, that is, the ones above, below, and on the red line (orange, blue, and red). In other words, if the three regions do not change under a new division, $\pi(\mathbf{x})$ is also established under this new division. The case of $\beta - \alpha = 1$ ($\alpha = 2, \ldots$) actually just increases the number of elements in $Q_1^y$ and $Q_2^z$ such that the division of the mapping region changes, namely a new division is generated. Suppose that $q_\tau^y$ is an added element of $Q_1^y$, and $q_\tau^z$ is an added element of $Q_1^z$ on the basis of the case of $\beta - \alpha = 1$ ($\alpha = 1$). If the intersection of the line $q_\tau^y$ and $-q_\tau^z$ is on the red line associated with the case of $\beta - \alpha = 1$ ($\alpha = 1$) as shown in Fig. 10.(b), namely Equation (49) holds, the three regions do not change, and thus $\pi(\mathbf{x})$ is still established in the case of $\beta - \alpha = 1$ ($\alpha = 2, \ldots$). Otherwise, as shown in Fig. 10.(c), the comparison criterion for the gray region changes so that $\pi(\mathbf{x})$ is not established in the new division.

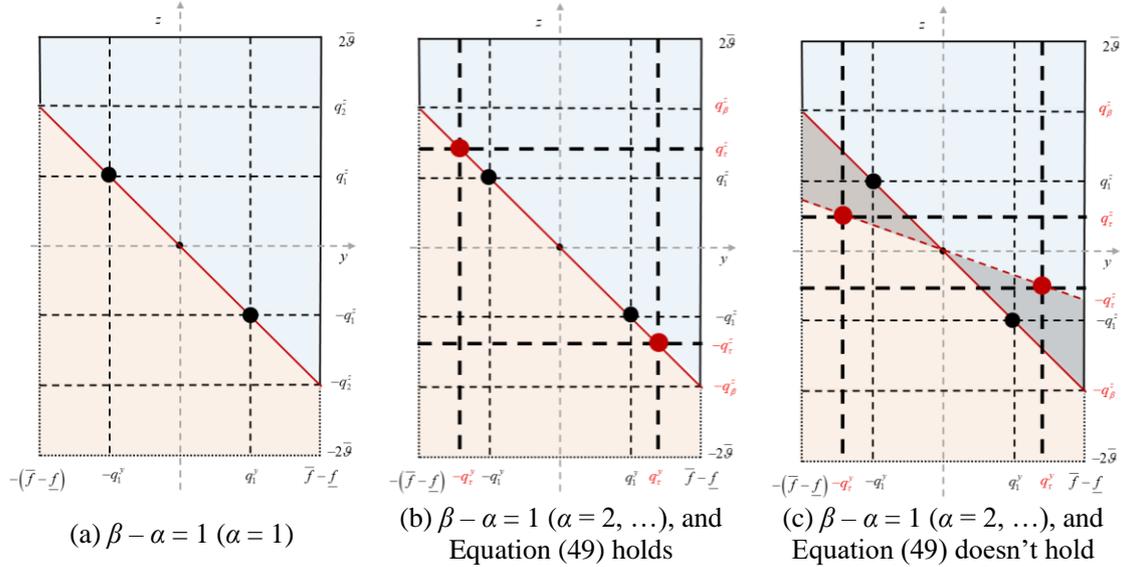

(a) $\beta - \alpha = 1$ ($\alpha = 1$)  
(b) $\beta - \alpha = 1$ ($\alpha = 2, \ldots$), and Equation (49) holds  
(c) $\beta - \alpha = 1$ ($\alpha = 2, \ldots$), and Equation (49) doesn't hold

Fig. 10. Division of the mapping region with different cases of $\alpha$ and $\beta$

Overall, this section presents the proposed CHT from both qualitative and quantitative aspects. Moreover, we analyze the parameters involved in the proposed CHT and derive guidelines for setting their values. The COP considered in the above analysis is single-objective, but the proposal can be easily extended to the case of multiple-objective COP. For more details, please see Section III of Supplement Material.

## 5. Computational experiments

This section describes the computational experiments conducted on two test datasets, namely the CEC 2006 (Liang et al., 2006) and the CEC 2017 (Wu et al., 2017). In the CEC 2006, 24 test cases (G01-G24) with decision variables varying from 2-24 are included, where G20 and G22 are not considered in this study since existing studies have not been able to find any feasible solutions for them. The CEC 2017 is specially designed for the era of Big Data in a way that the number of decision variables of its 28 test cases (C01-C28), denoted as $D$, is scalable.

5.1. Implementation details

All experiments were carried out on a PC with a Intel(R) Core(TM) i7-8700 CPU running at 3.20GHz, with 16.0 GB memory. We used Python 3.9 for implementing all the experiments. The raw data, source code, and detail results of the experiments are provided in a shared repository https://github.com/kobeIsMyYouth/Constraint-optimization.

Considering that the optimization algorithms used in this study have certain randomness, to eliminate the interference caused by such randomness as much as possible, the optimization algorithms ran repeatedly 25 times when they were used to solve COPs, and statistical values of the results of the 25 runs are considered as the evaluation criteria.

In a run if the optimization algorithm reaches the maximum number of function evaluations (FES), it will stop. The maximum number of FES for the experiments on the CEC 2006 was set to 500, 000, and that was $20,000 \times D$ for the experiments on the CEC 2017.

5.2. Experiments on the CEC 2006

In the algorithm (Takahama and Sakai, 2006) that won the first place in the competition held at 2006 IEEE Congress on Evolutionary Computation, a hybrid CHT combining the $\varepsilon$ level comparison, the feasible elites preserving strategy, and the gradient-based repair algorithm (Chootinan and Chen, 2006), was used to cope with constraints, and the differential evolution (DE) was used as the basic optimization tool. To follow such state-of-the-art algorithm, DE is also considered as the optimization tool in the experiments on the CEC 2006.

Other four CHTs shown in Table 2 are considered for comparing with the proposed one. Since it has been reported in existing studies that feasibility rules and the $\varepsilon$ level comparison perform significantly better than stochastic ranking (Mezura-Montes and Coello Coello, 2011), the latter is not included in the comparison. For a fair comparison, parameters settings associated with the basic optimization tool, DE algorithm, and those related to CHT3-CHT5, are the same as that in the study (Takahama and Sakai, 2006). As for CHT2, as analyzed in Section 3.2 there is no parameter to be set. Three parameters need to be set in advance for the use of CHT1, that is, $\bar{\xi}$, $\underline{\xi}$, and $p$ in Equation (38). It should be noted that $\underline{\xi}$ is always set to 0 and $p$ is 5 in this study. A detail experimental analysis on the setting of $\bar{\xi}$ is presented in Section 5.4.

In what follows, we discuss the performance of the five CHTs on the CEC 2006 in terms of accuracy and computational efficiency.

Table 2. The five CHTs to be compared

| Name | Description |
| --- | --- |
| CHT1 | The proposed CHT |
| CHT2 | The feasibility rules |
| CHT3 | The $\varepsilon$ level comparison |
| CHT4 | The hybrid CHT combining the $\varepsilon$ level comparison and the feasible elites preserving |

| | strategy |
|---|---|
| CHT5 | The hybrid CHT used in (Takahama and Sakai, 2006) |

5.2.1. Accuracy

For evaluating the accuracy of the optimization algorithms with different CHTs, the feasible rate (FR) and success rate (SR) are calculated after the 25 runs. In a run if the optimization algorithms find a feasible solution, the run is called a feasible one and FR rate refers to the proportion of the number of feasible runs to 25. If the obtained feasible solution $\mathbf{x}$ satisfies $f(\mathbf{x}) - f(\mathbf{x}^*) \leq 0.0001$ where $\mathbf{x}^*$ denotes the best-known solution reported in prior studies (Liang et al., 2006), such run is called a successful one and the proportion of the number of successful runs to 25 is regarded as SR. Table 3 shows the specific FR and SR of 22 test cases under the five considered CHTs. It can be seen that CHT1 and CHT5 achieved 100% FR and SR on all the 22 test cases, which outperforms the other three CHTs. By comparing FR and SR of CHT3-CHT5, we can see that a hybrid CHT (such as CHT4 and CHT5) is always superior to a pure CHT (such as CHT3), and the hybrid CHT seems to be able to take advantage of different components such that an accuracy improvement is expected to be obtained.

Table 3. FR and SR of 24 test cases under the five considered CHTs

| COP | CHT | FR | SR | COP | CHT | FR | SR |
|---|---|---|---|---|---|---|---|
| G01 | CHT1 | 100% | 100% | G12 | CHT1 | 100% | 100% |
| | CHT2 | 100% | 100% | | CHT2 | 100% | 100% |
| | CHT3 | 100% | 100% | | CHT3 | 100% | 100% |
| | CHT4 | 100% | 100% | | CHT4 | 100% | 100% |
| | CHT5 | 100% | 100% | | CHT5 | 100% | 100% |
| G02 | CHT1 | 100% | 100% | G13 | CHT1 | 100% | 100% |
| | CHT2 | 100% | 100% | | CHT2 | 100% | **40%** |
| | CHT3 | 100% | 100% | | CHT3 | 100% | 100% |
| | CHT4 | 100% | **96%** | | CHT4 | 100% | 100% |
| | CHT5 | 100% | 100% | | CHT5 | 100% | 100% |
| G03 | CHT1 | 100% | 100% | G14 | CHT1 | 100% | 100% |
| | CHT2 | 100% | 100% | | CHT2 | 100% | 100% |
| | CHT3 | 100% | 100% | | CHT3 | 100% | 100% |
| | CHT4 | 100% | 100% | | CHT4 | 100% | 100% |
| | CHT5 | 100% | 100% | | CHT5 | 100% | 100% |
| G04 | CHT1 | 100% | 100% | G15 | CHT1 | 100% | 100% |
| | CHT2 | 100% | 100% | | CHT2 | 100% | 100% |
| | CHT3 | 100% | 100% | | CHT3 | **52%** | **52%** |
| | CHT4 | 100% | 100% | | CHT4 | 100% | 100% |
| | CHT5 | 100% | 100% | | CHT5 | 100% | 100% |
| G05 | CHT1 | 100% | 100% | G16 | CHT1 | 100% | 100% |
| | CHT2 | 100% | 100% | | CHT2 | 100% | 100% |
| | CHT3 | 100% | 100% | | CHT3 | **96%** | **92%** |
| | CHT4 | 100% | 100% | | CHT4 | 100% | 100% |
| | CHT5 | 100% | 100% | | CHT5 | 100% | 100% |
| G06 | CHT1 | 100% | 100% | G17 | CHT1 | 100% | 100% |
| | CHT2 | 100% | 100% | | CHT2 | 100% | **32%** |
| | CHT3 | 100% | 100% | | CHT3 | 100% | **96%** |
| | CHT4 | 100% | 100% | | CHT4 | 100% | 100% |
| | CHT5 | 100% | 100% | | CHT5 | 100% | 100% |
| G07 | CHT1 | 100% | 100% | G18 | CHT1 | 100% | 100% |
| | CHT2 | 100% | 100% | | CHT2 | 100% | 100% |
| | CHT3 | 100% | 100% | | CHT3 | 100% | 100% |
|       | CHT4  | 100%  | 100%  |       | CHT4  | 100%  | 100%  |
|       | CHT5  | 100%  | 100%  |       | CHT5  | 100%  | 100%  |
| G08   | CHT1  | 100%  | 100%  | G19   | CHT1  | 100%  | 100%  |
|       | CHT2  | 100%  | 100%  |       | CHT2  | 100%  | 100%  |
|       | CHT3  | 100%  | **24%** |     | CHT3  | 100%  | **92%** |
|       | CHT4  | 100%  | 100%  |       | CHT4  | 100%  | **96%** |
|       | CHT5  | 100%  | 100%  |       | CHT5  | 100%  | 100%  |
| G09   | CHT1  | 100%  | 100%  | G21   | CHT1  | 100%  | 100%  |
|       | CHT2  | 100%  | 100%  |       | CHT2  | 100%  | **88%** |
|       | CHT3  | 100%  | 100%  |       | CHT3  | 100%  | 100%  |
|       | CHT4  | 100%  | 100%  |       | CHT4  | 100%  | 100%  |
|       | CHT5  | 100%  | 100%  |       | CHT5  | 100%  | 100%  |
| G10   | CHT1  | 100%  | 100%  | G23   | CHT1  | 100%  | 100%  |
|       | CHT2  | 100%  | 100%  |       | CHT2  | 100%  | 100%  |
|       | CHT3  | 100%  | 100%  |       | CHT3  | 100%  | 100%  |
|       | CHT4  | 100%  | 100%  |       | CHT4  | 100%  | **96%** |
|       | CHT5  | 100%  | 100%  |       | CHT5  | 100%  | 100%  |
| G11   | CHT1  | 100%  | 100%  | G24   | CHT1  | 100%  | 100%  |
|       | CHT2  | 100%  | 100%  |       | CHT2  | 100%  | 100%  |
|       | CHT3  | 100%  | 100%  |       | CHT3  | 100%  | 100%  |
|       | CHT4  | 100%  | 100%  |       | CHT4  | 100%  | 100%  |
|       | CHT5  | 100%  | 100%  |       | CHT5  | 100%  | 100%  |

5.2.2. Computational efficiency

We evaluate the computational efficiency of different CHTs in terms of two aspects, the minimum FES required to obtain a successful solution and time consumption for solving COPs. To this end, the value of $f(\mathbf{x}) - f(\mathbf{x}^*)$ is recorded in the successive FES, and the minimum FES required to satisfy $f(\mathbf{x}) - f(\mathbf{x}^*) \leq 0.0001$ in a run is subsequently obtained. The average minimum FES is calculated after all the 25 runs are completed, and the results are shown in Table 4. For more details about the minimum FES, one can refer to Section IV of Supplementary Material. The smaller the average minimum FES, the faster the corresponding CHT guides the algorithm to find a successful solution, and thus the higher the computational efficiency. From Table 4 we can see that CHT1 achieve better average minimum FES in 10 out of 22 COPs compared with other CHTs. CHT5 obtain an optimal average minimum FES in 9 out of 22 COPs, which achieves a close performance to CHT1. Moreover, as can be seen in the last row of Table 4, CHT1 is significantly better than the other four CHTs from the average performance on the 22 COPs.

Table 4. Minimum FES required to obtain a successful solution for different CHTs

| FES | CHT1 | CHT2 | CHT3 | CHT4 | CHT5 |
|-----|------|------|------|------|------|
| G01 | 58504 | **58211** | 124000 | 129148 | 59308 |
| G02 | **148382** | 150736 | 150972 | 148939 | 149825 |
| G03 | **27231** | 27242 | 94386 | 95585 | 89407 |
| G04 | 28042 | 27939 | 50047 | 48536 | **26216** |
| G05 | **130534** | 89984 | 98729 | 98709 | 97431 |
| G06 | 7778 | 7692 | 99853 | 101012 | **7381** |
| G07 | 79758 | 82129 | 144931 | 149880 | **74303** |
| G08 | 1219 | **1064** | 105988 | 55105 | 1139 |
| G09 | 24113 | 24328 | 99566 | 102178 | **23121** |
| G10 | 114261 | 113681 | 176284 | 176284 | **105234** |
| G11 | **2027** | 12466 | 99564 | 99291 | 16420 |
| G12 | 2746 | 4162 | 2396 | **2310** | 4124 |
| G13 | **14980** | 266791 | 97754 | 97358 | 34738 |
| G14 | **72183** | 78072 | 135104 | 135883 | 113439 |
| G15 | **4733** | 73723 | 98981 | 98770 | 84216 |
| G16 | 13580 | 13459 | 100320 | 100210 | **12986** |
| G17 | **34802** | 254324 | 109262 | 107941 | 98861 |
| G18 | 64159 | 63998 | 130636 | 138159 | **59153** |
| G19 | 378675 | 382478 | 379839 | 383742 | **356350** |
| G21 | **62648** | 110630 | 169144 | 164349 | 135143 |
| G23 | **145508** | 257916 | 275238 | 252456 | 200765 |
| G24 | 3162 | 3185 | 3091 | 3082 | **2952** |
| Avg | **64501** | 95646 | 124822 | 122224 | 79660 |

The core computing task of using CHTs to solve COPs is to find an "optimal solution" from $l$ alternative solutions. In this case, different CHTs' time consumption used to complete the core computing task is similar since their computational complexity is the same as $O(l)$. However, in the design of some metaheuristic optimization algorithms, it is necessary to sort the performance of the $l$ solutions, such as the non-dominated sorting in NSGA-II (Deb et al., 2002). In such a case, the computational complexity of the pairwise comparison-based CHT (CHT2-CHT5) becomes $O(l\log l)$ or $O(l^2)$ (depending on the sorting algorithm used), while that of CHT1 is still $O(l)$ because the evaluation function value itself can be used for sorting. Fig. 11 shows the total time consumption for sorting $l$ solutions in the considered 22 COPs, where $l$ is set to a value taken from 10 to 200 in steps of 10. The values in the interval [10, 200] is considered mainly because the size of the set of alternative solutions in metaheuristic algorithms generally falls within [10, 200]. It can be seen that the evaluation function-based CHT (CHT1) exhibits a linear time growth with size ($l$), while the time consumption growth rate of the evaluation function-based CHTs increases continuously with the increase of size.

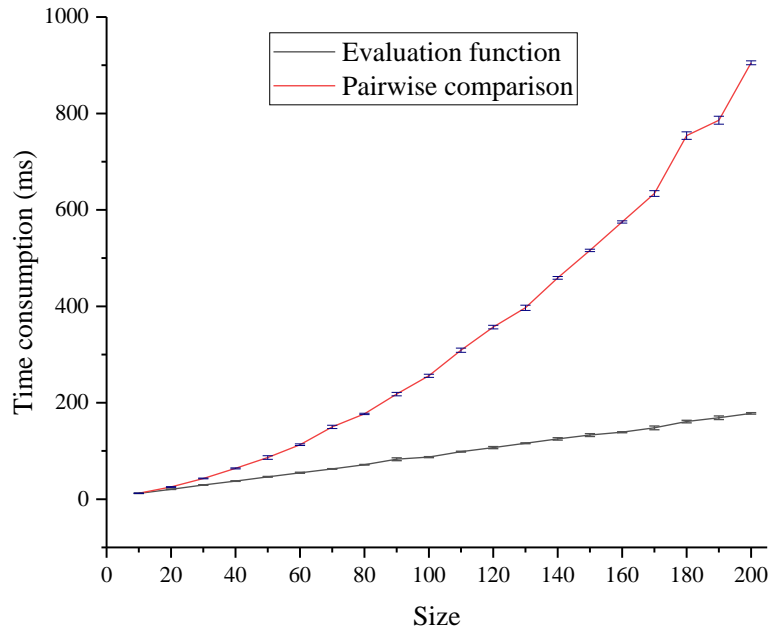

Fig. 11. Time consumption for evaluation function-based and pairwise comparison-based CHTs

5.3. Experiments on the CEC 2017

In the competition held at 2017 IEEE Congress on Evolutionary Computation, the LSHADE44

algorithm (Poláková, 2017) won the first place. The LSHADE44 algorithm employed four different DE strategies, labelled as current-to-pbest/1/bin, current-to-pbest/1/exp, randrl/1/bin, and randrl/1/exp. A competition mechanism was used to select a strategy for computing a trial point. Besides, a self-adaptive parameter determination method based on four pairs of historical circle memories and an archive strategy were also adopted. Regarding the CHT, the simple feasible rules (CHT2) was employed. Nevertheless, the LSHADE44 algorithm is much more powerful than the standard DE, and it has strong global optimization capabilities on high-dimensional COPs. In this study, we intentionally reduce the complexity of the optimization algorithm by discarding some designs, including the strategies competition, the parameter adaptation, and the archive strategy, resulting in a DE variant that only uses randrl/1/exp.

Four combinations of optimization algorithms and CHTs, that is, simplified DE+CHT1 (A1), simplified DE+CHT2 (A2), simplified DE+CHT3 (A3), and LSHADE44 + CHT2 (A4), were used for comparison on CEC 2017. In particular, the last one is what the study (Poláková, 2017) used, and it is considered mainly to investigate whether simple optimization algorithms can achieve or exceed the performance of complex one under the premise of using the proposed CHT. This is crucial if it is true, since we no longer need to pursue complex algorithm designs, which can greatly improve the computational efficiency of existing metaheuristic optimization algorithms.

5.3.1. Accuracy

The CEC 2017 dataset does not provide an optimal solution for each COP, so we cannot obtain SR. Besides FR, other five statistical indicators, that is, the best, median, worst objective function value obtained in 25 runs, as well as the mean and standard deviation, are introduced for accuracy evaluation. When FR = 100%, all the obtained results in 25 runs are involved in the calculation of the five statistical indicators. When 0 < FR < 100%, if a run is not feasible the obtained results are not involved in the calculation. When FR = 0, all the obtained results in 25 runs are involved in the calculation, and the constraint violations are also recorded for follow-up analysis. In this way, for each COP we can obtain the FR, Best, Median, Worst, Mean and STD of A1-A4, and the results are shown in Section V of Supplementary Material. A pairwise comparison between A1 and A2/A3/A4 is further conducted on the six indicators for a comprehensive evaluation on the 28 COPs, and the comparison results are recorded as "=", "+", "–", respectively. More specifically, if A1 has the same value as A2/A3/A4 on a certain indicator in a COP, then the corresponding count of "=" will be incremented by 1. Similarly, if A1 has a better (worse) value as A2/A3/A4 on a certain indicator in a COP, the count of "+" ("–") will be incremented by 1. Specificity, in the case of FR = 0, only when both the objective function value and constraint violation are better (worse) can the count of "+" ("–") be incremented by 1, otherwise the count of "=" is incremented by 1. The counting results of COPs with $D$ = 10, 30, 50, and 100 are presented in Tables 5-8.

It can be seen that in all the case of $D$ = 10, 30, 50, and 100 the counting results of A1 are better than those of A2 and A3 in terms of all the six indicators (FR, Best, Median, Worst, Mean and STD), which indicates that the proposed CHT significantly outperforms the feasibility rules and the $\varepsilon$ level comparison. In the comparison between A1 and A4, except for $D$ = 50 and $D$ = 100, the counting results of A1 in terms of FR are slightly worse than that of A4 (3vs5, and 2vs5), A1 is superior in other cases, which demonstrates that simple optimization algorithms combined with the proposed CHT can achieve similar even better performance to the complex one.

Table 5. Counting results of COPs with $D = 10$

|  |  | FR | Best | Median | Worst | Mean | STD |
|---|---|---|---|---|---|---|---|
| A1 | = | 20 | 16 | 13 | 12 | 10 | 9 |
| Vs | + | 6 | 10 | 13 | 12 | 13 | 13 |
| A2 | − | 2 | 2 | 2 | 4 | 5 | 6 |
| A1 | = | 22 | 22 | 15 | 11 | 10 | 9 |
| Vs | + | 5 | 4 | 9 | 9 | 10 | 11 |
| A3 | − | 1 | 2 | 4 | 8 | 8 | 8 |
| A1 | = | 21 | 13 | 10 | 7 | 6 | 6 |
| Vs | + | 5 | 12 | 14 | 15 | 15 | 14 |
| A4 | − | 2 | 3 | 4 | 6 | 7 | 8 |

Table 6. Counting results of COPs with $D = 30$

|  |  | FR | Best | Median | Worst | Mean | STD |
|---|---|---|---|---|---|---|---|
| A1 | = | 21 | 14 | 10 | 6 | 8 | 3 |
| Vs | + | 6 | 14 | 17 | 18 | 17 | 16 |
| A2 | − | 1 | 0 | 1 | 4 | 3 | 9 |
| A1 | = | 21 | 12 | 11 | 7 | 5 | 6 |
| Vs | + | 6 | 9 | 11 | 14 | 15 | 15 |
| A3 | − | 1 | 7 | 6 | 7 | 8 | 7 |
| A1 | = | 20 | 6 | 4 | 2 | 2 | 3 |
| Vs | + | 5 | 17 | 16 | 16 | 15 | 13 |
| A4 | − | 3 | 5 | 8 | 10 | 11 | 12 |

Table 7. Counting results of COPs with $D = 50$

|  |  | FR | Best | Median | Worst | Mean | STD |
|---|---|---|---|---|---|---|---|
| A1 | = | 20 | 9 | 7 | 7 | 7 | 3 |
| Vs | + | 7 | 18 | 19 | 20 | 19 | 20 |
| A2 | − | 1 | 1 | 2 | 1 | 2 | 5 |
| A1 | = | 21 | 10 | 8 | 6 | 5 | 3 |
| Vs | + | 4 | 10 | 13 | 15 | 15 | 16 |
| A3 | − | 3 | 8 | 7 | 7 | 8 | 9 |
| A1 | = | 21 | 3 | 1 | 2 | 1 | 1 |
| Vs | + | 3 | 19 | 17 | 17 | 18 | 16 |
| A4 | − | 4 | 6 | 10 | 9 | 9 | 11 |

Table 8. Counting results of COPs with $D = 100$

|  |  | FR | Best | Median | Worst | Mean | STD |
|---|---|---|---|---|---|---|---|
| A1 | = | 25 | 10 | 6 | 6 | 6 | 3 |
| Vs | + | 3 | 17 | 19 | 19 | 19 | 15 |
| A2 | − | 0 | 1 | 3 | 3 | 3 | 10 |
| A1 | = | 19 | 12 | 9 | 6 | 6 | 5 |
| Vs | + | 5 | 8 | 12 | 16 | 13 | 17 |
| A3 | − | 4 | 8 | 7 | 6 | 9 | 6 |
| A1 | = | 21 | 1 | 0 | 1 | 0 | 0 |
| Vs | + | 2 | 18 | 17 | 15 | 17 | 15 |
| A4 | − | 5 | 9 | 11 | 12 | 11 | 13 |

5.3.2. Computational efficiency

Due to the lack of optimal solutions of COPs in the CEC 2017 dataset, the analysis on the minimum FES required to obtain a successful solution cannot be carried out. In what follows, we present the results of two types of computational experiments, 1) the comparison between the time consumption of pairwise comparison-based and evaluation function-based CHTs for sorting $l$ solutions in the considered 28 COPs, where $l$ takes values on the interval [10, 200] with a step size of 10, and 2) the comparison between the time consumption of pairwise comparison-based and evaluation function-based CHTs for sorting a fixed number of solutions with different $D$ (10, 30, 50, 100). The results of the two types of computational experiments are shown in Figs. 12 and 13, respectively.

From Fig. 12 we can see that evaluation function-based CHTs always have a lower time consumption than pairwise comparison-based ones under different size of solutions. Besides, as the size increases, the time consumption gap between evaluation function-based and pairwise comparison-based CHTs becomes larger. We can also see that $D$ of COPs does not seem to affect the above observations.

From Fig. 13 we can see that the time consumption of evaluation function-based and pairwise comparison-based CHTs shows a non-linear growth with the increase of $D$. Moreover, the time consumption of evaluation function-based CHTs is always lower than that of pairwise comparison-based CHTs under different $D$, and the gap increases with the increase of $D$. We can also see that size of solutions does not seem to affect the above observations.

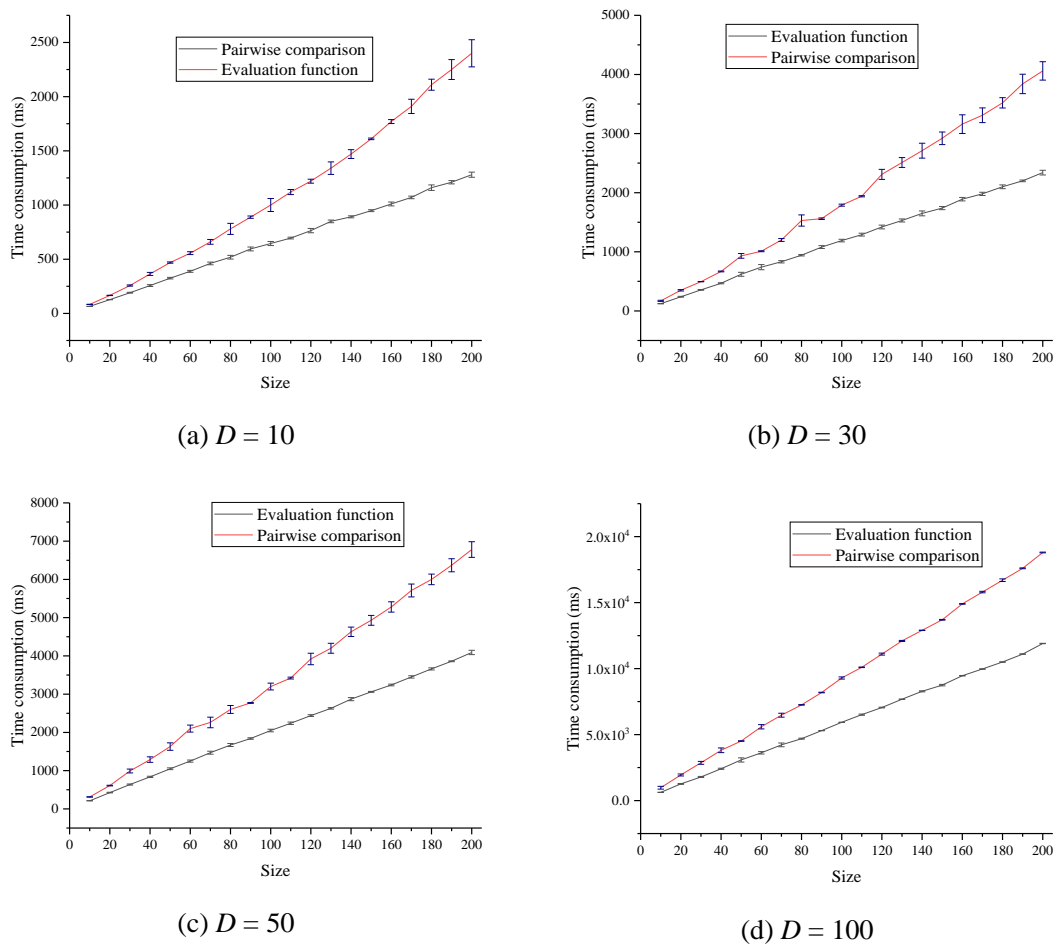

(a) $D = 10$     (b) $D = 30$

(c) $D = 50$     (d) $D = 100$

Fig. 12. Comparison results between time consumption of pairwise comparison-based and evaluation function-based CHTs for sorting $l$ solutions in the considered 28 COPs

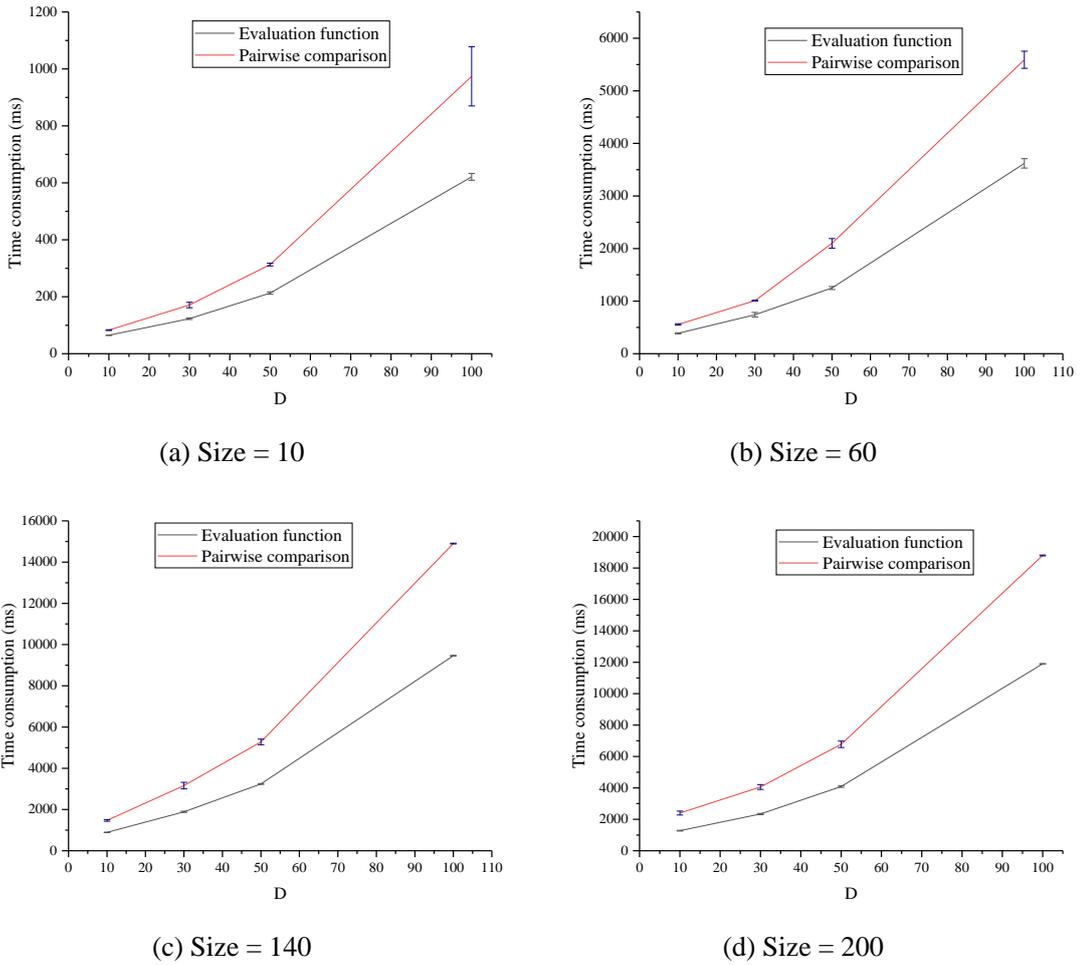

(a) Size = 10  (b) Size = 60

(c) Size = 140  (d) Size = 200

Fig. 13. Comparison results between time consumption of pairwise comparison-based and evaluation function-based CHTs for sorting a fixed number of solutions with different $D$

5.4. Experimental analysis for parameter $\xi$

Only two parameters are involved in the proposed CHT, that is, $p$ in Equation (38) and $\xi$ in Equation (35). In most of cases, setting $p = 5$ and $\xi = 1$ works very well, but in a few cases when metaheuristic optimization algorithms cannot obtain a feasible solution, we need to adjust the value of $\xi$. To this end, we fixed the value of $p$, and set $\xi$ as 1, 0.1, 0.01, 0.001, 0.0001, and 0.00001, respectively, for a parameter analysis. The COP C08 in the CEC 2017 is selected to conduct such an analysis. To solve C08, we repeatedly run A4 with different $\xi$ 25 times, and record the optimal objective function value and constraint violation obtained by each FES during each run. In this way, Fig. 14 that depicts the variations of the objective function value and constraint violation with FES in 25 runs, can be obtained. From Fig. 14 we can see that as $\xi$ is set smaller and smaller, the converged objective function value and constraint violation are also smaller and smaller. When $\xi$ is set small enough (e.g., $\xi = 0.0001$), the constraint violation eventually converges to 0, that is, we have obtained a feasible solution. Besides, from the comparison between Fig. 14.(e) and 14.(f), it can be seen that on the premise that a feasible solution can be obtained continuously reducing the value of $\xi$ can increase the convergence speed. But this is not always the case. During the experiments in some cases, especially when it is easy to obtain a feasible solution, reducing the value of $\xi$ may increase the value of the final converged objective function.

The above experimental analysis results give us the following guidance of setting $\xi$.
(1) On the premise that a feasible solution can be obtained, $\xi$ should be set as large as possible so that the converged objective function value would be better.
(2) When it is difficult to obtain a feasible solution, $\xi$ should be set as small as possible so that a feasible solution can be finally obtained.

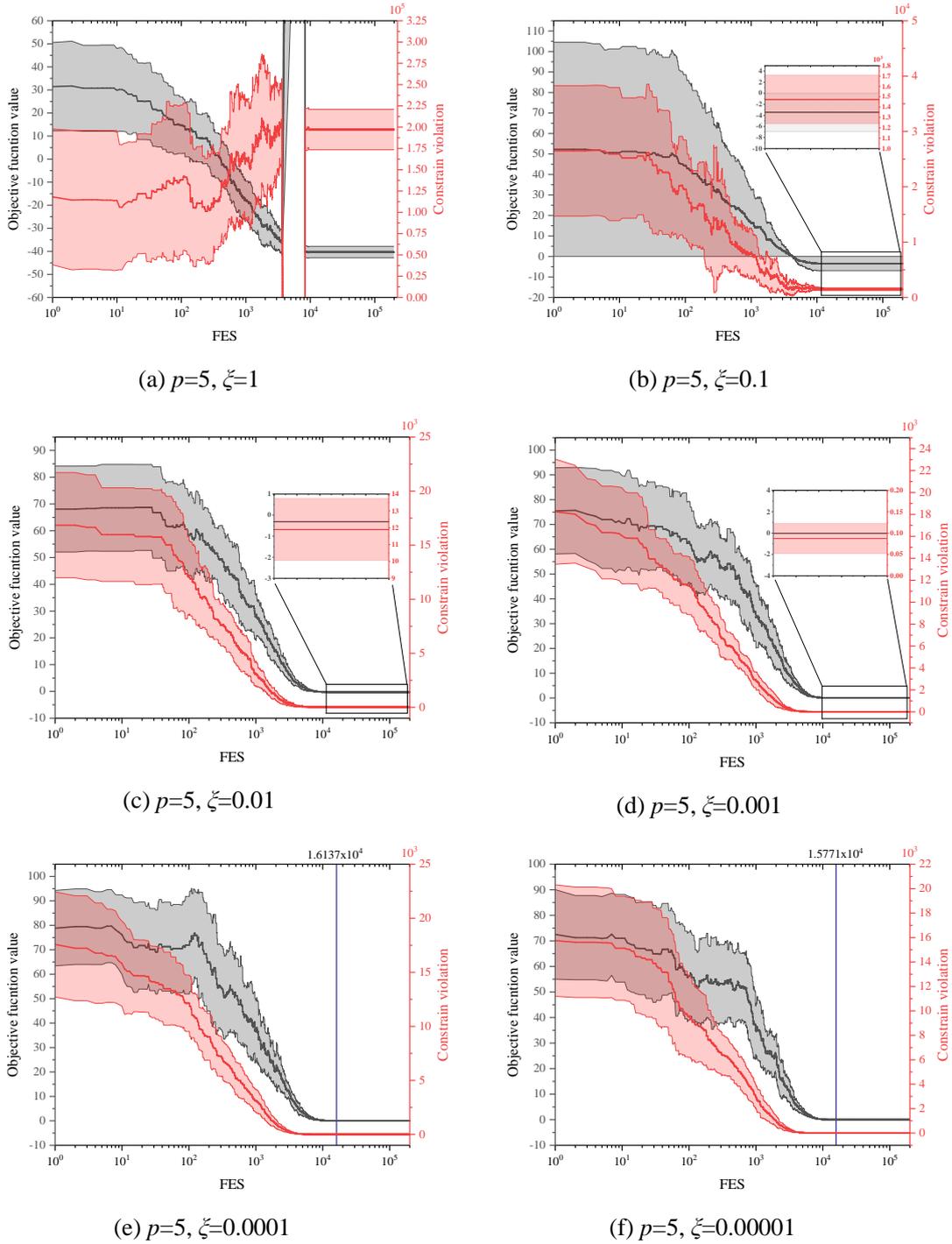

(a) $p=5$, $\xi=1$

(b) $p=5$, $\xi=0.1$

(c) $p=5$, $\xi=0.01$

(d) $p=5$, $\xi=0.001$

(e) $p=5$, $\xi=0.0001$

(f) $p=5$, $\xi=0.00001$

Fig. 14. Variations of objective functions values and constraint violations with different $\xi$ and a fixed $p$.

## 6. Conclusion

In this study we propose a quantitative pairwise comparison-based constraint handling technique for constrained optimization. Compared with existing techniques, the proposed one enjoy higher evaluation accuracy and computational efficiency, which help solve constrained optimization problems more effectively and efficiently. The versatility of our technique allows seamless integration with various metaheuristic optimization algorithms, such as, tabu search, evolutionary algorithms, swarm intelligence algorithms, and so on, enabling the resolution of diverse constrained optimization problems. Thorough theoretical and experimental analyses have been conducted on the involved parameters, facilitating the identification of their optimal values with ease. Besides, this study introduces a fresh theoretical framework for effectively balancing objective function values and constraint violations in constrained optimizations, with the proposed technique serving as a distinctive case within this framework. Furthermore, this theoretical foundation opens up possibilities for developing more effective and efficient similar techniques for solving constrained optimization problems.

## Acknowledgments

This research was supported by the National Natural Science Foundation of China (Nos. 72188101 and 72171069).

# Supplementary Material of
## "A quantitative pairwise comparison-based constrained handling technique for constrained optimization"

Ting Huang[1,2,3], Qiang Zhang[1,2,3*], Witold Pedrycz[4], Shanlin Yang[1,2,3]

[1]*School of Management, Hefei University of Technology, Hefei, Box 270, Hefei 230009, Anhui, P.R. China;*

[2]*Key Laboratory of Process Optimization and Intelligent Decision-making, Ministry of Education, Hefei, Box 270, Hefei 230009, Anhui, P.R. China;*

[3]*Ministry of Education Engineering Research Center for Intelligent Decision-Making & Information System Technologies, Hefei 230009, Anhui, P.R. China;*

[4]*Department of Electrical & Computer Engineering, University of Alberta, Edmonton, AB T6G 2R3, Canada*

## I. Introduction

This supplementary material provides (1) proofs of lemmas (Section II); (2) extension to multi-objective COP (Section III); (3) boxplots of the minimum required FES for CHT1-CHT4 (Section IV); (4) statistical results of A1-A4 (Section V).

## II. Proofs of lemmas

**Proof for Lemma 1**:

Let's first discuss the domains of $z$ in the following four cases.

(1) If $\mathbf{x}_1 \in F$, $\mathbf{x}_2 \in F$, then $\vartheta(\mathbf{x}_1) = \vartheta(\mathbf{x}_1) = 0$, and thus $\theta(\mathbf{x}_1) = \theta(\mathbf{x}_2) = -\bar{\vartheta}$, namely $z = 0$;

(2) If $\mathbf{x}_1 \in F$, $\mathbf{x}_2 \in S - F$, then $\vartheta(\mathbf{x}_1) = 0$, $\vartheta(\mathbf{x}_2) \in (0, \bar{\vartheta}]$, and thus $\theta(\mathbf{x}_1) = -\bar{\vartheta}$, $\theta(\mathbf{x}_2) = \vartheta(\mathbf{x}_2) \in (0, \bar{\vartheta}]$, namely $z \in (-\bar{\vartheta}, -2\bar{\vartheta}]$;

(3) If $\mathbf{x}_1 \in S - F$, $\mathbf{x}_2 \in F$, then $\vartheta(\mathbf{x}_1) \in (0, \bar{\vartheta}]$, $\vartheta(\mathbf{x}_2) = 0$, and thus $\theta(\mathbf{x}_1) = \vartheta(\mathbf{x}_1) \in (0, \bar{\vartheta}]$, $\theta(\mathbf{x}_2) = -\bar{\vartheta}$, namely $z \in (\bar{\vartheta}, 2\bar{\vartheta}]$;

(4) If $\mathbf{x}_1 \in S - F$, $\mathbf{x}_2 \in S - F$, then $\vartheta(\mathbf{x}_1) \in (0, \bar{\vartheta}]$, $\vartheta(\mathbf{x}_2) \in (0, \bar{\vartheta}]$, and thus $\theta(\mathbf{x}_1) = \vartheta(\mathbf{x}_1) \in (0, \bar{\vartheta}]$, $\theta(\mathbf{x}_2) = \vartheta(\mathbf{x}_2) \in (0, \bar{\vartheta}]$, namely $z \in [0, \bar{\vartheta}) \cup (-\bar{\vartheta}, 0]$

Based on the obtained domain of $z$ under the four cases, we then can obtain the corresponding domains of $s$ as follows.

$$G_1 = \left[-\left(\bar{f} - \underline{f}\right), \bar{f} - \underline{f}\right] \times \{0\} \quad \text{(S1)}$$

$$G_2 = \left[-\left(\bar{f} - \underline{f}\right), \bar{f} - \underline{f}\right] \times \left(-\bar{\vartheta}, -2\bar{\vartheta}\right] \quad \text{(S2)}$$

$$G_3 = \left[-\left(\bar{f} - \underline{f}\right), \bar{f} - \underline{f}\right] \times \left(\bar{\vartheta}, 2\bar{\vartheta}\right] \quad \text{(S3)}$$

$$G_4 = \left[-(\overline{f}-\underline{f}), \overline{f}-\underline{f}\right] \times \left[0, \overline{\vartheta}\right] \cup \left(-\overline{\vartheta}, 0\right] \tag{S4}$$

Given $G_1$, $G_2$, $G_3$, and $G_4$, it can be easily deduced that Property 1 and Property 2 are satisfied, which completes the proof.

**Proof for Lemma 2**:

Given two arbitrary solutions $\mathbf{x}_1$ and $\mathbf{x}_2$ of a COP, we can have that

$$\begin{aligned}
&\pi(\mathbf{x}_1) - \pi(\mathbf{x}_2) \\
&= (f(\mathbf{x}_2) - f(\mathbf{x}_1)) + (\sigma(\mathbf{x}_2) - \sigma(\mathbf{x}_1)) \\
&= (f(\mathbf{x}_2) - f(\mathbf{x}_1)) + (\mathbb{I}(\rho(\mathbf{x}_2)>0) - \mathbb{I}(\rho(\mathbf{x}_1)>0)) \cdot (\overline{f}-\underline{f}) + \frac{(\overline{f}-\underline{f})}{\xi \cdot \overline{\vartheta}} (\theta(\mathbf{x}_2) - \theta(\mathbf{x}_1)). \\
&= -y + (\mathbb{I}(\rho(\mathbf{x}_2)>0) - \mathbb{I}(\rho(\mathbf{x}_1)>0)) \cdot (\overline{f}-\underline{f}) - \frac{(\overline{f}-\underline{f})}{\xi \cdot \overline{\vartheta}} z
\end{aligned} \tag{S5}$$

In what followings, let's discuss Equation (S5) under the five sub-conditions of Equation (31).

(1) If $-2\overline{\vartheta} \leq z < -\overline{\vartheta}$, as analyzed in Section 4.1.2 $\mathbf{x}_1 \in F$ and $\mathbf{x}_2 \notin F$ such that $\rho(\mathbf{x}_1) = 0$, $\rho(\mathbf{x}_2) > 0$, that is,

$$\mathbb{I}(\rho(\mathbf{x}_2)>0) - \mathbb{I}(\rho(\mathbf{x}_1)>0) = 1. \tag{S6}$$

As such,

$$\begin{aligned}
\pi(\mathbf{x}_1) - \pi(\mathbf{x}_2) &= -y + (\overline{f}-\underline{f}) - \frac{(\overline{f}-\underline{f})}{\xi \cdot \overline{\vartheta}} z \\
&> -y + (\overline{f}-\underline{f}) - \frac{(\overline{f}-\underline{f})}{\xi \cdot \overline{\vartheta}} (-\overline{\vartheta}). \\
&= -y + (\overline{f}-\underline{f}) + \frac{(\overline{f}-\underline{f})}{\xi} \\
&> 0
\end{aligned} \tag{S7}$$

(2) If $-\overline{\vartheta} \leq z < -\xi \cdot \overline{\vartheta}$, as analyzed in Section 4.1.2 $\mathbf{x}_1 \notin F$, $\mathbf{x}_2 \notin F$. Consequently, $\rho(\mathbf{x}_1) > 0$, $\rho(\mathbf{x}_2) > 0$, that is,

$$\mathbb{I}(\rho(\mathbf{x}_2)>0) - \mathbb{I}(\rho(\mathbf{x}_1)>0) = 0. \tag{S8}$$

As such,

$$\begin{aligned}
\pi(\mathbf{x}_1) - \pi(\mathbf{x}_2) &= -y - \frac{(\overline{f}-\underline{f})}{\xi \cdot \overline{\vartheta}} z \\
&> -y - \frac{(\overline{f}-\underline{f})}{\xi \cdot \overline{\vartheta}} (-\xi \cdot \overline{\vartheta}). \\
&= -y + (\overline{f}-\underline{f}) \\
&> 0
\end{aligned} \tag{S9}$$

(3) If $-\xi \cdot \overline{\vartheta} \le z < -\xi \cdot \eta \cdot \overline{\vartheta}$, as analyzed in Section 4.1.2 $\mathbf{x}_1 \notin F$, $\mathbf{x}_2 \notin F$, which indicates that Equation (S8) is also fulfilled in this situation. As such,

$$\begin{aligned} \pi(\mathbf{x}_1) - \pi(\mathbf{x}_2) &= -y - \frac{(\overline{f} - \underline{f})}{\xi \cdot \overline{\vartheta}} z \\ &> -y - \frac{(\overline{f} - \underline{f})}{\xi \cdot \overline{\vartheta}} (-\xi \cdot \eta \cdot \overline{\vartheta}) \\ &= -y + \eta \cdot (\overline{f} - \underline{f}) \end{aligned} \tag{S10}$$

Moreover, from $-(\overline{f} - \underline{f}) \le y \le \eta \cdot (\overline{f} - \underline{f})$ we can deduce that $-y \ge -\eta \cdot (\overline{f} - \underline{f})$, which demonstrates $\pi(\mathbf{x}_1) > \pi(\mathbf{x}_2)$.

(4) If $-\xi \cdot \eta \cdot \overline{\vartheta} \le z \le 0$, Equation (S8) holds and thus

$$\pi(\mathbf{x}_1) - \pi(\mathbf{x}_2) = \begin{cases} -y - \dfrac{(\overline{f} - \underline{f})}{\xi \cdot \overline{\vartheta}} z > -y, & \text{if } -\xi \cdot \eta \cdot \overline{\vartheta} \le z < 0 \\ -y, & \text{if } z = 0 \end{cases} \tag{S11}$$

Besides, it can be derived from $-(\overline{f} - \underline{f}) \le y \le 0$ that $-y \ge 0$, and thus

$$\pi(\mathbf{x}_1) - \pi(\mathbf{x}_2) \begin{cases} > 0, & \text{if } -\xi \cdot \eta \cdot \overline{\vartheta} \le z < 0 \\ = -y, & \text{if } z = 0 \end{cases} \tag{S12}$$

Specially, when $z = 0$, $\mathbf{x}_1 \in F$, $\mathbf{x}_2 \in F$, then $y \ne 0$ since we have $z + y \ne 0$, therefore we can also obtain $\pi(\mathbf{x}_1) > \pi(\mathbf{x}_2)$ in the case of $z = 0$.

(5) If $0 < z \le \xi \cdot \eta \cdot \overline{\vartheta}$, $\mathbf{x}_1 \notin F$, $\mathbf{x}_2 \notin F$, and we can have that

$$\begin{aligned} \pi(\mathbf{x}_1) - \pi(\mathbf{x}_2) &= -y - \frac{(\overline{f} - \underline{f})}{\xi \cdot \overline{\vartheta}} z \\ &\ge -y - \frac{(\overline{f} - \underline{f})}{\xi \cdot \overline{\vartheta}} \xi \cdot \eta \cdot \overline{\vartheta} \\ &= -y - \eta \cdot (\overline{f} - \underline{f}) \end{aligned} \tag{S13}$$

Since $-(\overline{f} - \underline{f}) \le y < -\eta \cdot (\overline{f} - \underline{f})$, we can further deduce that

$$\pi(\mathbf{x}_1) - \pi(\mathbf{x}_2) > \eta \cdot (\overline{f} - \underline{f}) - \eta \cdot (\overline{f} - \underline{f}) = 0. \tag{S14}$$

In summary, this completes the proof.

**Proof for Lemma 3**:

Given two arbitrary solutions $\mathbf{x}_1$ and $\mathbf{x}_2$ of a COP, according to Lemma 2, we have that $\pi(\mathbf{x}_1) > \pi(\mathbf{x}_2)$ if $\varphi < 0$. Let $\mathbf{x}_1' = \mathbf{x}_2$, $\mathbf{x}_2' = \mathbf{x}_1$, then

$$z' = \theta(\mathbf{x}_1') - \theta(\mathbf{x}_2') = \theta(\mathbf{x}_2) - \theta(\mathbf{x}_1) = -z, \tag{S15}$$

$$y' = f(\mathbf{x}_1') - f(\mathbf{x}_2') = f(\mathbf{x}_2) - f(\mathbf{x}_1) = -y. \tag{S16}$$

As such, $z'$ and $y'$ meet Equation (32) when $\varphi < 0$. Furthermore,

$$\pi(\mathbf{x}_1') - \pi(\mathbf{x}_2') = \pi(\mathbf{x}_2) - \pi(\mathbf{x}_1) = -z < 0, \tag{S17}$$

which completes the proof.

**Proof for Lemma 4**:

When $\varphi = 0$, as analyzed in Section 4.1.2 $\mathbf{x}_1 \notin F$, $\mathbf{x}_2 \notin F$, and thus

$$\pi(\mathbf{x}_1) - \pi(\mathbf{x}_2) = -y - \frac{(\overline{f} - \underline{f})}{\xi \cdot \overline{\vartheta}} z \tag{S18}$$

Besides, we have that $z < -e \cdot y$, therefore

$$\begin{aligned}
\pi(\mathbf{x}_1) - \pi(\mathbf{x}_2) &= -y - \frac{(\overline{f} - \underline{f})}{\xi \cdot \overline{\vartheta}} z \\
&> -y - \frac{(\overline{f} - \underline{f})}{\xi \cdot \overline{\vartheta}} \cdot (-e \cdot y) \\
&= -y + \frac{(\overline{f} - \underline{f})}{\xi \cdot \overline{\vartheta}} \cdot \left( \frac{\xi \cdot \overline{\vartheta}}{(\overline{f} - \underline{f})} \cdot y \right) \\
&= 0
\end{aligned} \tag{S19}$$

Which completes the proof.

**Proof for Lemma 5**:

Given two arbitrary solutions $\mathbf{x}_1$ and $\mathbf{x}_2$ of a COP, if $\varphi = 0$ and $z = -e \cdot y$, then either $\mathbf{x}_1 \in F$, $\mathbf{x}_2 \in F$ or $\mathbf{x}_1 \notin F$, $\mathbf{x}_2 \notin F$ according to the analysis in Section 4.1.2. Consequently, from Equation (S5) we can have that

$$\begin{aligned}
\pi(\mathbf{x}_1) - \pi(\mathbf{x}_2) &= -y + \left( \mathbb{I}(\rho(\mathbf{x}_2) > 0) - \mathbb{I}(\rho(\mathbf{x}_1) > 0) \right) \cdot (\overline{f} - \underline{f}) - \frac{(\overline{f} - \underline{f})}{\xi \cdot \overline{\vartheta}} z \\
&= -y - \frac{(\overline{f} - \underline{f})}{\xi \cdot \overline{\vartheta}} z \\
&= -y - \frac{(\overline{f} - \underline{f})}{\xi \cdot \overline{\vartheta}} \cdot \left( -\frac{\xi \cdot \overline{\vartheta}}{\overline{f} - \underline{f}} \cdot y \right) \\
&= -y + y \\
&= 0
\end{aligned} \tag{S20}$$

which completes the proof.

**Proof for Lemma 6**:

Continued with the proof of Lemma 4, when $z > -e \cdot y$, we have that

$$\begin{aligned}\pi(\mathbf{x}_1) - \pi(\mathbf{x}_2) &= -y - \frac{(\overline{f} - \underline{f})}{\xi \cdot \overline{\vartheta}} z \\ &< -y - \frac{(\overline{f} - \underline{f})}{\xi \cdot \overline{\vartheta}} \cdot (-e \cdot y) \\ &= -y + \frac{(\overline{f} - \underline{f})}{\xi \cdot \overline{\vartheta}} \cdot \left( \frac{\xi \cdot \overline{\vartheta}}{(\overline{f} - \underline{f})} \cdot y \right) \\ &= 0 \end{aligned} \quad (S21)$$

Which completes the proof.

**Proof for Lemma 7:**

If $0 < z \leq \xi \cdot \eta \cdot \overline{\vartheta}$, $\mathbf{x}_1 \notin F$, $\mathbf{x}_2 \notin F$, and we can have that

$$\pi(\mathbf{x}_1) - \pi(\mathbf{x}_2) = -y - \frac{f}{\xi \cdot \vartheta} z. \quad (S22)$$

Since $z < -\hat{e} \cdot y$, we can further deduce that

$$\pi(\mathbf{x}_1) - \pi(\mathbf{x}_2) > -y - \frac{f}{\xi \cdot \vartheta} \left( -\frac{\xi \cdot \vartheta}{f} \right) y = 0, \quad (S23)$$

Which completes the proof.

**Proof for Lemma 8:**

If $-\xi \cdot \eta \cdot \overline{\vartheta} \leq z < 0$, $\mathbf{x}_1 \notin F$, $\mathbf{x}_2 \notin F$, and we can have that

$$\pi(\mathbf{x}_1) - \pi(\mathbf{x}_2) = -y - \frac{f}{\xi \cdot \vartheta} z. \quad (S24)$$

Since $z > -\hat{e} \cdot y$, we can further deduce that

$$\pi(\mathbf{x}_1) - \pi(\mathbf{x}_2) < -y - \frac{f}{\xi \cdot \vartheta} \left( -\frac{\xi \cdot \vartheta}{f} \right) y = 0, \quad (S25)$$

Which completes the proof.

**Proof for Lemma 9:**

Let's first discuss the case of Inequality (42).(1). When $\varphi = 0$ and $z > 0$, $\mathbf{x}_1 \notin F$, $\mathbf{x}_2 \notin F$, and we can have that

$$\pi(\mathbf{x}_1) - \pi(\mathbf{x}_2) = -y - \frac{f}{\xi \cdot \vartheta} z. \quad (S26)$$

If $z < -\hat{e} \cdot y$ and $z > 0$, we can have that $\pi(\mathbf{x}_1) - \pi(\mathbf{x}_2) > 0$. If $z > -e \cdot y$ and $z > 0$, we can have that

$$\pi(\mathbf{x}_1) - \pi(\mathbf{x}_2) = -y - \frac{f}{\xi \cdot \vartheta} z$$
$$< y - \frac{f}{\xi \cdot \vartheta}(\xi \cdot \eta \cdot \overline{\vartheta}). \quad \text{(S27)}$$
$$< y - \eta(\overline{f} - \underline{f})$$

Besides, $-(\overline{f} - \underline{f}) \leq y < -\eta \cdot (\overline{f} - \underline{f})$, and thus $\pi(\mathbf{x}_1) - \pi(\mathbf{x}_2) < 0$

The discussion on the case of Inequality (42).(2) is similar to the case of Inequality (42).(1). This completes the proof.

### III. Extension to multi-objective COP

In general, a multi-objective COP can be formalized as follows.

$$\min \mathbf{F}(\mathbf{x}) = (f_1(\mathbf{x}), f_2(\mathbf{x}), \ldots f_k(\mathbf{x}))^T$$
$$s.t. \begin{cases} g_i(\mathbf{x}) \leq 0, i = 1,\ldots,m \\ h_i(\mathbf{x}) = 0, i = m+1,\ldots,n \end{cases} \quad \text{(S28)}$$

Given any two solutions $\mathbf{x}_1$ and $\mathbf{x}_2$, a $k+1$ dimensional Cartesian coordinate system is constructed by taking

$$y_1 = f_1(\mathbf{x}_1) - f_1(\mathbf{x}_2),\ y_2 = f_2(\mathbf{x}_1) - f_1(\mathbf{x}_2),\ \ldots,\ y_k = f_k(\mathbf{x}_1) - f_k(\mathbf{x}_2),\ z = \theta(\mathbf{x}_1) - \theta(\mathbf{x}_2) \quad \text{(S29)}$$

as its axes. Based on Section 4.1.1, we further obtain the difference division and rank related to $y_1, y_2, \ldots, y_k$, and $z$, respectively,

$$D_{\alpha_1}^{y_1},\ R_{\alpha_1}^{y_1},\ D_{\alpha_2}^{y_2},\ R_{\alpha_2}^{y_2},\ \ldots,\ D_{\alpha_k}^{y_k},\ R_{\alpha_k}^{y_k},\ D_{\beta}^{z},\ R_{\beta}^{z}. \quad \text{(S30)}$$

Subsequently, we define the Cartesian product of $y_1, y_2, \ldots, y_k$, and $z$ as follows,
$$s = y_1 \times y_2 \times \cdots y_k \times z. \quad \text{(S31)}$$
We can further obtain the difference division and rank associated with $s$ below.

$$D_{\gamma}^{s} = \left\{ d_{\varphi}^{s} = \bigcup_{\chi_1 + \chi_2 + \cdots + \chi_k + \lambda = \varphi} d_{\chi_1}^{y_1} \times d_{\chi_2}^{y_2} \times \cdots \times d_{\chi_k}^{y_k} \times d_{\lambda}^{z} \middle| d_{\chi_1}^{y_1} \in D_{\alpha_1}^{y_1}, d_{\chi_2}^{y_2} \in D_{\alpha_2}^{y_2}, \cdots d_{\chi_k}^{y_k} \in D_{\alpha_k}^{y_k}, d_{\lambda}^{z} \in D_{\beta}^{z} \right\} \quad \text{(S32)}$$

$$\varphi \in R_{\gamma}^{s} = \left\{ \chi_1 + \chi_2 + \cdots + \chi_k + \lambda \middle| \chi_1 \in R_{\alpha_1}^{y_1}, \chi_2 \in R_{\alpha_2}^{y_2}, \cdots, \chi_k \in R_{\alpha_k}^{y_k}, \lambda \in R_{\beta} \right\}$$
$$= \left\{ -\left(k + 1 + \beta + \sum_{l=1}^{k} \alpha_l\right), \cdots, \left(k + 1 + \beta + \sum_{l=1}^{k} \alpha_l\right) \right\} \quad \text{(S33)}$$

In this way, the following comparison criterion similar to the one denoted as Equation (28) can be obtained as

$$\begin{cases} \pi(\mathbf{x}_1) > \pi(\mathbf{x}_2), & \forall \mathbf{x}_1 \in S, \mathbf{x}_2 \in S, \text{ if } \varphi < 0 \\ \pi(\mathbf{x}_1) < \pi(\mathbf{x}_2), & \forall \mathbf{x}_1 \in S, \mathbf{x}_2 \in S, \text{ if } \varphi > 0 \\ \pi(\mathbf{x}_1) = \pi(\mathbf{x}_2), & \forall \mathbf{x}_1 \in S, \mathbf{x}_2 \in S, \text{ if } \varphi = 0 \text{ and } z = -\sum_{l=1}^{k}(e_l \cdot y_l) \\ \pi(\mathbf{x}_1) > \pi(\mathbf{x}_2), & \forall \mathbf{x}_1 \in S, \mathbf{x}_2 \in S, \text{ if } \varphi = 0 \text{ and } z < -\sum_{l=1}^{k}(e_l \cdot y_l) \\ \pi(\mathbf{x}_1) < \pi(\mathbf{x}_2), & \forall \mathbf{x}_1 \in S, \mathbf{x}_2 \in S, \text{ if } \varphi = 0 \text{ and } z > -\sum_{l=1}^{k}(e_l \cdot y_l) \end{cases} \quad (S34)$$

where

$$e_l = \begin{cases} \dfrac{2\bar{\vartheta}}{q_{\alpha_l}^{y_l}}, & \alpha_l > \beta \\ \dfrac{2\bar{\vartheta}}{\bar{f}-\underline{f}}, & \alpha_l = \beta, l = 1,2,\ldots k \\ \dfrac{q_{\beta}^{z}}{\bar{f}-\underline{f}}, & \alpha_l < \beta \end{cases} \quad (S35)$$

**IV. Boxplots of the minimum required FES for CHT1-CHT4**

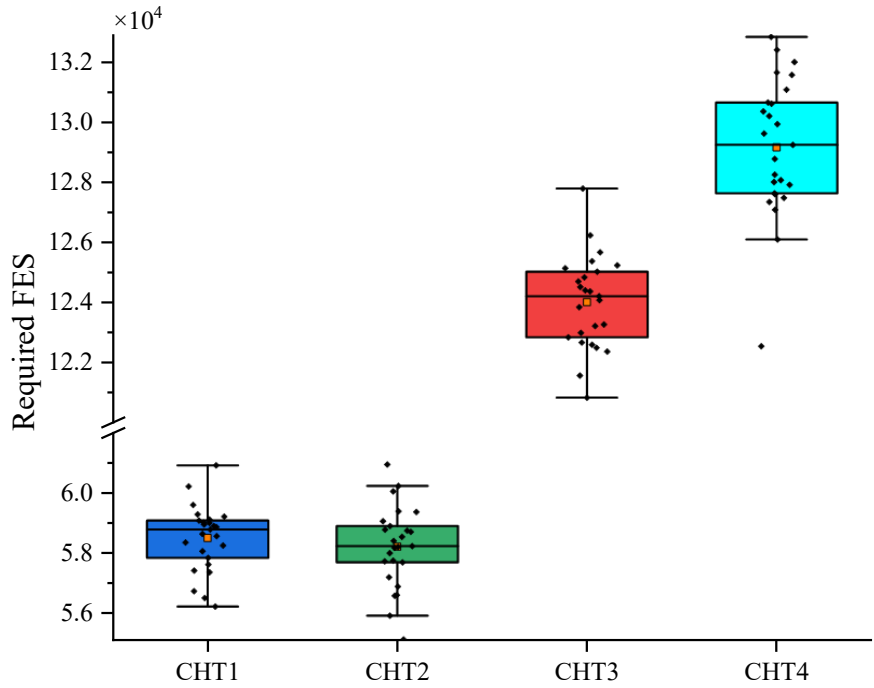

Fig. S1. Boxplot of the minimum required FES on G01

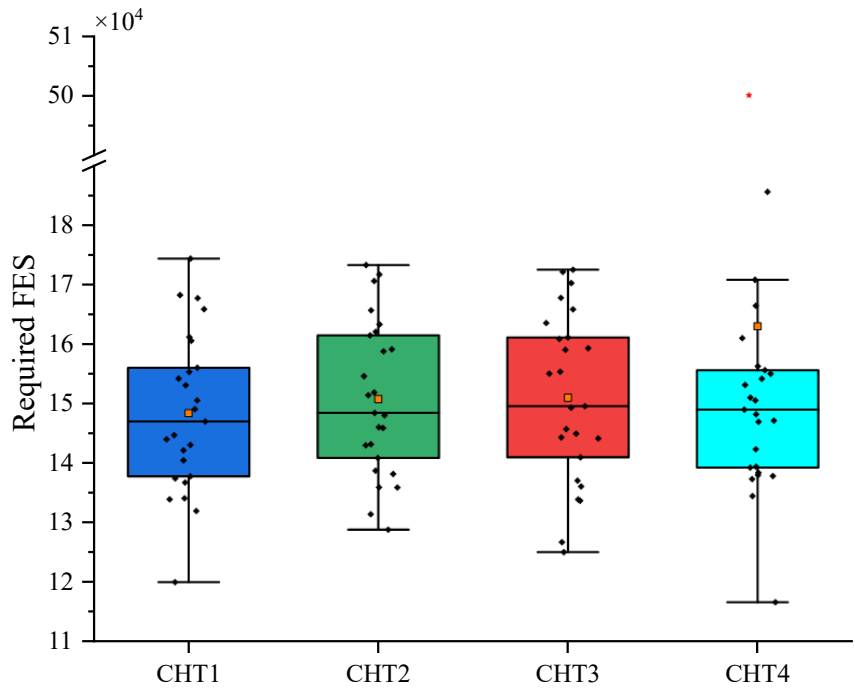

Fig. S2. Boxplot of the minimum required FES on G02

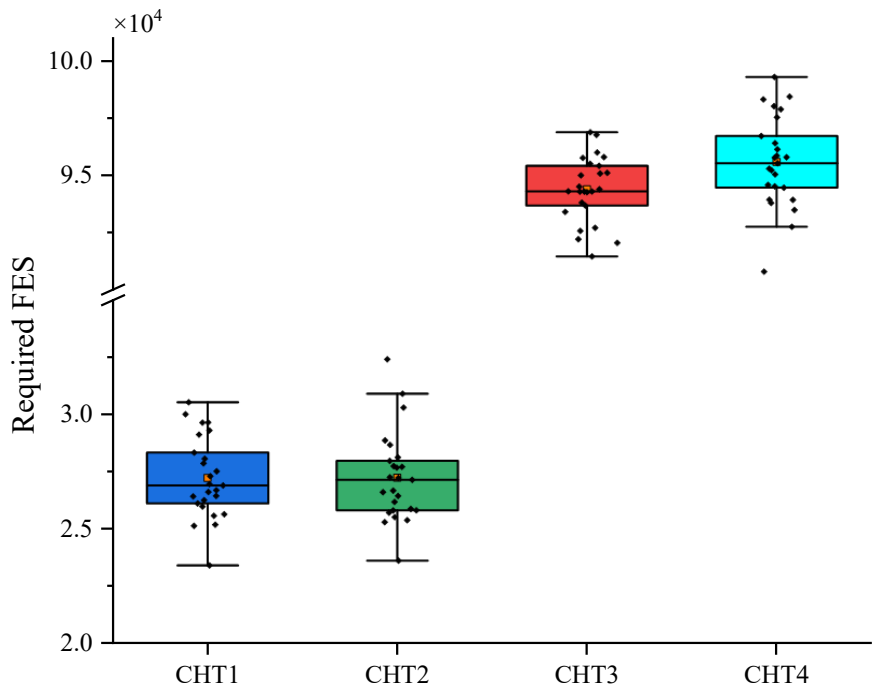

Fig. S3. Boxplot of the minimum required FES on G03

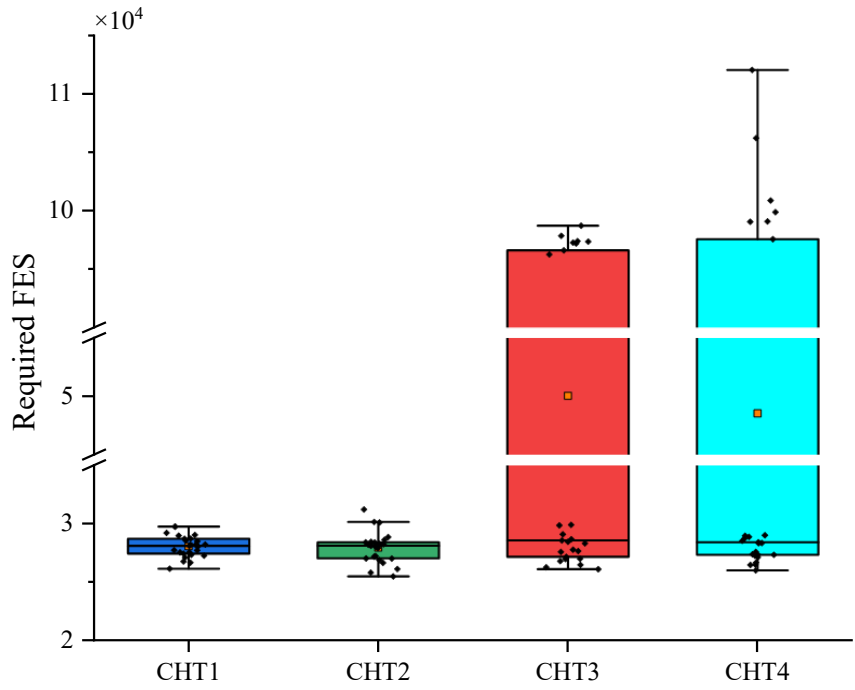

Fig. S4. Boxplot of the minimum required FES on G04

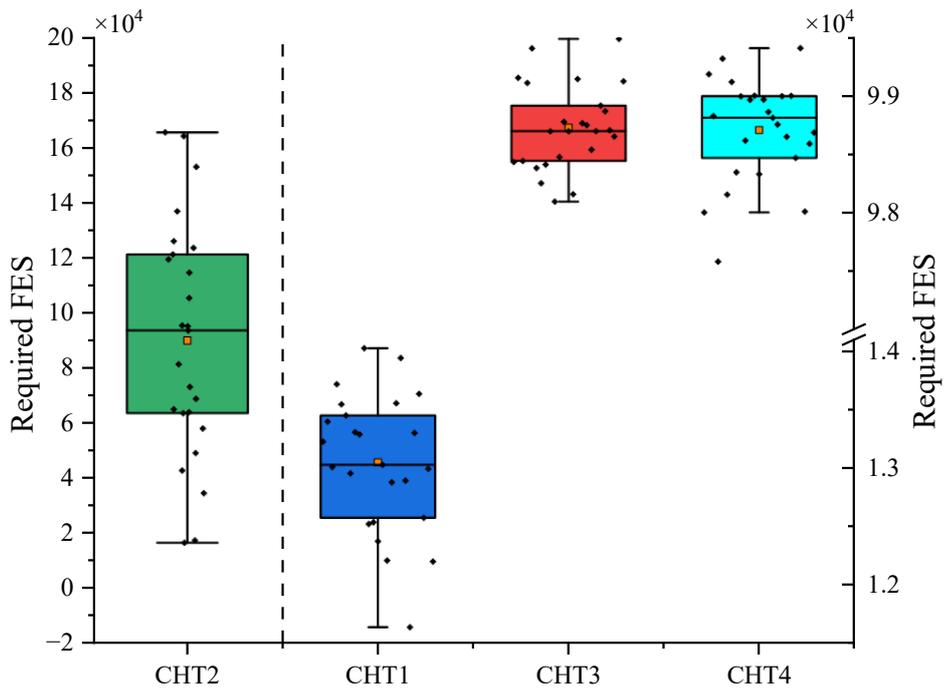

Fig. S5. Boxplot of the minimum required FES on G05

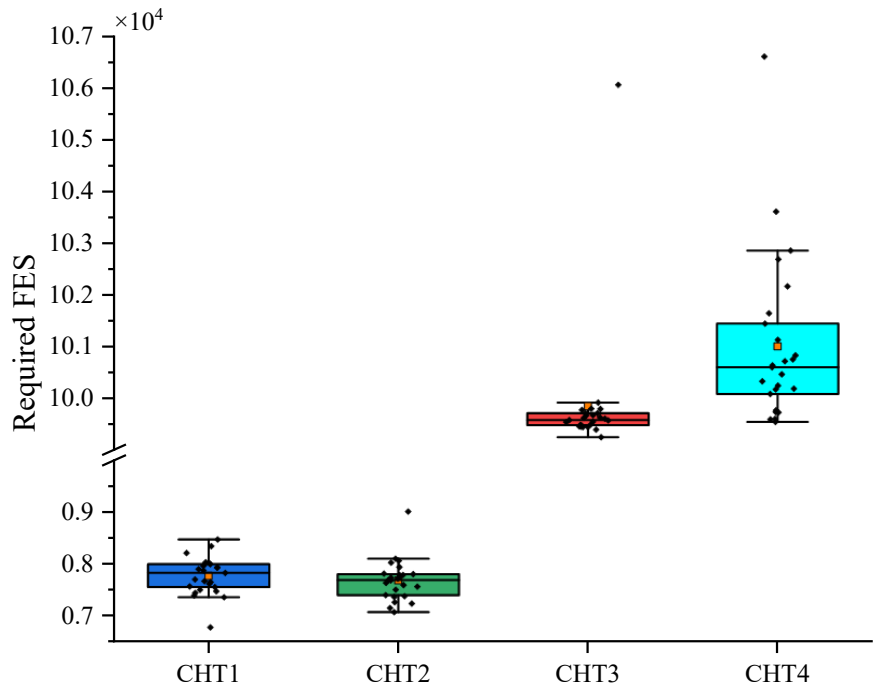

Fig. S6. Boxplot of the minimum required FES on G06

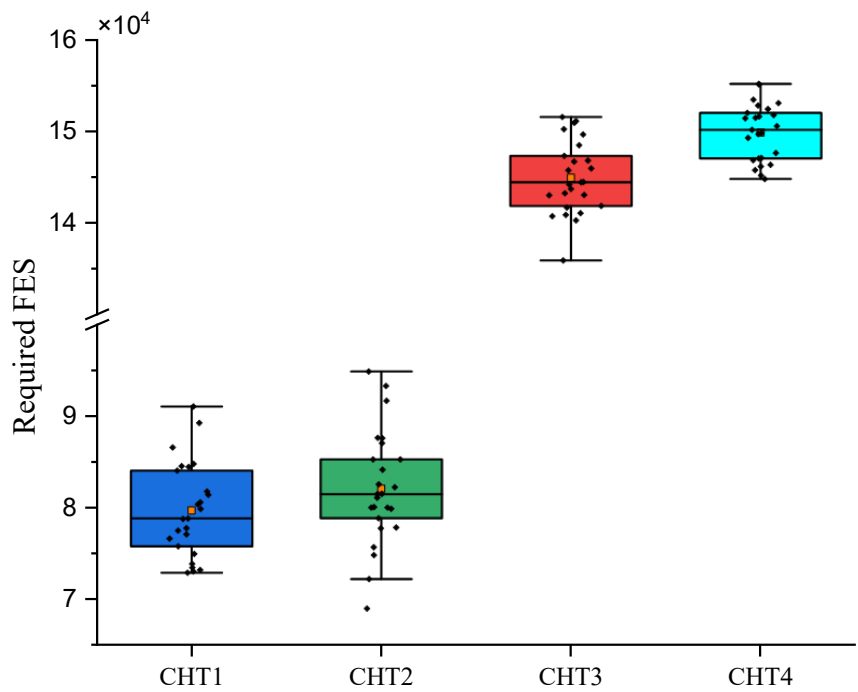

Fig. S7. Boxplot of the minimum required FES on G07

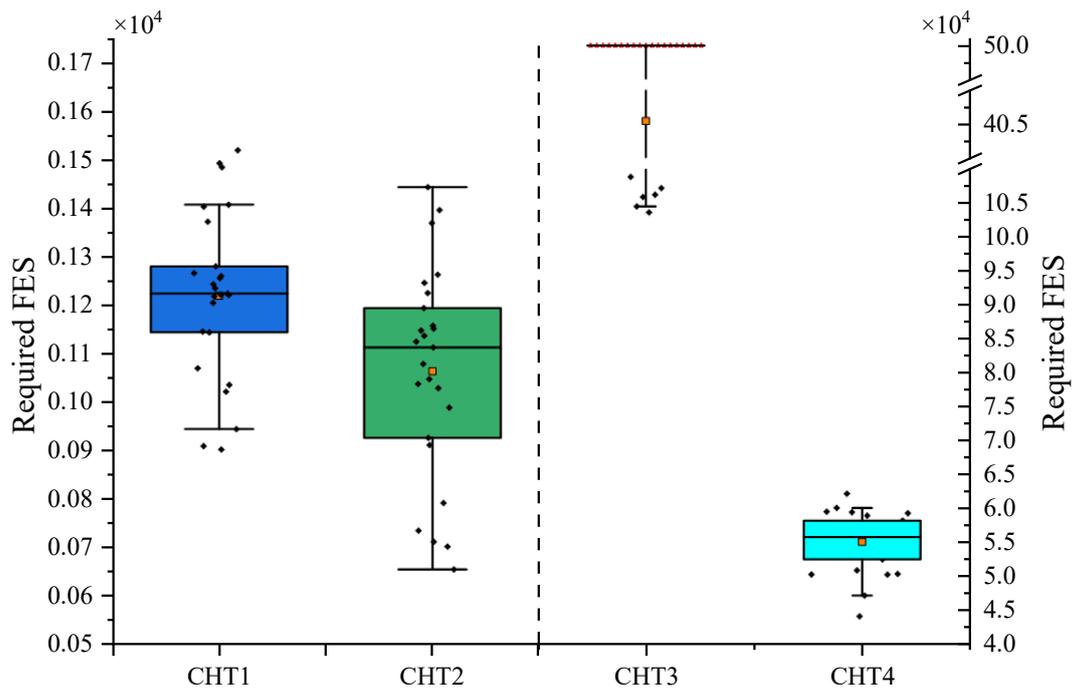

Fig. S8. Boxplot of the minimum required FES on G08

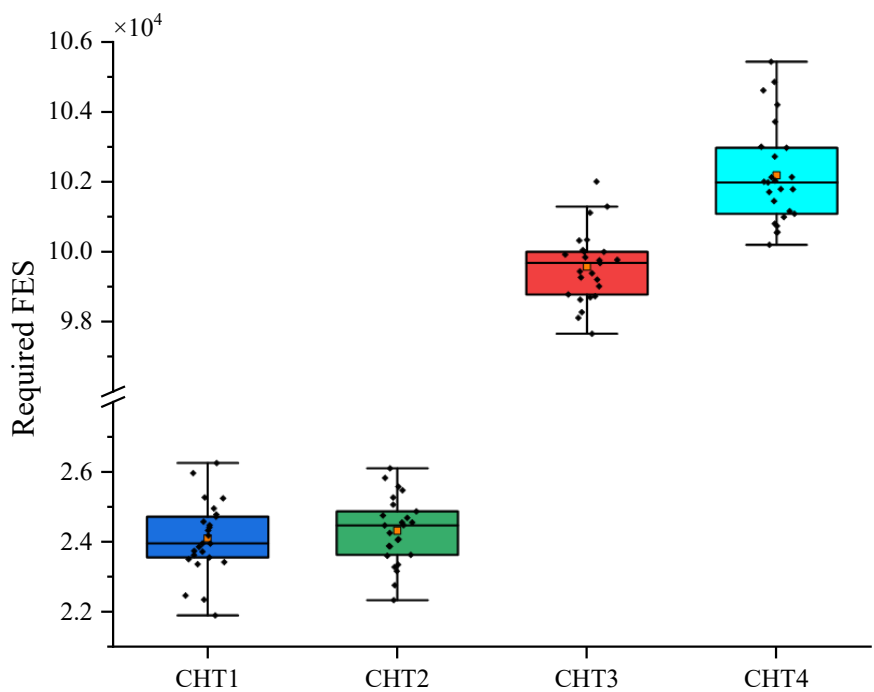

Fig. S9. Boxplot of the minimum required FES on G09

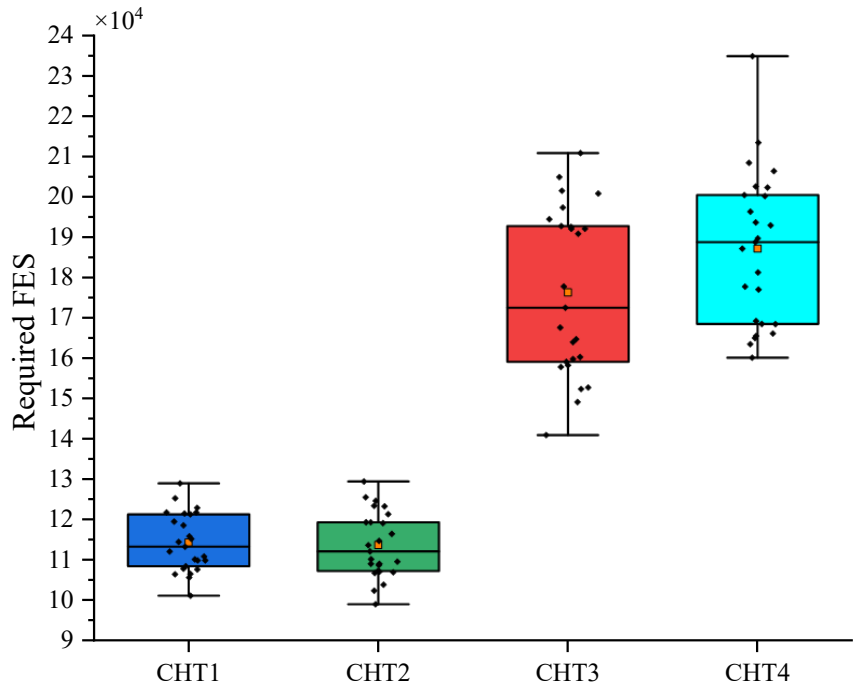

Fig. S10. Boxplot of the minimum required FES on G10

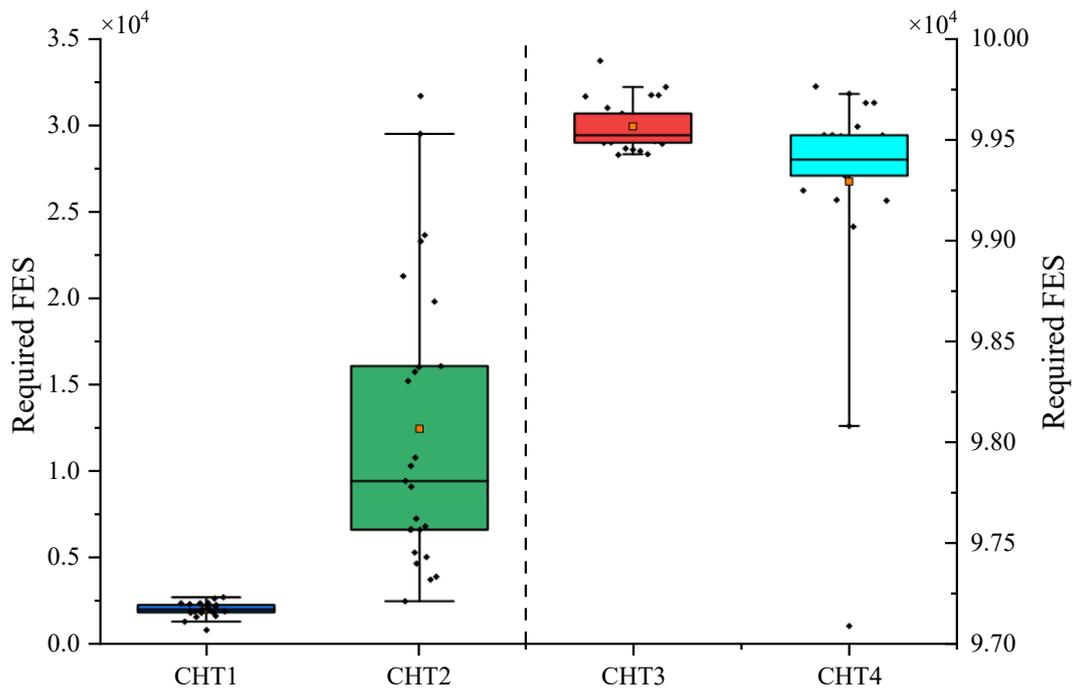

Fig. S11. Boxplot of the minimum required FES on G11

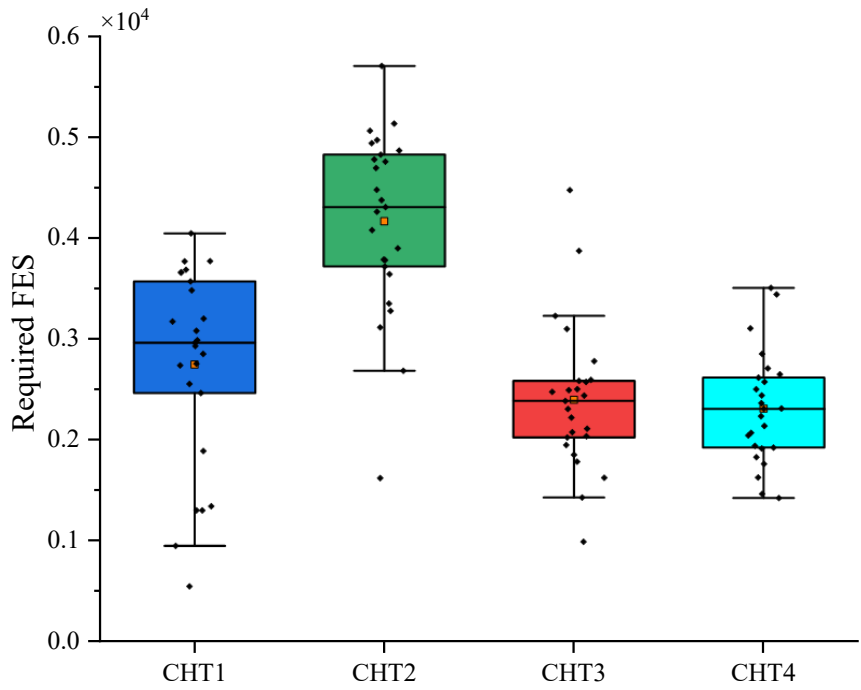

Fig. S12. Boxplot of the minimum required FES on G12

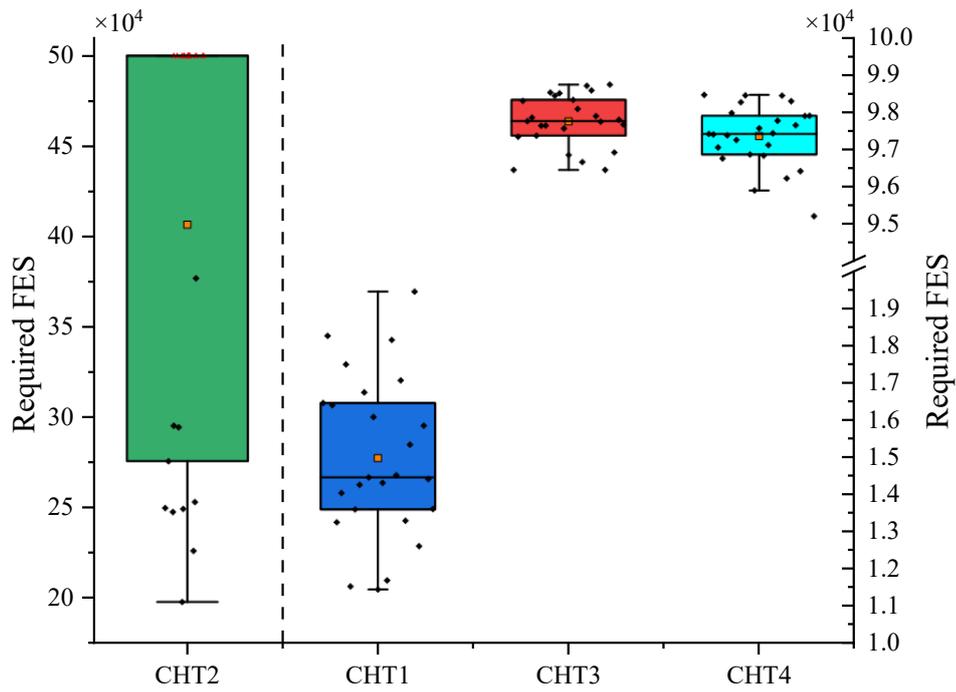

Fig. S13. Boxplot of the minimum required FES on G13

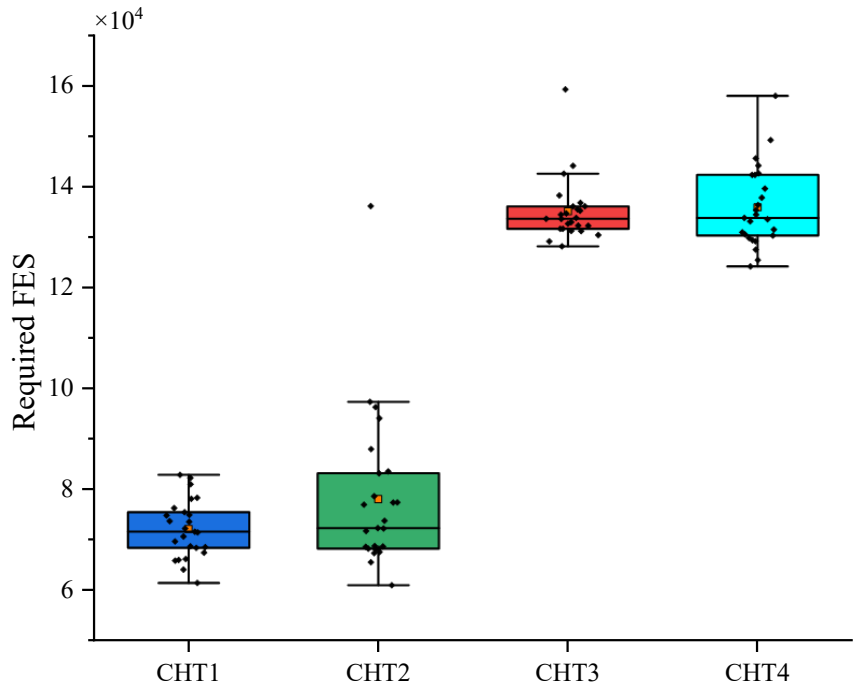

Fig. S14. Boxplot of the minimum required FES on G14

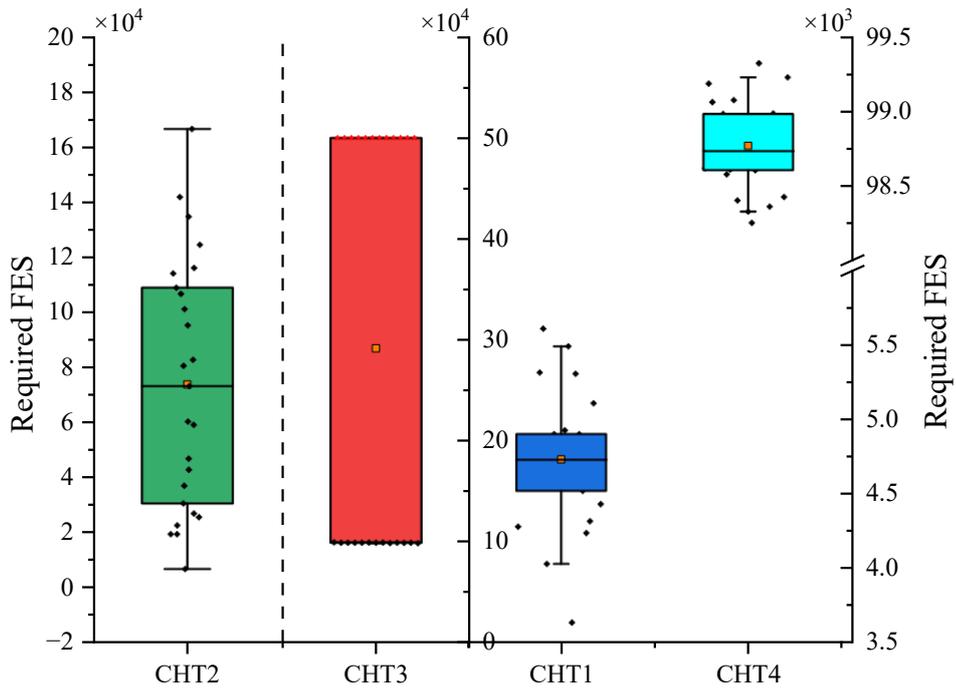

Fig. S15. Boxplot of the minimum required FES on G15

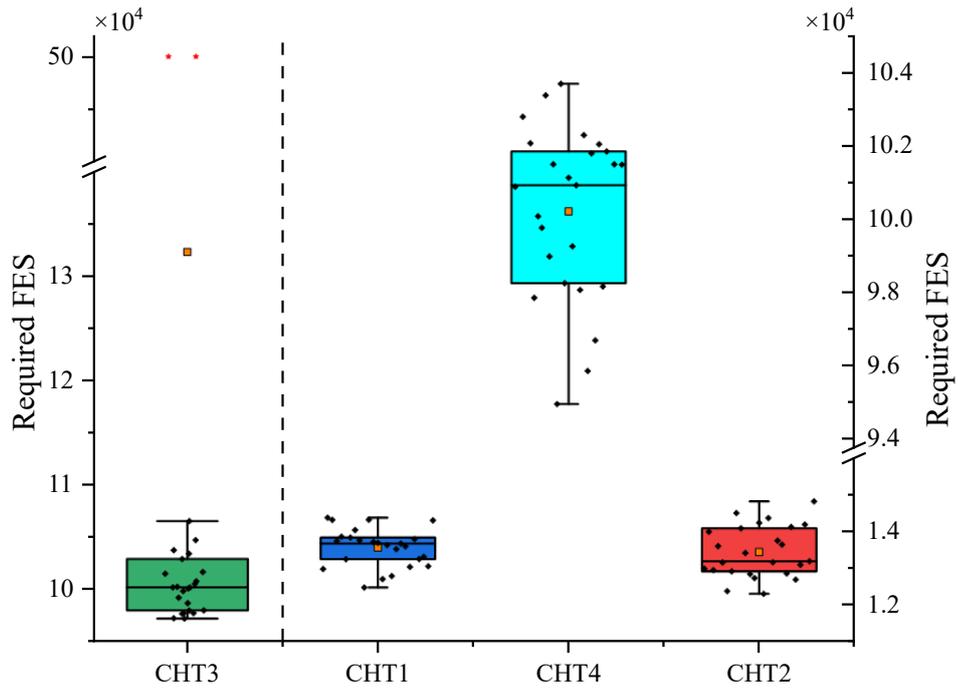

Fig. S16. Boxplot of the minimum required FES on G16

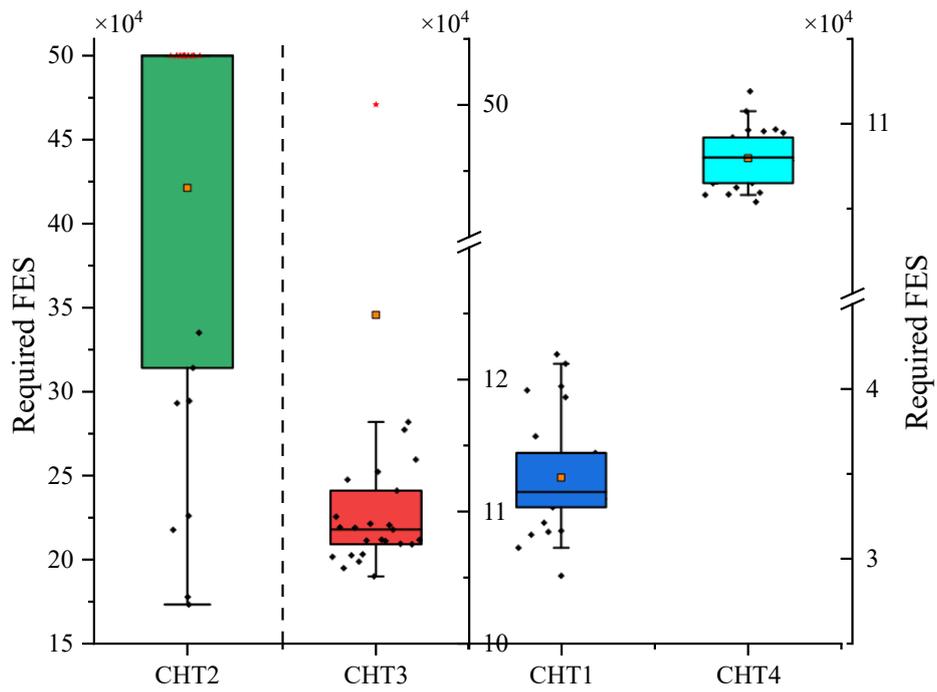

Fig. S17. Boxplot of the minimum required FES on G17

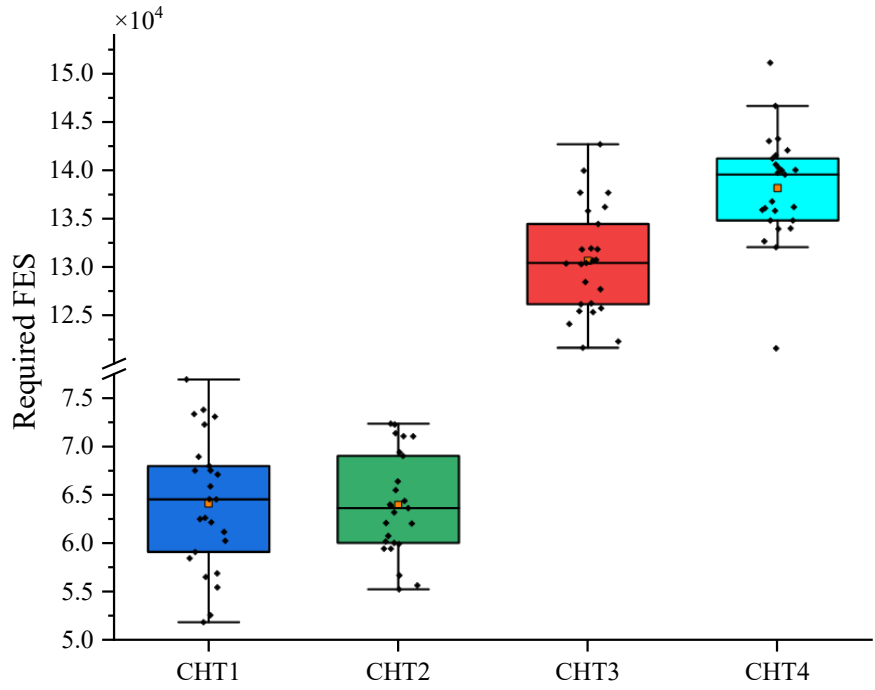

Fig. S18. Boxplot of the minimum required FES on G18

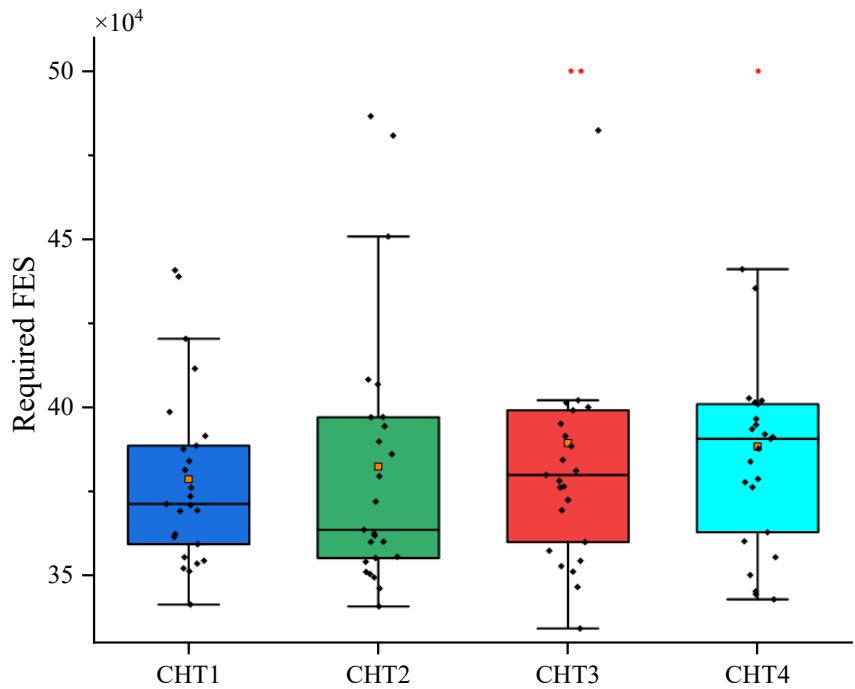

Fig. S19. Boxplot of the minimum required FES on G19

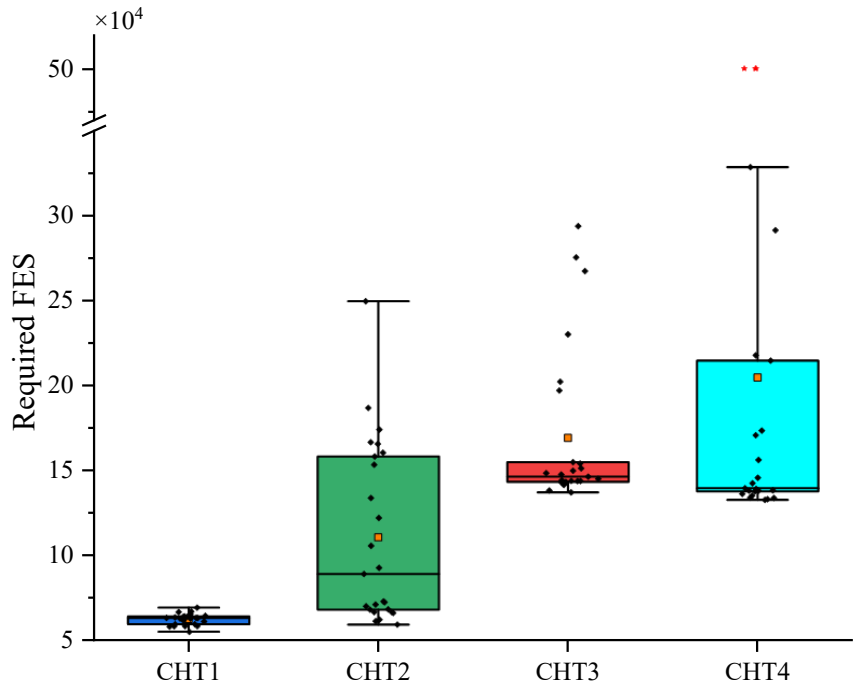

Fig. S21. Boxplot of the minimum required FES on G21

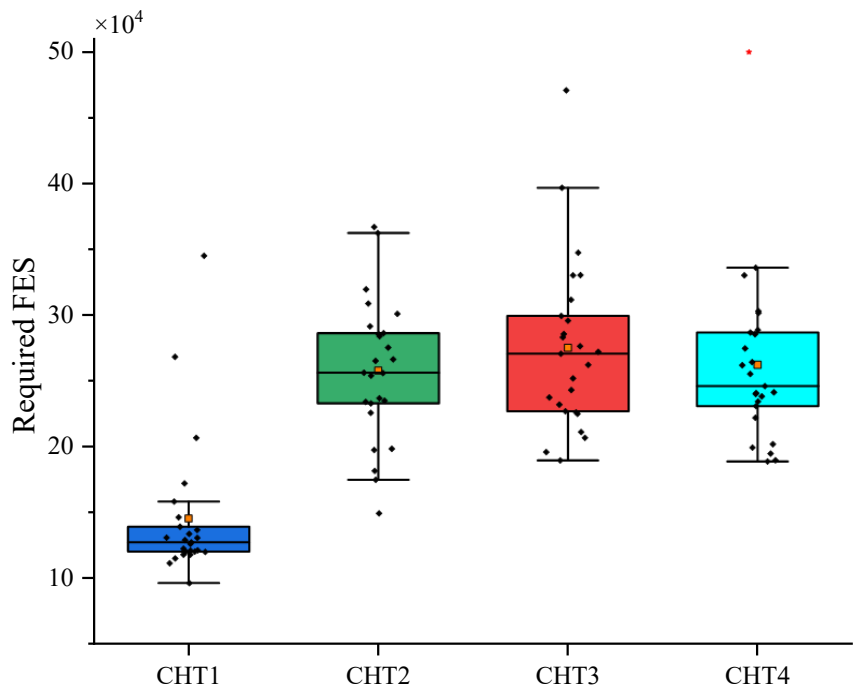

Fig. S23. Boxplot of the minimum required FES on G23

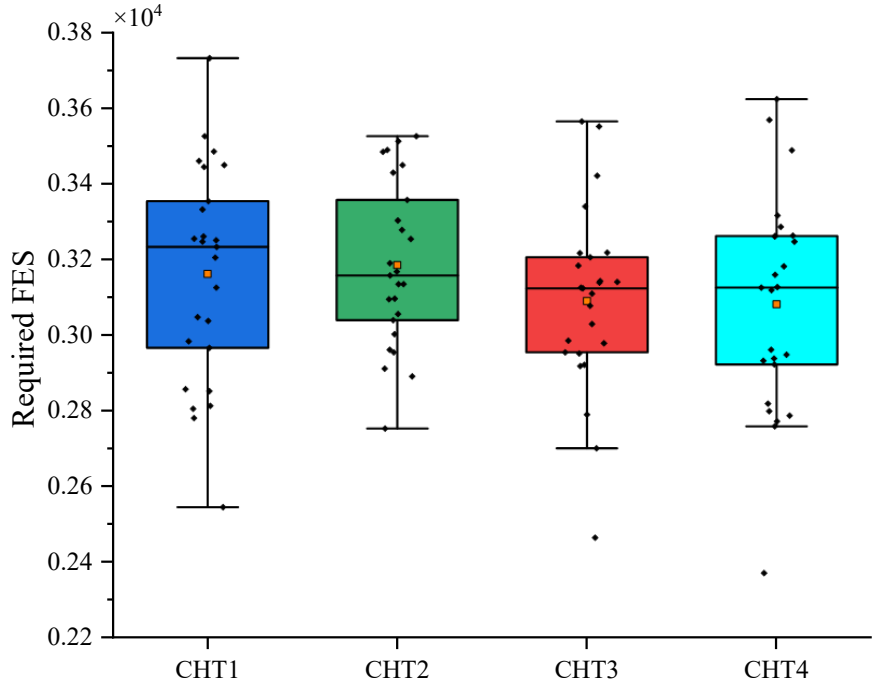

Fig. S24. Boxplot of the minimum required FES on G24

## V. Statistical results of A1-A4

Table S1. Statistical results of A1-A4 on C01-C07 with $D = 10$

|    |        | C01 | C02 | C03 | C04 | C05 | C06 | C07 |
|----|--------|-----|-----|-----|-----|-----|-----|-----|
| A1 | FR     | 100% | 100% | 100% | 100% | 100% | 100% | 100% |
|    | Best   | 0 | 0 | 0 | 1.357e+1 | 0 | 0 | -4.744e+2 |
|    | Median | 0 | 0 | 0 | 1.592e+1 | 0 | 0 | -4.744e+2 |
|    | Worst  | 0 | 0 | 0 | 1.691e+1 | 0 | 0 | -4.613e+2 |
|    | Mean   | 0 | 0 | 0 | 1.568e+1 | 0 | 0 | -4.719e+2 |
|    | STD    | 0 | 0 | 0 | 8.012e-1 | 0 | 0 | 5.080 |
| A2 | FR     | 100% | 100% | 12% | 100% | 100% | 0% | 0% |
|    | Best   | 0 | 0 | 1.051e+5 | 1.357e+1 | 0 | 3.377e+2/1.931 | -1.015e+2/5.531e-3 |
|    | Median | 0 | 0 | 1.160e+5 | 1.357e+1 | 0 | 1.356e+3/3.004 | 2.986e+1/1.993e-2 |
|    | Worst  | 0 | 0 | 2.378e+5 | 1.592e+1 | 0 | 1.997e+3/1.857 | 2.062e+2/1.343e-2 |
|    | Mean   | 0 | 0 | 1.530e+5 | 1.457e+1 | 0 | 1.205e+3/2.514 | 4.065e+1/1.004e-2 |
|    | STD    | 0 | 0 | 6.017e+4 | 1.132 | 0 | 4.464e+2/0.716 | 7.486e+1/1.147e-2 |
| A3 | FR     | 100% | 100% | 100% | 100% | 100% | 100% | 100% |
|    | Best   | 0 | 0 | 0 | 1.357e+1 | 0 | 0 | -4.744e+2 |
|    | Median | 0 | 0 | 0 | 1.357e+1 | 0 | 0 | -4.684e+2 |
|    | Worst  | 0 | 0 | 0 | 1.592e+1 | 0 | 0 | -4.660e+2 |
|    | Mean   | 0 | 0 | 0 | 1.470e+1 | 0 | 0 | -4.700e+2 |
|    | STD    | 0 | 0 | 0 | 1.172 | 0 | 0 | 2.846 |
| A4 | FR     | 100% | 100% | 92% | 100% | 100% | 100% | 88% |
|    | Best   | 0 | 0 | 3.534e+3 | 1.357e+1 | 0 | 3.490e+2 | -1.012e+2 |
|    | Median | 0 | 0 | 2.114e+4 | 1.359e+1 | 0 | 1.260e+3 | 1.278e+1 |
|    | Worst  | 0 | 0 | 1.180e+5 | 1.380e+1 | 0 | 1.369e+3 | 1.056e+2 |
|    | Mean   | 0 | 0 | 3.155e+4 | 1.361e+1 | 0 | 6.489e+2 | 3.743 |
|    | STD    | 0 | 0 | 3.7012e+4 | 6.15e-2 | 0 | 2.837e+2 | 6.957e+1 |

Table S2. Statistical results of A1-A4 on C08-C14 with $D = 10$

|    |      | C08 | C09 | C10 | C11 | C12 | C13 | C14 |
|----|------|-----|-----|-----|-----|-----|-----|-----|
| A1 | FR   | 100% | 100% | 100% | 76% | 100% | 100% | 100% |
|    | Best | -1.348e-3 | -4.975e-3 | -5.096e-4 | -1.688e-1 | 3.988 | 0 | 2.376 |
|    | Median | -1.348e-3 | -4.975e-3 | -5.096e-4 | -1.688e-1 | 3.988 | 0 | 2.376 |
|    | Worst  | -1.348e-3 | -4.975e-3 | -5.096e-4 | -9.219e-2 | 3.988 | 0 | 2.376 |
|    | Mean   | -1.348e-3 | -4.975e-3 | -5.096e-4 | -1.613e-1 | 3.988 | 0 | 2.376 |
|    | STD    | 0 | 0 | 0 | 2.200e-2 | 0 | 0 | 0 |
| A2 | FR     | 100% | 100% | 100% | 100% | 100% | 100% | 100% |
|    | Best   | -1.348e-3 | -4.975e-3 | -5.096e-4 | -1.688e-1 | 3.988 | 0 | 3.274 |
|    | Median | -1.348e-3 | -4.975e-3 | -5.096e-4 | -1.688e-1 | 3.988 | 0 | 3.553 |
|    | Worst  | -1.348e-3 | -4.975e-3 | -5.096e-4 | -1.688e-1 | 3.988 | 0 | 3.719 |
|    | Mean   | -1.348e-3 | -4.975e-3 | -5.096e-4 | -1.688e-1 | 3.988 | 0 | 3.545 |
|    | STD    | 0 | 0 | 0 | 0 | 0 | 0 | 1.173e-1 |
| A3 | FR     | 100% | 96% | 100% | 12% | 100% | 100% | 100% |
|    | Best   | -1.272e-3 | -4.975e-3 | -4.892e-4 | -1.687e-1 | 3.988 | 0 | 2.376 |
|    | Median | -1.042e-3 | -4.969e-3 | -4.180e-4 | -1.685e-1 | 3.988 | 0 | 2.376 |
|    | Worst  | -6.560e-4 | 3.513e-1 | 4.718e-5 | -1.681e-1 | 3.988 | 0 | 2.833 |
|    | Mean   | -1.021e-3 | 1.004e-2 | -3.848e-4 | -1.684e-1 | 3.988 | 0 | 2.415 |
|    | STD    | 1.902e-4 | 7.117e-2 | 1.139e-4 | 2.684e-4 | 0 | 0 | 1.149e-1 |
| A4 | FR     | 100% | 100% | 100% | 100% | 100% | 100% | 100% |
|    | Best   | -1.348e-3 | -4.975e-3 | -5.096e-4 | -1.688e-1 | 3.988 | 0 | 2.498 |
|    | Median | -1.348e-3 | -4.975e-3 | -5.096e-4 | -1.688e-1 | 3.993 | 0 | 2.905 |
|    | Worst  | -1.348e-3 | -4.855e-3 | -5.096e-4 | -1.688e-1 | 3.998 | 0 | 3.245 |
|    | Mean   | -1.348e-3 | -4.970e-3 | -5.096e-4 | -1.688e-1 | 3.993 | 0 | 2.882 |
|    | STD    | 0 | 2.44e-5 | 0 | 0 | 0.232e-2 | 0 | 0.204 |

Table S3. Statistical results of A1-A4 on C15-C21 with $D = 10$

|    |        | C15 | C16 | C17 | C18 | C19 | C20 | C21 |
|----|--------|-----|-----|-----|-----|-----|-----|-----|
| A1 | FR     | 100% | 100% | 0% | 100% | 0% | 100% | 100% |
|    | Best   | 2.356 | 0 | 1.241e-2/9.023 | 3.660e+1 | -9.252/1.329e+4 | 2.712e-1 | 3.988 |
|    | Median | 2.356 | 0 | 1.919e-2/9.013 | 3.660e+1 | -9.252/1.329e+4 | 4.326e-1 | 3.988 |
|    | Worst  | 5.498 | 1.414e+1 | 2.191e-1/1.139e+1 | 1.413e+2 | -9.032/1.329e+4 | 5.894e-1 | 12.032 |
|    | Mean   | 2.482 | 1.508 | 7.815e-2/1.072e+1 | 4.078e+1 | -9.243/1.329e+4 | 4.468e-1 | 5.183 |
|    | STD    | 6.156e-1 | 3.093 | 8.320e-2/1.523 | 2.051e+1 | 4.301e-2/8.979e-3 | 6.752e-2 | 2.396 |
| A2 | FR     | 84% | 92% | 0% | 100% | 0% | 100% | 100% |
|    | Best   | 1.178e+1 | 3.770e+1 | 9.929e-1/1.100e+1 | 3.660e+1 | 0/1.327e+4 | 5.885e-1 | 3.988 |
|    | Median | 1.806e+1 | 6.283e+1 | 1.021/1.100e+1 | 3.660e+1 | 0/1.327e+4 | 7.386e-1 | 4.544 |
|    | Worst  | 2.121e+1 | 8.796e+1 | 1.160/1.100e+1 | 3.722e+1 | 0/1.327e+4 | 8.438e-1 | 12.369 |
|    | Mean   | 1.776e+1 | 6.160e+1 | 1.036/1.100e+1 | 3.662e+1 | 0/1.327e+4 | 7.173e-1 | 6.139 |
|    | STD    | 2.548 | 1.150e+1 | 3.914e-2/2.475e-4 | 1.214e+1 | 0/0 | 7.220e-2 | 2.955 |
| A3 | FR     | 100% | 100% | 0% | 72% | 0% | 100% | 100% |
|    | Best   | 2.356 | 0 | 1.034e-2/1.100e+1 | 3.660e+1 | 0/1.327e+4 | 5.584e-1 | 3.988 |
|    | Median | 5.498 | 0 | 4.226e-2/9.000 | 3.660e+1 | 0/1.327e+4 | 6.886e-1 | 3.988 |
|    | Worst  | 1.806e+1 | 0 | 5.510e-1/9.000 | 3.662e+1 | 0/1.327e+4 | 8.094e-1 | 2.298e+1 |
|    | Mean   | 6.629 | 0 | 9.633e-2/9.080 | 3.660e+1 | 0/1.327e+4 | 6.942e-1 | 6.567 |
|    | STD    | 4.857 | 0 | 1.527e-1/3.919e-1 | 4.669e-3 | 0/0 | 7.171e-2 | 5.551 |
| A4 | FR     | 48% | 100% | 0% | 0% | 0% | 100% | 100% |
|    | Best   | 1.178e+1 | 3.142e+1 | 0.835/- | 5.355e+2/- | 7.58e-7/- | 0.839e-1 | 3.988 |
|    | Median | 1.178e+1 | 3.927e+1 | 0.991/- | 2.926e+3/- | 1.47e-6/- | 0.191 | 3.988 |
|    | Worst  | 1.492e+1 | 5.655e+1 | 1.011/- | 3.326e+3/- | 1.87e-6/- | 0.193 | 3.990 |
|    | Mean   | 1.455e+1 | 4.072e+1 | 0.911/- | 2.106e+3/- | 1.45e-6/- | 0.284 | 3.988 |
|    | STD    | 3.774 | 6.725 | 0.177/- | 5.905e+3/- | 3.12e-7/- | 0.579e-1 | 0.520e-3 |

Table S4. Statistical results of A1-A4 on C22-C28 with $D = 10$

|    |        | C22 | C23 | C24 | C25 | C26 | C27 | C28 |
|----|--------|-----|-----|-----|-----|-----|-----|-----|
| A1 | FR     | 100% | 100% | 68% | 100% | 0% | 100% | 0% |
|    | Best   | 3.937e-27 | 2.376 | 5.498 | 0 | 1.827e-2/11.009 | 3.660e+1 | 2.654e+1/1.335e+4 |
|    | Median | 3.937e-27 | 2.376 | 8.639 | 1.571 | 1.843e-1/12.513 | 3.660e+1 | 3.721e+1/1.335e+4 |
|    | Worst  | 3.987 | 2.376 | 18.064 | 4.398 | 2.403e-1/11.066 | 3.660e+1 | 4.633e+1/1.336e+4 |
|    | Mean   | 7.973e-1 | 2.376 | 9.933 | 8.796 | 1.318e-1/11.172 | 3.660e+1 | 3.777e+1/1.335e+4 |
|    | STD    | 1.595 | 0 | 3.243 | 13.607 | 8.884e-2/1.175 | 0 | 5.011/4.204 |
| A2 | FR     | 100% | 100% | 96% | 92% | 0% | 100% | 0% |
|    | Best   | 3.937e-27 | 3.321 | 8.639 | 45.553 | 0.343/9.422 | 3.660e+1 | 0/1.327e+4 |
|    | Median | 3.978e-27 | 3.589 | 14.923 | 64.403 | 1.008/11.000 | 3.660e+1 | 0/1.327e+4 |
|    | Worst  | 3.987 | 4.015 | 21.206 | 76.969 | 1.098/11.000 | 3.722e+1 | 2.466e+1/1.329e+4 |
|   | Mean | 4.780e-1 | 3.658 | 15.237 | 63.583 | 0.971/10.914 | 3.663e+1 | 3.417/1.327e+4 |
|   | STD | 1.295 | 0.192 | 3.568 | 9.704 | 0.164/0.325 | 0.1419 | 7.202/7.800 |
| A3 | FR | 100% | 92% | 100% | 100% | 0% | 44% | 0% |
|   | Best | 3.937e-27 | 2.376 | 2.356 | 0 | 0.116/11.000 | 3.660e+1 | 0/1.327e+4 |
|   | Median | 3.937e-27 | 2.377 | 2.356 | 0 | 1.003/11.000 | 3.661e+1 | 1.570e+1/1.329e+4 |
|   | Worst | 1.186e-26 | 2.377 | 8.639 | 0 | 1.024/11.000 | 5.618e+1 | 4.867e+1/1.334e+4 |
|   | Mean | 4.695e-27 | 2.377 | 3.018 | 0 | 0.862/11.000 | 3.940e+1 | 1.287e+1/1.328e+4 |
|   | STD | 1.895e-27 | 2.474e-4 | 1.637 | 0 | 2.515e-1/0 | 6.406 | 1.062e+1/1.420e+1 |
| A4 | FR | 100% | 100% | 100% | 100% | 0% | 0% | 0% |
|   | Best | 3.44e-27 | 2.376 | 5.498 | 3.142e+1 | 0.584/- | 8.050e+2/- | 1.030e-6/- |
|   | Median | 3.49e-27 | 2.403 | 8.639 | 3.927e+1 | 1.056/- | 1.845e+3/- | 2.727e+1/- |
|   | Worst | 3.987 | 3.346 | 1.178e+1 | 4.555e+1 | 1.123/- | 3.547e+3/- | 2.757e+1/- |
|   | Mean | 0.638 | 2.541 | 8.639 | 3.870e+1 | 0.972/- | 1.782e+4/- | 1.970e+1/- |
|   | STD | 1.492 | 0.249 | 0.907 | 5.188 | 0.323/- | 4.811e+3/- | 1.155e+1/- |

Table S5. Statistical results of A1-A4 on C01-C07 with $D = 30$

|   |   | C01 | C02 | C03 | C04 | C05 | C06 | C07 |
|---|---|---|---|---|---|---|---|---|
| A1 | FR | 100% | 100% | 92% | 100% | 100% | 100% | 100% |
|   | Best | 0 | 0 | 7.820e+2 | 1.357e+1 | 0 | 0 | -1.197e+3 |
|   | Median | 0 | 0 | 3.371e+3 | 1.592e+1 | 0 | 0 | -1.058e+3 |
|   | Worst | 3.550e-29 | 5.049e-29 | 3.228e+4 | 1.691e+1 | 0 | 0 | -9.044e+2 |
|   | Mean | 4.677e-30 | 5.569e-30 | 5.829e+3 | 1.502e+1 | 0 | 0 | -1.065e+3 |
|   | STD | 8.166e-30 | 1.169e-29 | 6.640e+3 | 1.198 | 0 | 0 | 7.096e+1 |
| A2 | FR | 100% | 100% | 24% | 100% | 100% | 0% | 80% |
|   | Best | 0 | 0 | 5.316e+5 | 1.357e+1 | 0 | 4.517e+3/1.130 | 38.254 |
|   | Median | 4.177e-29 | 4.720e-29 | 2.443e+6 | 1.357e+1 | 0 | 5.617e+3/1.771 | 53.537 |
|   | Worst | 6.910e-28 | 4.161e-27 | 1.074e+7 | 1.592e+1 | 0 | 7.839e+3/1.103 | 68.820 |
|   | Mean | 1.289e-28 | 4.163e-28 | 3.737e+6 | 1.470e+1 | 0 | 5.653e+3/1.513 | 53.537 |
|   | STD | 2.022e-28 | 9.078e-28 | 3.486e+6 | 1.172e+1 | 0 | 8.031e+2/0.548 | 15.283 |
| A3 | FR | 100% | 100% | 100% | 100% | 100% | 100% | 100% |
|   | Best | 0 | 0 | 0 | 1.357e+1 | 0 | 0 | -1.386e+3 |
|   | Median | 0 | 3.550e-29 | 4.577e-28 | 1.592e+1 | 0 | 0 | -1.386e+3 |
|   | Worst | 1.664e-28 | 1.264e-27 | 5.510e+2 | 1.691e+1 | 0 | 0 | -1.376e+3 |
|   | Mean | 5.064e-29 | 1.160e-28 | 7.700e+1 | 1.510e+1 | 0 | 0 | -1.385e+3 |
|   | STD | 5.089e-29 | 2.520e-28 | 1.386e+2 | 1.136 | 0 | 0 | 2.821 |
| A4 | FR | 100% | 100% | 100% | 100% | 100% | 0% | 96% |
|   | Best | 0 | 0 | 3.906e+4 | 1.357e+1 | 0 | 2.406e+3/- | -2.457e+2 |
|   | Median | 0 | 0 | 2.059e+5 | 1.357e+1 | 0 | 3.122e+3/- | -1.344e+2 |
|   | Worst | 2.080e-29 | 3.340e-29 | 2.180e+6 | 1.357e+1 | 0 | 4.071e+3/- | 8.163e+1 |
|   | Mean | 3.870e-30 | 5.260e-30 | 3.551e+5 | 1.357e+1 | 0 | 5.803e+3/- | -1.094e+2 |
|   | STD | 6.100e-30 | 8.390e-30 | 4.468e+5 | 0 | 0 | 9.815e+2/- | 8.874e+1 |

Table S6. Statistical results of A1-A4 on C08-C14 with $D = 30$

|   |   | C08 | C09 | C10 | C11 | C12 | C13 | C14 |
|---|---|---|---|---|---|---|---|---|
| A1 | FR | 100% | 100% | 100% | 20% | 100% | 100% | 100% |
|   | Best | -2.840e-4 | -2.666e-3 | -1.028e-4 | -9.244e-1 | 3.983 | 0 | 1.409 |
|   | Median | -2.839e-4 | -2.666e-3 | -1.028e-4 | -9.212e-1 | 3.983 | 0 | 1.409 |
|   | Worst | -2.835e-4 | -2.666e-3 | -1.028e-4 | -5.063e-1 | 3.983 | 8.166e+1 | 1.409 |
|   | Mean | -2.839e-4 | -2.666e-3 | -1.028e-4 | -7.714e-1 | 3.983 | 3.563e+1 | 1.409 |
|   | STD | 1.166e-7 | 0 | 0 | 1.868e-1 | 0 | 4.020e+1 | 0 |
| A2 | FR | 100% | 100% | 100% | 16% | 100% | 100% | 100% |
|   | Best | -2.798e-4 | -2.666e-3 | -1.028e-4 | -9.091e-1 | 3.983 | 0 | 2.018 |
|   | Median | -2.590e-4 | -2.666e-3 | -1.026e-4 | -7.194e-1 | 3.983 | 0 | 2.109 |
|   | Worst | -2.158e-4 | -2.666e-3 | -1.022e-4 | -5.209e-1 | 3.985 | 8.061e+1 | 2.225 |
|   | Mean | -2.588e-4 | -2.666e-3 | -1.026e-4 | -7.172e-1 | 3.983 | 9.673 | 2.112 |
|   | STD | 1.665e-5 | 0 | 1.626e-7 | 1.519e-1 | 6.299e-4 | 2.619e+1 | 5.824e-2 |
| A3 | FR | 100% | 56% | 100% | 0% | 100% | 100% | 96% |
|   | Best | 7.160e-5 | -2.663e-3 | 3.257e-5 | -7.118e+2/6.039e+1 | 3.983 | 0 | 1.409 |
|   | Median | 7.480e-4 | -2.616e-3 | 6.854e-4 | -5.838e+2/4.072e+1 | 3.983 | 0 | 1.409 |
|   | Worst | 2.622e-3 | -1.651e-3 | 1.410e-3 | -2.947e+1/1.213e-1 | 1.028e+1 | 5.838e-28 | 1.661 |
|   | Mean | 8.570e-4 | -2.501e-3 | 6.172e-4 | -4.988e+2/3.503e+1 | 5.060 | 4.670e-29 | 1.421 |
|   | STD | 6.015e-4 | 2.700e-4 | 3.421e-4 | 2.110e+2/1.789e+1 | 2.216 | 1.584e-28 | 5.073e-2 |

|  |  |  |  |  |  |  |  |
|---|---|---|---|---|---|---|---|
| A4 | FR | 100% | 100% | 100% | 100% | 100% | 100% | 100% |
|  | Best | -2.800e-4 | -2.666e-3 | -0.100e-3 | -0.925 | 3.983 | 0.234 | 1.565 |
|  | Median | -2.800e-4 | -2.666e-3 | -0.100e-3 | -0.925 | 3.985 | 5.907 | 1.839 |
|  | Worst | -2.800e-4 | -2.666e-3 | -0.100e-3 | -0.475 | 3.989 | 1.683e+1 | 1.937 |
|  | Mean | -2.800e-4 | -2.666e-3 | -0.100e-3 | -0.875 | 3.985 | 6.145 | 1.830 |
|  | STD | 0 | 0 | 0 | 0.110 | 1.443e-3 | 4.723 | 0.083 |

Table S7. Statistical results of A1-A4 on C15-C21 with $D = 30$

|  |  | C15 | C16 | C17 | C18 | C19 | C20 | C21 |
|---|---|---|---|---|---|---|---|---|
| A1 | FR | 96% | 100% | 0% | 100% | 0% | 100% | 100% |
|  | Best | 2.356 | 1.257e+1 | 9.417e-1/3.128e+1 | 3.652e+1 | -2.776e+1/4.283e+4 | 1.280 | 3.983 |
|  | Median | 2.356 | 1.257e+1 | 1.023/3.103e+1 | 3.652e+1 | -2.776e+1/4.283e+4 | 1.657 | 5.843 |
|  | Worst | 1.806e+1 | 1.257e+1 | 1.031/3.103e+1 | 3.652e+1 | -2.776e+1/4.283e+4 | 2.001 | 1.026e+1 |
|  | Mean | 3.403 | 1.257e+1 | 1.014/3.111e+1 | 3.652e+1 | -2.776e+1/4.283e+4 | 1.637 | 6.255 |
|  | STD | 3.353 | 0 | 2.153e-2/8.668e-2 | 0 | 0/0 | 1.767 | 2.151 |
| A2 | FR | 64% | 100% | 0% | 100% | 0% | 100% | 100% |
|  | Best | 1.492e+1 | 1.759e+2 | 1.026/3.100e+1 | 3.652e+1 | 0/4.275e+4 | 2.081 | 1.162e+1 |
|  | Median | 2.435e+1 | 1.963e+2 | 1.030/3.100e+1 | 3.652e+1 | 0/4.275e+4 | 2.311 | 1.706e+1 |
|  | Worst | 2.749e+1 | 2.340e+2 | 1.066/3.100e+1 | 3.653e+1 | 0/4.275e+4 | 2.607 | 2.533e+1 |
|  | Mean | 2.297e+1 | 2.001e+2 | 1.033/3.100e+1 | 3.652e+1 | 0/4.275e+4 | 2.300 | 1.683e+1 |
|  | STD | 2.932 | 1.266e+1 | 7.810e-3/0 | 1.870e-3 | 0/0 | 1.297 | 3.273 |
| A3 | FR | 88% | 100% | 0% | 76% | 0% | 100% | 100% |
|  | Best | 2.356 | 0 | 3.310e-2/3.100e+1 | 3.652e+1 | 0/4.275e+4 | 2.091 | 3.983 |
|  | Median | 7.069 | 0 | 1.078/2.900e+1 | 3.652e+1 | 0/4.275e+4 | 2.284 | 1.009e+1 |
|  | Worst | 2.435e+1 | 0 | 1.308/2.900e+1 | 3.652e+1 | 0/4.275e+4 | 2.453 | 1.977e+1 |
|  | Mean | 9.782 | 0 | 9.294e-1/2.996e+1 | 3.652e+1 | 0/4.275e+4 | 2.278 | 9.978 |
|  | STD | 6.931 | 0 | 3.758e-1/9.992e-1 | 0 | 0/0 | 1.008e-1 | 4.709 |
| A4 | FR | 84% | 100% | 0% | 0% | 0% | 100% | 100% |
|  | Best | 1.178e+1 | 1.147e+2 | 0.948/- | 2.912e+2/- | 5.600e-6/- | 1.275 | 9.775 |
|  | Median | 1.806e+1 | 1.508e+2 | 1.001/- | 2.082e+3/- | 6.160e-6/- | 1.451 | 2.833e+1 |
|  | Worst | 1.920e+1 | 1.696e+2 | 1.018/- | 1.494e+4/- | 6.920e-6/- | 1.553 | 3.209e+1 |
|  | Mean | 1.806e+1 | 1.463e+2 | 0.997/- | 3.046e+3/- | 6.200e-6/- | 1.417 | 2.251e+1 |
|  | STD | 3.614 | 1.403e+1 | 0.167e-1/- | 4.348e+3/- | 3.400e-7/- | 0.991e-1 | 8.908 |

Table S8. Statistical results of A1-A4 on C22-C28 with $D = 30$

|  |  | C22 | C23 | C24 | C25 | C26 | C27 | C28 |
|---|---|---|---|---|---|---|---|---|
| A1 | FR | 100% | 100% | 100% | 100% | 0% | 72% | 0% |
|  | Best | 1.325e+1 | 1.415 | 2.356 | 1.885e+1 | 7.901e-1/3.247e+1 | 3.653e+1 | 1.303e+2/4.304e+4 |
|  | Median | 2.784e+2 | 1.444 | 5.498 | 2.513e+1 | 1.018/3.117e+1 | 3.653e+1 | 1.506e+2/4.303e+4 |
|  | Worst | 1.466e+3 | 1.508 | 8.639 | 4.556e+1 | 1.038/3.108e+1 | 3.882e+1 | 1.687e+2/4.304e+4 |
|  | Mean | 4.642e+2 | 1.448 | 5.105 | 2.765e+1 | 1.002/3.123e+1 | 3.666e+1 | 1.485e+2/4.303e+4 |
|  | STD | 4.736e+2 | 1.905e-2 | 1.651 | 6.796 | 4.812e-2/3.073e-1 | 5.251e-1 | 9.994/4.139 |
| A2 | FR | 0% | 100% | 100% | 100% | 0% | 92% | 0% |
|  | Best | 1.029e+5/7.665e+1 | 2.131 | 1.806e+1 | 1.696e+2 | 1.017/3.100e+1 | 3.653e+1 | 1.475e+2/4.298e+4 |
|  | Median | 2.332e+5/6.774e+1 | 2.229 | 2.121e+1 | 2.073e+2 | 1.030/3.100e+1 | 3.653e+1 | 1.804e+2/4.299e+4 |
|  | Worst | 3.740e+5/8.522e+1 | 2.324 | 2.121e+1 | 2.278e+2 | 1.047/3.100e+1 | 3.686e+1 | 1.991e+2/4.299e+4 |
|  | Mean | 2.321e+5/7.854e+1 | 2.229 | 1.970e+1 | 2.061e+2 | 1.030/3.100e+1 | 3.655e+1 | 1.762e+2/4.299e+4 |
|  | STD | 6.982e+4/1.269e+1 | 5.118e-2 | 1.570 | 1.453e+1 | 4.596e-3/0 | 6.584e-2 | 1.486e+1/4.224 |
| A3 | FR | 100% | 100% | 100% | 100% | 0% | 28% | 0% |
|  | Best | 1.503e+1 | 1.418 | 5.498 | 0 | 9.669e-2/3.100e+1 | 3.692e+1 | 1.256e+2/4.294e+4 |
|  | Median | 1.762e+1 | 1.463 | 1.492e+1 | 5.308e-14 | 1.027/3.100e+1 | 4.841e+1 | 1.466e+2/4.297e+4 |
|  | Worst | 8.598e+1 | 1.884 | 2.749e+1 | 7.218e-12 | 1.031/3.100e+1 | 5.553e+1 | 1.784e+2/4.298e+4 |
|  | Mean | 2.037e+1 | 1.523 | 1.367e+1 | 6.529e-13 | 9.728e-1/3.100e+1 | 4.690e+1 | 1.493e+2/4.296e+4 |
|  | STD | 1.352e+1 | 1.183e-1 | 4.785 | 1.553e-12 | 1.833e-1/0 | 6.696 | 1.351e+1/1.721e+1 |
| A4 | FR | 100% | 100% | 100% | 100% | 0% | 0% | 0% |
|  | Best | 3.156 | 1.500 | 1.178e+1 | 1.131e+2 | 0.955/- | 1.186e+3/- | 4.521e+1/- |
|  | Median | 7.959e+2 | 1.741 | 1.178e+1 | 1.461e+2 | 0.999/- | 3.987e+3/- | 1.433e+2/- |
|  | Worst | 4.865e+3 | 1.879 | 1.492e+1 | 1.587e+2 | 1.020/- | 1.236e+4/- | 1.441e+2/- |
|  | Mean | 1.268e+3 | 1.715 | 1.216e+1 | 1.422e+2 | 0.997/- | 9.318e+3/- | 1.335e+2/- |
|  | STD | 1.408e+3 | 0.127 | 1.042 | 1.393e+1 | 1.640e-2/- | 8.788e+3/- | 2.447e+1/- |

Table S9. Statistical results of A1-A4 on C01-C07 with $D = 50$

|    |        | C01       | C02       | C03       | C04       | C05       | C06           | C07       |
|----|--------|-----------|-----------|-----------|-----------|-----------|---------------|-----------|
| A1 | FR     | 100%      | 100%      | 68%       | 100%      | 100%      | 100%          | 72%       |
|    | Best   | 1.255e-27 | 5.567e-28 | 1.713e+4  | 1.357e+1  | 0         | 0             | -8.957e+2 |
|    | Median | 4.785e-26 | 2.775e-26 | 5.212e+4  | 1.592e+1  | 0         | 0             | -6.569e+2 |
|    | Worst  | 4.923e-24 | 8.692e-24 | 1.033e+6  | 1.691e+1  | 0         | 0             | -3.139e+2 |
|    | Mean   | 4.109e-25 | 5.170e-25 | 1.368e+5  | 1.483e+1  | 0         | 0             | -6.432e+2 |
|    | STD    | 9.648e-25 | 1.724e-24 | 2.484e+5  | 1.226     | 0         | 0             | 1.260e+2  |
| A2 | FR     | 100%      | 100%      | 40%       | 100%      | 100%      | 0%            | 12%       |
|    | Best   | 2.837e-22 | 8.149e-22 | 2.369e+6  | 1.357e+1  | 0         | 7.839e+3/1.026| -6.797e+1 |
|    | Median | 2.399e-20 | 1.292e-20 | 5.161e+6  | 1.592e+1  | 0         | 9.329e+3/0.973| 1.470e+2  |
|    | Worst  | 4.204e-18 | 1.810e-18 | 1.526e+7  | 1.691e+1  | 0         | 1.086e+4/0.565| 3.459e+2  |
|    | Mean   | 2.260e-19 | 1.152e-19 | 6.576e+6  | 1.502e+1  | 0         | 9.347e+3/0.924| 1.416e+2  |
|    | STD    | 8.186e-19 | 3.528e-19 | 4.180e+6  | 1.198     | 0         | 7.922e+2/0.247| 1.690e+2  |
| A3 | FR     | 100%      | 100%      | 100%      | 100%      | 100%      | 100%          | 100%      |
|    | Best   | 8.441e-23 | 1.000e-22 | 1.961e-17 | 1.357e+1  | 0         | 0             | -2.305e+3 |
|    | Median | 4.420e-21 | 3.536e-21 | 6.521e+2  | 1.592e+1  | 0         | 0             | -2.298e+3 |
|    | Worst  | 2.080e-18 | 4.807e-18 | 4.900e+3  | 1.691e+1  | 0         | 0             | -2.272e+3 |
|    | Mean   | 1.986e-19 | 2.024e-19 | 1.245e+3  | 1.502e+1  | 0         | 0             | -2.297e+3 |
|    | STD    | 4.589e-19 | 9.400e-19 | 1.247e+3  | 1.198     | 0         | 0             | 7.246     |
| A4 | FR     | 100%      | 100%      | 100%      | 100%      | 100%      | 0%            | 100%      |
|    | Best   | 8.680e-30 | 2.500e-29 | 2.867e+5  | 1.357e+1  | 0         | 6.637e+3/-    | -3.620e+2 |
|    | Median | 7.730e-29 | 1.020e-28 | 6.337e+5  | 1.357e+1  | 1.300e-28 | 7.515e+3/-    | -1.738e+2 |
|    | Worst  | 1.420e-28 | 1.780e-28 | 3.870e+6  | 1.357e+1  | 6.400e-28 | 6.709e+3/-    | 1.437e+1  |
|    | Mean   | 7.790e-29 | 9.790e-29 | 8.945e+5  | 1.357e+1  | 1.680e-28 | 8.637e+3/-    | -1.791e+2 |
|    | STD    | 3.080e-29 | 4.600e-29 | 7.405e+5  | 0         | 1.590e-28 | 1.418e+3/-    | 8.975e+1  |

Table S10. Statistical results of A1-A4 on C08-C14 with $D = 50$

|    |        | C08       | C09       | C10       | C11                | C12      | C13       | C14     |
|----|--------|-----------|-----------|-----------|--------------------|----------|-----------|---------|
| A1 | FR     | 100%      | 100%      | 100%      | 16%                | 100%     | 100%      | 100%    |
|    | Best   | -1.290e-4 | -2.037e-3 | -4.814e-5 | -1.824             | 3.981    | 0         | 1.100   |
|    | Median | -1.142e-4 | -2.037e-3 | -4.766e-5 | -1.184             | 3.981    | 9.935e+1  | 1.100   |
|    | Worst  | -9.383e-5 | -2.037e-3 | -4.617e-5 | -1.002             | 3.982    | 1.005e+2  | 1.100   |
|    | Mean   | -1.143e-4 | -2.037e-3 | -4.752e-5 | -1.298             | 3.981    | 5.585e+1  | 1.100   |
|    | STD    | 1.017e-5  | 0         | 4.479e-7  | 3.330e-1           | 1.255e-5 | 4.950e+1  | 0       |
| A2 | FR     | 100%      | 100%      | 100%      | 0%                 | 100%     | 100%      | 100%    |
|    | Best   | 6.309e-4  | -2.037e-3 | -4.575e-5 | -3.607e+1/8.118e-8 | 3.982    | 0         | 1.367   |
|    | Median | 1.954e-3  | -2.037e-3 | -4.193e-5 | 1.919e+1/1.613e-21 | 3.982    | 0         | 1.593   |
|    | Worst  | 4.269e-3  | -2.037e-3 | -3.228e-5 | 3.776e+1/7.944e-7  | 3.986    | 1.009e+2  | 1.671   |
|    | Mean   | 2.022e-3  | -2.037e-3 | -4.127e-5 | 1.320e+1/3.639e-8  | 3.982    | 4.405e+1  | 1.584   |
|    | STD    | 1.005e-3  | 0         | 3.0569e-6 | 1.897e+1/1.556e-7  | 7.948e-4 | 4.969e+1  | 6.312e-2|
| A3 | FR     | 100%      | 100%      | 100%      | 0%                 | 100%     | 100%      | 96%     |
|    | Best   | 1.393e-3  | -2.033e-3 | 2.037e-4  | -2.594e+3/1.677e+2 | 3.981    | 0         | 1.100   |
|    | Median | 3.009e-3  | -1.947e-3 | 5.956e-4  | -2.494e+3/1.553e+2 | 3.982    | 0         | 1.100   |
|    | Worst  | 1.330e-2  | -1.529e-3 | 4.428e-3  | -5.339e+2/2.995e+1 | 7.889    | 1.275e-27 | 1.168   |
|    | Mean   | 3.632e-3  | -1.926e-3 | 7.836e-4  | -2.416e+3/1.505e+2 | 4.512    | 5.099e-29 | 1.109   |
|    | STD    | 2.493e-3  | 1.112e-4  | 7.699e-4  | 3.909e+2/2.606e+1  | 1.224    | 2.498e-28 | 1.764e-2|
| A4 | FR     | 100%      | 100%      | 100%      | 100%               | 100%     | 100%      | 100%    |
|    | Best   | -1.300e-5 | -2.037e-3 | -4.830e-5 | -2.011             | 2.049e+1 | 3.039e+1  | 1.320   |
|    | Median | -1.300e-5 | -2.037e-3 | -4.830e-5 | -2.011             | 6.231e+1 | 3.510e+1  | 1.414   |
|    | Worst  | -1.300e-5 | -2.037e-3 | -4.830e-5 | -1.026             | 6.478e+1 | 3.900e+1  | 1.453   |
|    | Mean   | -1.300e-5 | -2.037e-3 | -4.830e-5 | -1.765             | 4.976e+1 | 2.670e+1  | 1.404   |
|    | STD    | 0         | 0         | 0         | 0.334              | 1.994e+1 | 1.359e+1  | 0.375e-1|

Table S11. Statistical results of A1-A4 on C15-C21 with $D = 50$

|    |        | C15   | C16      | C17             | C18     | C19               | C20      | C21     |
|----|--------|-------|----------|-----------------|---------|-------------------|----------|---------|
| A1 | FR     | 100%  | 100%     | 0%              | 88%     | 0%                | 100%     | 100%    |
|    | Best   | 2.356 | 1.885e+1 | 1.030/5.131e+1  | 3.647e+1| -4.604e+1/7.237e+4| 2.558    | 4.798   |
|    | Median | 2.356 | 2.042e+1 | 1.050/5.101e+1  | 3.647e+1| -4.604e+1/7.237e+4| 3.029    | 5.402   |
|    | Worst  | 8.639 | 3.142e+1 | 1.052/5.115e+1  | 3.722e+1| -4.604e+1/7.237e+4| 3.339    | 8.472   |
|    | Mean   | 3.110 | 2.161e+1 | 1.048/5.113e+1  | 3.651e+1| -4.604e+1/7.237e+4| 2.969    | 5.576   |
|    | STD    | 1.609 | 3.421    | 5.425e-3/1.091e-1| 0.155  | 0/0               | 2.017e-1 | 7.461e-1|

|   |        | C22         | C23        | C24             | C25         | C26               | C27        | C28              |
|---|--------|-------------|------------|-----------------|-------------|-------------------|------------|------------------|
| A2 | FR     | 64%         | 100%       | 0%              | 100%        | 0%                | 100%       | 100%             |
|    | Best   | 1.492e+1    | 3.267e+2   | 1.044/5.100e+1  | 3.647e+1    | 0/7.223e+4        | 3.571      | 8.243            |
|    | Median | 2.749e+1    | 3.519e+2   | 1.050/5.100e+1  | 3.648e+1    | 0/7.223e+4        | 3.839      | 1.074e+1         |
|    | Worst  | 3.691e+1    | 3.770e+2   | 1.050/5.100e+1  | 3.722e+1    | 0/7.223e+4        | 4.093      | 1.806e+1         |
|    | Mean   | 2.611e+1    | 3.517e+2   | 1.049/5.100e+1  | 3.651e+1    | 0/7.223e+4        | 3.833      | 1.174e+1         |
|    | STD    | 5.874       | 1.209e+1   | 1.372e-3/0      | 0.145       | 0/0               | 1.276      | 2.333            |
| A3 | FR     | 100%        | 100%       | 0%              | 48%         | 0%                | 100%       | 100%             |
|    | Best   | 8.639       | 0          | 6.642e-2/5.100e+1 | 3.647e+1  | 0/7.223e+4        | 3.534      | 7.405            |
|    | Median | 1.492e+1    | 0          | 9.905e-1/5.100e+1 | 3.647e+1  | 0/7.223e+4        | 3.812      | 1.228e+1         |
|    | Worst  | 3.063e+1    | 0          | 1.045/5.100e+1  | 3.649e+1    | 0/7.223e+4        | 4.026      | 2.042e+1         |
|    | Mean   | 1.580e+1    | 0          | 8.256e-1/5.084e+1 | 3.648e+1  | 0/7.223e+4        | 3.794      | 1.283e+1         |
|    | STD    | 6.440       | 0          | 2.808e-1/5.426e-1 | 4.690e-3  | 0/0               | 1.261e-1   | 3.600            |
| A4 | FR     | 100%        | 100%       | 0%              | 0%          | 0%                | 100%       | 100%             |
|    | Best   | 1.178e+1    | 2.215e+2   | 1.026/-         | 7.830e+2/-  | 1.070e-5/-        | 2.8862     | 6.217e+1         |
|    | Median | 1.806e+1    | 2.529e+2   | 1.038/-         | 5.667e+3/-  | 1.190e-5/-        | 3.215      | 6.217e+1         |
|    | Worst  | 2.435e+1    | 2.890e+2   | 1.044/-         | 1.823e+4/-  | 1.680e-5/-        | 3.448      | 6.717e+1         |
|    | Mean   | 1.781e+1    | 2.531e+2   | 1.037/-         | 8.031e+3/-  | 1.210e-5/-        | 3.203      | 6.286e+1         |
|    | STD    | 2.997       | 1.621e+1   | 4.681e-3/-      | 6.844e+3/-  | 1.120e-6/-        | 0.143      | 1.590            |

Table S12. Statistical results of A1-A4 on C22-C28 with $D = 50$

|   |        | C22                | C23       | C24       | C25        | C26                | C27       | C28                  |
|---|--------|--------------------|-----------|-----------|------------|--------------------|-----------|----------------------|
| A1 | FR     | 80%                | 100%      | 100%      | 100%       | 0%                 | 80%       | 0%                   |
|    | Best   | 2.106e+3           | 1.121     | 8.639     | 4.398e+1   | 1.022/5.102e+1     | 4.036e+1  | 2.464e+2/7.271e+4    |
|    | Median | 8.572e+3           | 1.135     | 1.178e+1  | 6.440e+1   | 1.049/5.110e+1     | 4.418e+1  | 2.704e+2/7.272e+4    |
|    | Worst  | 1.444e+4           | 1.159     | 1.492e+1  | 8.796e+1   | 1.051/5.101e+1     | 4.801e+1  | 2.948e+2/7.272e+4    |
|    | Mean   | 7.991e+3           | 1.136     | 1.040e+1  | 6.987e+1   | 1.047/5.112e+1     | 4.418e+1  | 2.718e+2/7.271e+4    |
|    | STD    | 3.677e+3           | 9.825e-3  | 1.795     | 1.168e+1   | 6.328e-3/1.269e-1  | 3.822     | 1.130e+1/5.054       |
| A2 | FR     | 0%                 | 100%      | 100%      | 100%       | 0%                 | 56%       | 0%                   |
|    | Best   | 3.712e+5/3.259e+2  | 1.592     | 1.806e+1  | 2.953e+2   | 1.040/5.100e+1     | 3.647e+1  | 2.699e+2/7.266e+4    |
|    | Median | 8.657e+5/3.644e+2  | 1.688     | 2.121e+1  | 3.471e+2   | 1.049/5.100e+1     | 3.649e+1  | 2.989e+2/7.266e+4    |
|    | Worst  | 1.454e+6/3.597e+2  | 1.730     | 2.435e+1  | 3.723e+2   | 1.051/5.100e+1     | 5.252e+1  | 3.328e+2/7.266e+4    |
|    | Mean   | 8.725e+5/3.435e+2  | 1.679     | 2.083e+1  | 3.468e+2   | 1.049/5.100e+1     | 3.763e+1  | 2.985e+2/7.266e+4    |
|    | STD    | 2.535e+5/3.602e+1  | 3.535e-2  | 1.353     | 1.656e+1   | 2.301e-3/0         | 4.129e+1  | 1.411e+1/5.451       |
| A3 | FR     | 100%               | 100%      | 100%      | 100%       | 0%                 | 16%       | 0%                   |
|    | Best   | 3.550e+1           | 1.115     | 1.806e+1  | 5.543e-12  | 1.033/5.100e+1     | 3.717e+1  | 2.539e+2/7.263e+4    |
|    | Median | 3.968e+1           | 1.408     | 2.121e+1  | 1.885e+1   | 1.049/5.100e+1     | 4.694e+1  | 2.913e+2/7.264e+4    |
|    | Worst  | 5.878e+2           | 1.713     | 2.435e+1  | 3.157e+2   | 1.051/5.100e+1     | 5.663e+1  | 3.173e+2/7.265e+4    |
|    | Mean   | 6.866e+1           | 1.420     | 2.007e+1  | 5.894e+1   | 1.047/5.100e+1     | 4.692e+1  | 2.883e+2/7.263e+4    |
|    | STD    | 1.116e+2           | 1.709     | 1.963     | 7.686e+1   | 5.064e-3/9.370e-5  | 8.605     | 1.612e+1/1.752e+1    |
| A4 | FR     | 96%                | 100%      | 100%      | 100%       | 0%                 | 0%        | 0%                   |
|    | Best   | 2.669e+3           | 1.163     | 1.178e+1  | 2.199e+2   | 1.029/-            | 1.647e+3/-| 2.283e+2/-           |
|    | Median | 4.797e+3           | 1.342     | 1.492e+1  | 2.450e+2   | 1.040/-            | 2.173e+4/-| 2.659e+2/-           |
|    | Worst  | 2.607e+4           | 1.439     | 1.492e+1  | 2.780e+2   | 1.046/-            | 1.085e+5/-| 2.817e+2/-           |
|    | Mean   | 8.393e+3           | 1.338     | 1.429e+1  | 2.485e+2   | 1.040/-            | 1.157e+4/-| 2.650e+2/-           |
|    | STD    | 6.738e+3           | 0.616e-1  | 1.283     | 1.584e+1   | 3.809e-3/-         | 2.554e+4/-| 1.585e+1/-           |

Table S13. Statistical results of A1-A4 on C01-C07 with $D = 100$

|   |        | C01        | C02        | C03       | C04       | C05       | C06               | C07        |
|---|--------|------------|------------|-----------|-----------|-----------|-------------------|------------|
| A1 | FR     | 100%       | 100%       | 64%       | 100%      | 100%      | 92%               | 32%        |
|    | Best   | 5.182e-15  | 4.345e-14  | 1.172e+5  | 1.357e+1  | 0         | 0                 | -8.810e+2  |
|    | Median | 7.301e-14  | 2.429e-13  | 2.967e+5  | 1.592e+1  | 0         | 0                 | -6.782e+2  |
|    | Worst  | 1.302e-11  | 1.712e-11  | 1.060e+6  | 1.691e+1  | 0         | 0                 | -2.149e+2  |
|    | Mean   | 9.852e-13  | 1.388e-12  | 3.843e+5  | 1.515e+1  | 0         | 0                 | -6.417e+2  |
|    | STD    | 2.935e-12  | 3.357e-12  | 2.782e+5  | 1.215     | 0         | 0                 | 1.951e+2   |
| A2 | FR     | 100%       | 100%       | 64%       | 100%      | 100%      | 0%                | 20%        |
|    | Best   | 9.986e-12  | 1.935e-11  | 5.989e+6  | 1.357e+1  | 0         | 1.537e+4/9.123e-1 | 1.477e+2   |
|    | Median | 1.914e-10  | 1.374e-9   | 1.844e+7  | 1.534e+1  | 0         | 1.874e+4/7.209e-1 | 2.248e+2   |
|    | Worst  | 9.599e-9   | 5.176e-8   | 7.874e+7  | 1.691e+1  | 1.183e-28 | 2.094e+4/8.357e-1 | 3.731e+2   |
|    | Mean   | 8.114e-10  | 7.194e-9   | 3.250e+7  | 1.485e+1  | 7.425e-30 | 1.842e+4/6.383e-1 | 2.348e+2   |
|    | STD    | 1.889e-9   | 1.392e-8   | 2.446e+7  | 1.261     | 2.377e-29 | 1.251e+3/2.011e-1 | 7.559e+1   |
| A3 | FR     | 100%       | 100%       | 96%       | 100%      | 100%      | 96%               | 100%       |
|    | Best   | 2.4738e-11 | 2.932e-12  | 1.514e+3  | 1.357e+1  | 0         | 0                 | -4.653e+3  |
|    | Median | 1.539e-10  | 1.421e-10  | 5.608e+3  | 1.592e+1  | 0         | 0                 | -4.508e+3  |

|   |        | Worst  | 4.338e-9   | 4.606e-9   | 2.495e+5 | 1.691e+1 | 9.880e-29 | 0        | -4.233e+3 |
|---|--------|--------|------------|------------|----------|----------|-----------|----------|-----------|
|   |        | Mean   | 4.782e-10  | 5.499e-10  | 2.248e+4 | 1.491e+1 | 5.810e-30 | 0        | -4.500e+3 |
|   |        | STD    | 9.303e-10  | 9.688e-10  | 4.942e+4 | 1.250    | 2.000e-29 | 0        | 9.690e+1  |
| A4 |       | FR     | 100%       | 100%       | 100%     | 100%     | 100%      | 100%     | 0%        |
|   |        | Best   | 1.300e-26  | 1.310e-26  | 1.340e+6 | 1.357e+1 | 4.340e-7  | 1.716e+4 | -6.832e+2/- |
|   |        | Median | 4.500e-26  | 4.590e-26  | 2.470e+6 | 1.357e+1 | 4.900e-6  | 1.580e+4 | -3.012e+2/- |
|   |        | Worst  | 6.510e-25  | 6.380e-25  | 4.790e+6 | 1.557e+1 | 0.416e-3  | 1.672e+4 | -1.385e+1/- |
|   |        | Mean   | 1.030e-25  | 8.470e-26  | 2.730e+6 | 1.371e+1 | 3.280e-5  | 1.556e+4 | -3.027e+2/- |
|   |        | STD    | 1.680e-25  | 1.250e-25  | 9.656e+5 | 0.463    | 9.250e-5  | 1.595e+3 | 1.350e+2/- |

Table S14. Statistical results of A1-A4 on C08-C14 with $D = 100$

|     |        | C08                  | C09      | C10       | C11                   | C12      | C13      | C14      |
|-----|--------|----------------------|----------|-----------|-----------------------|----------|----------|----------|
| A1  | FR     | 0%                   | 100%     | 100%      | 0%                    | 100%     | 100%     | 100%     |
|     | Best   | 6.266e+1/8.201e+4    | 0        | -7.807e-6 | -1.252e+3/4.513e+1    | 3.981    | 3.375e+1 | 7.844e-1 |
|     | Median | 7.973e+1/7.512e+4    | 0        | -2.018e-6 | -9.581e+1/1.715       | 3.981    | 1.229e+2 | 7.869e-1 |
|     | Worst  | 9.652e+1/1.139e+5    | 6.521e-8 | 1.035e-5  | 9.762e+2/6.536        | 5.146    | 1.479e+2 | 7.908e-1 |
|     | Mean   | 7.961e+1/9.028e+4    | 2.608e-9 | -1.091e-6 | -7.814e+1/1.340e+1    | 4.119    | 9.211e+1 | 7.871e-1 |
|     | STD    | 7.872/1.379e+4       | 1.278e-8 | 4.701e-6  | 6.682e+2/1.159e+1     | 3.735e-1 | 5.555e+1 | 1.706e-3 |
| A2  | FR     | 0%                   | 100%     | 100%      | 0%                    | 100%     | 100%     | 100%     |
|     | Best   | 1.250/6.256          | 0        | 8.053e-6  | -3.167e+1/6.341e-2    | 3.981    | 3.390e+1 | 9.650e-1 |
|     | Median | 3.239/2.345e+1       | 2.018e-9 | 2.593e-5  | 5.356e+2/2.654        | 5.363    | 3.452e+1 | 1.072    |
|     | Worst  | 8.056/5.088e+1       | 3.047e-5 | 5.432e-5  | 1.016e+3/6.023        | 3.365e+1 | 1.478e+2 | 1.112    |
|     | Mean   | 3.723/2.321e+1       | 1.756e-6 | 2.617e-5  | 4.412e+2/1.723        | 9.036    | 5.275e+1 | 1.065    |
|     | STD    | 1.711/1.584e+1       | 6.000e-6 | 1.060e-5  | 2.896e+2/1.926        | 6.602    | 4.132e+1 | 3.501e-2 |
| A3  | FR     | 0%                   | 96%      | 96%       | 0%                    | 100%     | 100%     | 56%      |
|     | Best   | 1.295e-2/4.563e-4    | 0        | 1.765e-4  | -7.657e+3/2.509e+2    | 3.981    | 0        | 7.842e-1 |
|     | Median | 7.427e-2/1.016e-2    | 0        | 6.338e-4  | -7.136e+3/1.781e+2    | 3.982    | 0        | 7.842e-1 |
|     | Worst  | 3.878e-1/3.510e-1    | 1.131e-5 | 1.111e-3  | -5.067e+3/2.489e+2    | 5.548    | 2.066e+2 | 7.912e-1 |
|     | Mean   | 1.101e-1/3.99e-2     | 5.975e-7 | 6.409e-4  | -6.555e+3/2.373e+2    | 4.140    | 3.909e+1 | 7.847e-1 |
|     | STD    | 9.005e-2/7.154e-2    | 2.254e-6 | 2.464e-4  | 9.329e+2/7.911e+1     | 4.324e-1 | 7.829e+1 | 1.791e-3 |
| A4  | FR     | 100%                 | 100%     | 100%      | 88%                   | 100%     | 100%     | 100%     |
|     | Best   | -4.830e-5            | -1.430e-3 | -1.72e-5 | -3.935e+1             | 3.159e+1 | 6.688e+1 | 0.922    |
|     | Median | -4.820e-5            | -1.430e-3 | -1.72e-5 | -5.725                | 3.244e+1 | 7.758e+1 | 0.975    |
|     | Worst  | -4.770e-5            | -1.430e-3 | -1.71e-5 | -3.499                | 3.490e+1 | 9.391e+1 | 1.003    |
|     | Mean   | -4.810e-5            | -1.430e-3 | -1.72e-5 | -3.651                | 3.251e+1 | 8.070e+1 | 0.972    |
|     | STD    | 1.330e-7             | 0        | 1.29e-8   | 8.939                 | 0.819    | 7.379    | 1.939e-2 |

Table S15. Statistical results of A1-A4 on C15-C21 with $D = 100$

|     |        | C15      | C16      | C17                  | C18              | C19                 | C20      | C21      |
|-----|--------|----------|----------|----------------------|------------------|---------------------|----------|----------|
| A1  | FR     | 100%     | 100%     | 0%                   | 100%             | 0%                  | 100%     | 100%     |
|     | Best   | 5.498    | 4.398e+1 | 4.139e-2/3.399e+2    | 3.639e+1         | -9.230e+1/1.462e+5  | 5.639    | 6.131    |
|     | Median | 8.639    | 5.184e+1 | 9.830e-1/1.010e+2    | 3.639e+1         | -9.230e+1/1.462e+5  | 6.323    | 6.510    |
|     | Worst  | 2.435e+1 | 6.283e+1 | 1.029/1.010e+2       | 3.640e+1         | -9.230e+1/1.462e+5  | 6.746    | 1.239e+1 |
|     | Mean   | 9.142    | 5.278e+1 | 9.396e-1/1.106e+2    | 3.639e+1         | -9.230e+1/1.462e+5  | 6.303    | 7.240    |
|     | STD    | 4.232    | 5.066    | 1.869e-1/4.682e+1    | 3.906e-3         | 1.946e-8/0          | 2.743e-1 | 1.688    |
| A2  | FR     | 68%      | 100%     | 0%                   | 100%             | 0%                  | 100%     | 100%     |
|     | Best   | 2.435e+1 | 6.912e+2 | 1.099/1.010e+2       | 3.640e+1         | 0/1.459e+5          | 7.247    | 1.010e+1 |
|     | Median | 3.063e+1 | 7.367e+2 | 1.100/1.010e+2       | 3.641e+1         | 0/1.459e+5          | 7.625    | 1.472e+1 |
|     | Worst  | 3.691e+1 | 7.556e+2 | 1.100/1.010e+2       | 3.648e+1         | 0/1.459e+5          | 8.033    | 1.951e+1 |
|     | Mean   | 3.082e+1 | 7.332e+2 | 1.100/1.010e+2       | 3.642e+1         | 0/1.459e+5          | 7.611    | 1.461e+1 |
|     | STD    | 3.805    | 1.575e+1 | 1.951e-4/0           | 2.162e-2         | 0/0                 | 1.823e-1 | 2.311    |
| A3  | FR     | 92%      | 100%     | 0%                   | 28%              | 0%                  | 100%     | 100%     |
|     | Best   | 1.492e+1 | 0        | 1.244e-1/1.010e+2    | 3.639e+1         | 0/1.459e+5          | 7.213    | 1.355e+1 |
|     | Median | 2.121e+1 | 0        | 1.022/1.010e+2       | 3.641e+1         | 0/1.459e+5          | 7.643    | 1.822e+1 |
|     | Worst  | 4.006e+1 | 0        | 1.100/1.010e+2       | 3.642e+1         | 0/1.459e+5          | 8.044    | 3.128e+1 |
|     | Mean   | 2.244e+1 | 0        | 9.282e-1/1.010e+2    | 3.641e+1         | 0/1.459e+5          | 7.616    | 1.909e+1 |
|     | STD    | 6.728    | 0        | 2.427e-1/0           | 1.036e-2         | 0/0                 | 2.340e-1 | 3.722    |
| A4  | FR     | 100%     | 100%     | 0%                   | 0%               | 0%                  | 100%     | 100%     |
|     | Best   | 1.492e+1 | 4.775e+2 | 1.093/-              | 2.284e+3/-       | 3.500e-5/-          | 8.726    | 3.158e+1 |
|     | Median | 1.806e+1 | 5.341e+2 | 1.097/-              | 3.342e+3/-       | 4.560e-5/-          | 9.440    | 3.158e+1 |
|     | Worst  | 2.121e+1 | 6.032e+2 | 1.098/-              | 4.198e+3/-       | 5.710e-5/-          | 1.007e+1 | 3.159e+1 |
|     | Mean   | 1.806e+1 | 5.345e+2 | 1.096/-              | 3.436e+3/-       | 4.680e-5/-          | 9.363    | 3.158e+1 |
|     | STD    | 1.283    | 3.076e+1 | 1.201e-3/-           | 9.386e+2/-       | 5.780e-6/-          | 0.378    | 2.987e-3 |

Table S16. Statistical results of A1-A4 on C22-C28 with $D = 100$

|    |        | C22 | C23 | C24 | C25 | C26 | C27 | C28 |
|----|--------|-----|-----|-----|-----|-----|-----|-----|
| A1 | FR     | 0% | 100% | 100% | 100% | 0% | 0% | 0% |
|    | Best   | 2.752e+4/4.231e+2 | 8.186e-1 | 1.178e+1 | 2.340e+2 | 1.096/1.010e+2 | 1.002e+3/9.260e+5 | 5.635e+2/1.469e+5 |
|    | Median | 4.715e+4/3.292e+2 | 8.503e-1 | 1.492e+1 | 2.953e+2 | 1.100/1.010e+2 | 1.270e+3/1.794e+6 | 6.067e+2/1.469e+5 |
|    | Worst  | 7.453e+4/2.659e+2 | 9.381e-1 | 1.806e+1 | 3.707e+2 | 1.101/1.012e+2 | 1.495e+3/2.986e+6 | 6.332e+2/1.469e+5 |
|    | Mean   | 4.737e+4/3.356e+2 | 8.556e-1 | 1.505e+1 | 3.028e+2 | 1.100/1.012e+2 | 1.255e+3/1.970e+6 | 6.049e+2/1.469e+5 |
|    | STD    | 1.203e+4/1.609e+2 | 2.401e-2 | 1.081 | 3.825e+1 | 8.877e-4/1.638e-1 | 9.876e+1/5.325e+5 | 1.848e+1/4.700 |
| A2 | FR     | 0% | 100% | 100% | 100% | 0% | 0% | 0% |
|    | Best   | 3.274e+6/1.597e+3 | 1.115 | 2.121e+1 | 7.241e+2 | 1.099/1.010e+2 | 6.660e+3/1.314e+4 | 5.812e+2/1.469e+5 |
|    | Median | 5.701e+6/1.606e+3 | 1.138 | 2.435e+1 | 7.414e+2 | 1.100/1.010e+2 | 1.772e+4/2.593e+4 | 6.466e+2/1.468e+5 |
|    | Worst  | 9.646e+6/1.634e+3 | 1.155 | 2.435e+1 | 7.681e+2 | 1.100/1.010e+2 | 3.551e+4/3.513e+4 | 7.055e+2/1.469e+5 |
|    | Mean   | 5.804e+6/1.591e+3 | 1.136 | 2.322e+1 | 7.425e+2 | 1.100/1.010e+2 | 1.715e+4/2.335e+4 | 6.455e+2/1.468e+5 |
|    | STD    | 1.488e+6/1.053e+2 | 1.061e-2 | 1.508 | 1.052e+1 | 2.498e-4/0 | 7.199e+3/1.049e+4 | 2.566e+1/6.580 |
| A3 | FR     | 0% | 100% | 100% | 100% | 0% | 28% | 0% |
|    | Best   | 5.199e+4/2.433e+2 | 8.165e-1 | 1.806e+1 | 5.105e+2 | 1.088/1.010e+2 | 3.816e+1 | 6.252e+2/1.468e+5 |
|    | Median | 1.365e+5/3.414e+2 | 8.649e-1 | 2.121e+1 | 6.990e+2 | 1.100/1.010e+2 | 4.664e+1 | 6.729e+2/1.468e+5 |
|    | Worst  | 2.553e+5/3.084e+2 | 1.028 | 2.435e+1 | 7.603e+2 | 1.100/1.010e+2 | 5.548e+1 | 7.014e+2/1.469e+5 |
|    | Mean   | 1.323e+5/3.169e+2 | 8.909e-1 | 2.171e+1 | 6.778e+2 | 1.100/1.010e+2 | 4.648e+1 | 6.655e+2/1.468e+5 |
|    | STD    | 5.089e+4/5.265e+1 | 7.175e-2 | 1.705 | 5.647e+1 | 2.314e-3/0 | 5.098 | 2.312e+1/6.264 |
| A4 | FR     | 4% | 100% | 100% | 100% | 0% | 0% | 0% |
|    | Best   | 2.219e+4 | 0.849 | 1.492e+1 | 4.917e+2 | 1.093/- | 1.959e+4/- | 5.434e+2/- |
|    | Median | 4.115e+4 | 0.976 | 1.806e+1 | 5.482e+2 | 1.096/- | 2.994e+4/- | 5.691e+2/- |
|    | Worst  | 6.959e+4 | 1.024 | 1.806e+1 | 5.859e+2 | 1.098/- | 4.302e+4/- | 6.305e+2/- |
|    | Mean   | 5.036e+4 | 0.969 | 1.718e+1 | 5.442e+2 | 1.096/- | 3.688e+4/- | 5.840e+2/- |
|    | STD    | 1.182e+4 | 4.257e-2 | 1.440 | 2.860e+1 | 1.309e-3/- | 1.460e+4/- | 2.798e+1/- |